\documentclass[a4paper,11pt]{article}
\usepackage{amsmath,amsthm,amssymb}
\usepackage{booktabs}
\usepackage[latin1]{inputenc}
\usepackage{hyperref}
\usepackage{xcolor}

\def\seq#1#2#3{#1_{#2},\,\ldots,#1_{#3}}
\def\abs#1{\vert#1\vert}
\def\w{\widetilde}
\def\b{\overline}

\def\tt{{\underline{t}}}

\newcommand{\st} {\,|\,}
\newcommand{\ve}[3] {({#1}_{#2},\dotsc,{#1}_{#3})}
\newcommand{\range}[3]{#1 = #2, \dotsc, #3}
\newcommand{\Sii} {\Leftrightarrow}
\newcommand{\ord}{\textsf{ord}}
\newcommand{\sseq}[3]{#1^#2,\dotsc,#1^#3}
\newcommand{\e}{{\underline{1}}}

\def\0{\underline{0}}
\def\R{\mathbb R}
\def\P{\mathbb P}
\def\N{\mathbb N}
\def\Z{\mathbb Z}
\def\C{\mathbb C}
\def\OO{{\mathcal O}}
\def\PP{{\mathcal P}}
\newcommand{\Zr}[1]{{\mathbb Z}_{\ge 0}^{\,#1}}

\def\Fin{\textsf{Fin\,}}

\newcommand{\tD}{\widetilde{\Delta}}
\newcommand{\ul}{\underline}

\newtheorem{theorem}{Theorem}[section]
\newtheorem{corollary}[theorem]{Corollary}
\newtheorem{lemma}[theorem]{Lemma}
\newtheorem{proposition}[theorem]{Proposition}

\theoremstyle{definition}
\newtheorem{definition}[theorem]{Definition}
\newtheorem{remark}[theorem]{Remark}
\newtheorem{remarks}[theorem]{Remarks}
\newtheorem{example}[theorem]{Example}

\newtheorem{notations}[theorem]{Notations}
\numberwithin{equation}{section}

\newenvironment{listin}{%
\begin{list}{}{\itemsep=0pt \topsep=0pt}}%
{\end{list}}

\title{Completely Symmetric Good Semigroups and Poincar\'e series
}

\author{
F.~Delgado \thanks{Partially supported by the grant
PID2022-138906NB-C21 funded by
MICIU/AEI/10.13039/501100011033 and by ERDF/EU ``A way of making Europe''.
\newline
Address:
IMUVA (Instituto de Investigaci\'on en Matem\'aticas).
University of Valladolid. Valladolid, Spain.
E-mail: fdelgado\symbol{'100}uva.es.
}
}

\date{}

\begin{document}

\maketitle

\begin{abstract}
Good semigroups of $\Z^r_{\ge0}$ are the combinatorial
objects which generalize the value semigroups of reduced curve singularities with
$r$ branches. We develop, for an arbitrary good semigroup $S$, a purely
combinatorial theory of the Poincar\'e series $P_S(\tt)$.
The coefficients of $P_S$ are Euler characteristics of some natural simplicial
complexes attached to the
elements of $S$, so that the support of $P_S$ is contained in the set of maximal
elements of $S$. For symmetric good semigroups (the counterpart of Gorenstein
curves) these complexes satisfy an Alexander duality, which yields the functional
equation $P_S(\tt)=(-1)^r\tt^{\tau}P_S(\tt^{-1})$
and the exchange of absolute and
relative maximals. We then introduce completely symmetric semigroups, those whose
projections are all symmetric, compute their conductor in terms of the numbers
$\xi_{ij}$ of maximal elements of the two-dimensional projections, and prove a
projection formula for $P_S$, the analogue of the formula of Torres for the
Alexander polynomial, which holds as soon as the semigroup and its projection are
symmetric and fails for symmetric semigroups in general.
We also show that complete symmetry characterizes the
semigroups of plane curves among the good semigroups whose one-dimensional
projections are equal to the semigroup of a plane branch and whose two-dimensional
projections are diagonal (in particular, every completely symmetric semigroup
whose one-dimensional projections are all equal to $\Zr{}$ is the semigroup of a plane
curve). Finally, we show that a
completely symmetric semigroup is determined by its Poincar\'e series whenever the
numbers $\xi_{ij}$ are, which is the case for $r\le3$: a combinatorial counterpart of
the theorem of Yamamoto.
\end{abstract}

\section*{Introduction}\label{intro}

Let $C$ be a reduced curve singularity (an algebroid curve) with $r$ branches
$C_1,\dots,C_r$, let $\OO$ be its local ring and let
$\b\OO=\prod_{i=1}^r\C[[t_i]]$ be its normalization. The \emph{semigroup of values}
of $C$ is the set  $S=\{v(g)=(\ord_{t_1}g,\dots,\ord_{t_r}g) : g\in\OO \text{ a
non-zero-divisor}\}\subset\Z^r_{\ge0}$. It is a basic invariant of the singularity: for
plane curves it is equivalent to the topological (equisingularity) type of $C$ (see
\cite{Waldi}), and in general it encodes a large amount of information on the
ring $\OO$; for instance, by the theorem of Kunz \cite{Kunz} for $r=1$, and by
\cite{D-88} (see also \cite{CDK}) for any $r$, $\OO$ is Gorenstein if and only if $S$
is symmetric, a purely combinatorial condition. The arithmetical and combinatorial
properties of the semigroups of values were isolated in \cite{D-87,D-88} and in
Barucci, D'Anna and
Fr\"oberg \cite{Barucci-2000}, where the notion of \emph{good semigroup} was
introduced; since then good semigroups have been studied in their own right, see
for instance
\cite{Good,DGM,GMM-21, GMM-23,Korell-2019,Marcelo-2020}
and the references given there.
A key fact, proved in \cite{Barucci-2000}, is that not every good semigroup is the
semigroup of values of a curve. This is a permanent incentive to understand which
properties of the semigroups of values are of a purely combinatorial nature, that is,
consequences of the axioms of good semigroups, and which ones really depend on the
ring.
In this sense, several of the proofs in this paper (especially in
Sections~\ref{sec:Good} and \ref{Poincare S})
are not far from the proofs of similar results in the context of semigroups of
values. Our intention in including them is twofold: first, we aim to make the
article reasonably self-contained and thus accessible to a broad audience. Second,
and above all, we want to clearly identify those results that are purely
combinatorial and do not require algebraic or geometric conditions.

The second ingredient of the paper is the Poincar\'e series. In \cite{Duke} the
Poincar\'e series $P_C(\tt)$ of a curve singularity was defined from the multi-index
filtration of $\OO$ by the ideals $J(\alpha)=\{g\in\OO : v(g)\ge\alpha\}$, $\alpha\in\Z^r$,
and it was proved that, for a plane curve with $r\ge2$ branches, $P_C(\tt)$ coincides
with the Alexander polynomial of the link of $C$. Since the dimensions
$\dim\OO/J(\alpha)$ are determined by the semigroup $S$ (\cite{D-87}, see also the
distance function of \cite{Barucci-2000}), the series $P_C$ depends only on $S$, and
Yamamoto's theorem \cite{Yamamoto} (the Alexander polynomial determines the
topological type) can be read as: \emph{for plane curves, $P_S$ determines $S$}.
Two tools are essential in that circle of ideas: the projection formula of Torres
\cite{Torres}, which computes the Alexander polynomial of a sublink, and the
knowledge of the intersection multiplicities of the branches.

\medskip

For a good semigroup one can define its
Poincar\'e series $P_S(\tt)$ in a purely combinatorial way (see
Section~\ref{Poincare S}). However, unlike in the case of plane
curves, in general it contains less information than the semigroup (see
Example~\ref{ex:spacecusps} and Remark~\ref{rem:PS-sym-not-S}). The natural question
we address in this article is therefore under which conditions these two objects are
equivalent. Since the description of the semigroup (see~\cite{D-87}) relies heavily on
inductive arguments based on its projections, a natural approach is to determine the
conditions on the semigroup under which Torres' projection formula holds.

This idea leads to a strengthened notion of symmetry of a good
semigroup, which we call \emph{complete symmetry}. Consequently, a central part of
the paper is an in-depth study of the notion of symmetry. For semigroups of values,
complete symmetry means that all the curves formed by a subset of the branches of the
original curve are Gorenstein. To the best of our knowledge, the meaning of this
condition in geometric or algebraic terms has not been explored yet.

Roughly speaking, a good semigroup is completely symmetric when
there is an combinatorial notion of intersection multiplicity between
pairs of branches ---here, between pairs of numerical semigroups obtained as
projections of the given one--- which behaves, from a combinatorial point of view,
like the intersection multiplicity of two plane branches
(Theorem~\ref{thm:supersymmetry} and Proposition~\ref{prop:Apery-fibres}). It is
worth pointing out that this notion is still far from characterizing the semigroups
of plane curves: already in the irreducible case, the (well-known) characterization
is expressed in terms of the minimal system of generators, and nothing resembling the
intersection multiplicity of two branches plays a part.

\medskip

The main results and the organization of the paper are
as follows. Section~\ref{sec:Good} recalls the definition and
the
basic properties of good semigroups (conductor, fibres of the projections, maximal
elements and the generation theorem of \cite{D-87,D-88}).
In Section~\ref{Poincare S}
we define $P_S(\tt)$ for any good semigroup through two combinatorial
functions, the colength $\ell(\alpha)$ and the dimension $\lambda(\alpha)$ (these
functions are the combinatorial version of the colength of the valuation ideals
$J(\alpha)$ and the dimension of $J(\alpha)/J(\alpha+\e)$)
and we prove that the coefficient of $\tt^\alpha$ in
$P_S$ is the Euler characteristic of a simplicial complex $\nabla_{\alpha,i}$ built
from the elements of $S$ which are involved in $\lambda(\alpha)$
(Theorem~\ref{coef-PS}). This is
the combinatorial counterpart of the description of the coefficients of $P_C$ as the
Euler characteristic of the fibres of the extended semigroup and can be seen also as
an integral with respect
to the Euler characteristic (see~\cite{Duke, IJM-03}). As a consequence, the support
of
$P_S$ is contained in the set of maximal elements of $S$
(Proposition~\ref{prop:coef}),
but it may be strictly smaller (Example~\ref{ex:p-zero}).

Section~\ref{sec:Symmetry} is devoted to
symmetric semigroups.
We show that the symmetry of a good semigroup is
equivalent to an Alexander duality between the complexes
$\nabla_{\alpha,i}$ and $\nabla_{\tau-\alpha,i}$, $\tau$ being the Frobenius vector
(Proposition~\ref{teo-dualidad}). This gives at once the functional equation
$P_S(\tt)=(-1)^r\tt^{\tau}P_S(\tt^{-1})$ and the fact that $\alpha\mapsto\tau-\alpha$
exchanges absolute and relative maximals (Theorem~\ref{thm:duality-maximals}), so
that a symmetric semigroup is generated, in the sense of \cite{D-87}, by its
irreducible absolute maximals (Theorem~\ref{AbsoluteMaximals}).

In Section~\ref{Supersymmetric} we consider \emph{completely
symmetric}
semigroups, i.e.\ good semigroups all of whose projections are symmetric; the
semigroups of plane curves are the main examples. We prove that the conductor
$\delta$ of such a semigroup is given by $\delta_i=c_i+\sum_{j\neq i}\xi_{ij}$, where
$c_i$ is the conductor of the projection $S_i$ and $\xi_{ij}$ is the number of
maximal elements of the projection $S_{\{i,j\}}$, which plays the role of the
intersection multiplicity of two branches (Theorem~\ref{thm:supersymmetry}),
and that the
elements of a projection $S_J$ with finite fibre form the Ap\'ery set of $S_J$ with
respect to a suitable element (Proposition~\ref{prop:Apery-fibres}). The relation with
the semigroups of plane curves is studied in Subsection~\ref{sec:CSplane}: the class
of completely symmetric semigroups is strictly larger (Example~\ref{ex:spacecusps}),
but Theorem~\ref{thm:diagonal-plane} shows that a completely symmetric semigroup whose
one-dimensional projections are all equal to the semigroup $M$ of a plane branch and
whose two-dimensional projections are diagonal (Definition~\ref{def:diagonal}) is the
semigroup of a plane curve; this is always the case when all the projections $S_i$
are equal to $\Zr{}$ (Corollary~\ref{cor:smooth-plane}).
We then
prove a projection formula for
the Poincar\'e series (Theorem~\ref{projection} and Corollary~\ref{cor:projection1}),
$P_S(\tt)|_{t_i=1}=(1-\prod_{j\neq i}t_j^{\xi_{ij}})P_{S^i}$, which is the analogue of the
formula of Torres. More generally, a formula of this type (with exponents given
by the conductors) holds as soon as $S$ and the projection $S^i$ are
symmetric; a projection formula for Gorenstein curves was stated in \cite{Moyano},
but we show by an example that the Gorenstein property of $C$ alone is not
sufficient, and that for $r=2$ the formula holds exactly when the remaining branch is
Gorenstein (Proposition~\ref{thm:two}).

Finally (Section~\ref{sec:recovering}) we prove that a completely symmetric
semigroup is determined by its Poincar\'e series as soon as the numbers $\xi_{ij}$
are determined by $P_S$ (Theorem~\ref{thm:recons}), and that this is the case for
$r\le3$ (Theorem~\ref{thm:xi-r3}, Corollary~\ref{cor:PS-semigroup}). This is the
combinatorial counterpart of Yamamoto's theorem; whether it holds for every $r$ is
an open problem, which we discuss at the end of the paper.

\section{Good semigroups}\label{sec:Good}

Throughout  the paper $r$ is a positive integer and $I=\{1,\ldots, r\}$; unless
otherwise specified, we assume that $r>1$. We consider
$\Z^r$ (partially) ordered by the product order, that is,
$$
\alpha= \ve{\alpha}1r\le \beta=\ve{\beta}1r \Sii
\alpha_i\le \beta_i \mbox{ for all }
\range i1r\; .
$$
Let $\seq {\e}1r$ denote the canonical
basis of $\Z^r$, $\e = \sum \e_i = (1,\ldots,1)$ and, if $J\subset I=\{1,\ldots,r\}$,
$\e_J=\sum_{i\in J}\e_i$. For $a\in \Z$ we denote
$\ul{a}= a\cdot \e = (a,\ldots,a)\in \Z^r$.

For $J\subset I=\{1,\ldots, r\}$, let
$pr_J:\Z^r\to \Z^{\# J}$ be the projection onto the coordinates indexed by $J$.
Given $\alpha\in \Z^r$, $\alpha_J$ will denote $pr_J(\alpha)$.
For the sake of simplicity we will also use $pr^J$, $\alpha^J$ to denote the
projection onto the coordinates indexed by $I\setminus J$.

\begin{definition}
A subsemigroup $S\subset (\Zr r, +)$ will be called a {\it good semigroup} when
it satisfies the following properties:
\end{definition}
\begin{listin}
\item[P1)]
There exists $\delta\in S$ such that $\delta+\Zr r\subset S$.
\item[P2)]
If $\alpha, \beta\in S$ then
$
\alpha\wedge \beta = \inf (\alpha, \beta) := (\min(\alpha_1,\beta_1), \ldots,
\min(\alpha_r,\beta_r)) \in S
$
\item[P3)]
If $\alpha, \beta\in S$  and  $\alpha_i=\beta_i$ for some index
$i\in I$, then there exists
$\gamma \in S$ such that
$\gamma_j \ge \min \{\alpha_j, \beta_j\}$ for all $j=1,\ldots,r$;
$\gamma_j=\min \{\alpha_j, \beta_j\}$ provided $\alpha_j\neq \beta_j$ and, moreover,
$
\gamma_i > \alpha_i = \beta_i\; .
$
\item[P4)]
$\0=(0,\ldots,0)\in S$. If  $\alpha = \ve{\alpha}1r \in S$ with
$\alpha_i=0$ for some $i\in I$
then $\alpha=\0$.
\end{listin}

The definition of a good semigroup usually includes only properties P1 to P3 (see
\cite{Barucci-2000} and
\cite{Good}); when
condition P4 is also fulfilled, the semigroup is said to be a {\it local good
semigroup}. In our
case, all the semigroups we will deal with are local, therefore we incorporate
condition P4 into the definition.
In the sequel all the semigroups $S\subset \Zr r$ considered will be (local) good
semigroups.
Note that if $r = 1$,  then saying that $S$ is a good semigroup simply means that $S$ is
a numerical semigroup.

\begin{remark}\label{rem:conductor}
The element
$\gamma$ in property P3) is, in general, not unique. However, sometimes we will
denote by
$\alpha
\vee_{i} \beta \in S$ any  element $\gamma$
with these properties.

Let $J\subset I$ be a subset of $I$ and let
$S_J = pr_{J}(S)\subset \Zr {\# J}$ be the projection of $S$ onto $J$.

It is obvious
that  $S_J$ is also a good semigroup (for the property P3, the element $\gamma$
given by P3 for two lifts of $\alpha_J,\beta_J$ with $\alpha_i=\beta_i$ ($i\in J$)
projects onto an element with the required properties). In particular, for $i\in I$,
$S_i :=S_{\{i\}}\subset \Zr {}$ is a numerical semigroup and $S\subset
S_1\times\cdots
\times S_r$. We will also use $S^J$ (and $S^i$) instead of $S_{I\setminus J}$
(instead of $S_{I\setminus \{i\}}$).

If $\delta,\delta'$
are two elements with $\delta+\Zr r\subset S$ and $\delta'+\Zr r\subset S$, then every
$\gamma\ge\delta\wedge\delta'$ can be written as
$\gamma=(\gamma\vee\delta)\wedge(\gamma\vee\delta')$ with
$\gamma\vee\delta\in\delta+\Zr r$ and $\gamma\vee\delta'\in\delta'+\Zr r$ (here
$\vee$ is the coordinatewise maximum), so $\delta\wedge\delta'+\Zr r\subset S$.
Thus, one can see that there exists a minimum element $\delta\in S$
such that $\delta+\Zr r\subset S$. This element is called the {\bf conductor} of the
semigroup $S$ and we will denote it by $\delta(S)$, or simply by $\delta$ when the
semigroup is clear from the context. The element $\tau :=\delta -\e $ is
usually called the
{\bf
Frobenius vector} of $S$ (it is the analogue of the Frobenius number of a numerical
semigroup). As above, we will use $\tau(S)$ when we want to indicate the semigroup
explicitly. In the sequel we denote by
$c_i$ the conductor of the numerical semigroup $S_i$, $(1\le i \le r)$, and by
$\delta = \delta(S)=(c_1+\mu^1, \ldots, c_r+\mu^r)\in S$ the conductor of $S$. Note
that $\mu^i\ge 0$; however,  it may happen that  $\mu^i\notin S_i$.
\end{remark}

\begin{proposition}\label{prop:small}
The semigroup $S$ is never finitely generated if $r>1$. However,
it is finitely determined. More precisely,
the semigroup is uniquely determined by the
elements $\alpha\in S$ such that $0\le \alpha\le \delta$.
\end{proposition}

\begin{proof}
Let us assume that $S$ has a finite number of generators
$\beta^1, \ldots, \beta^k\in S\setminus \{\0\}$, i.e.
$$
S=\langle\beta^1,\ldots,\beta^k\rangle_{\Z_{\ge 0}} =
\{\lambda_1 \beta^1 + \cdots + \lambda_k \beta^k \st \lambda=(\seq{\lambda}1k)\in \Zr
k\}\;
.
$$
As
$\beta^j_i>0$ for all $j$ and for all $i$,
$\langle\beta^1,\ldots,\beta^k\rangle_{\R_{\ge 0}}$
is a (closed) cone in
$\R^r_{\ge 0}$ which only meets the coordinate hyperplanes at the origin.
As
$S=\langle\beta^1,\ldots,\beta^k\rangle_{\Z_{\ge 0}} \subset
\langle\beta^1,\ldots,\beta^k\rangle_{\R_{\ge 0}}$,
this contradicts the existence of the conductor: if $\delta$ is the conductor, then
$\delta+n \e_1\in S$ for all $n\ge 0$; now, $(\delta+n\e_1)/n\to \e_1$ when $n\to \infty$,
so $\delta+n\e_1$ in the cone for all $n$ would give
$\e_1$ in the cone, which is impossible.

\medskip

Let $\alpha\in S$ and assume that there is $i\in I$ with
$\alpha_i\ge \delta_i = c_i+\mu^i$.
Let $\beta\ge \delta$ be such that
$\beta_i=\alpha_i$ and  $\beta_j> \alpha_j$ for all $j\neq i$.
By P3, there exists an element
$\gamma = \alpha\vee_i \beta\in S$ such that
$\gamma_i>\alpha_i$ and  $\gamma_j=\alpha_j$ for all
$j \neq i$.
Iterating this argument one gets elements $\gamma^{(n)}\in S$ with $\gamma^{(n)}_i\to\infty$ and
$\gamma^{(n)}_j=\alpha_j$ for $j\ne i$. Now for a given $m\ge\delta_i$, the
element with $i$-th coordinate $m$ and $j$-th coordinate $\alpha_j$ is
$\gamma^{(n)}\wedge\varepsilon$ for $n$ large and $\varepsilon\ge\delta$ with
$\varepsilon_i=m$, $\varepsilon_j>\alpha_j$ ($j\ne i$). That is,
all the elements $\beta$ with $\beta_i\ge \delta_i$ and
$\beta_j=\alpha_j$ for all $j\neq i$ belong to the semigroup $S$.

Repeating the same argument for all the indices
of
$J:=\{j\in I\st \alpha_j\ge \delta_j\}$ we find that any
$\beta\in \Zr r$ with $\beta_j\ge \delta_j$ for  $j\in J$ and
$\beta_i=\alpha_i$ for  $i\notin J$ belongs to the semigroup.
As a consequence, $S$ is determined by the so-called \emph{small elements}, that is, by
the  finite set
$$
(\0, \delta]\cap S =\{ \alpha\in S \st \0 < \alpha\le
\delta \}\; .
$$
\end{proof}

\begin{notations}\label{not1}
Let $J\subset I$ and
$\alpha_J\in \Zr {\# J}$. Given
$\gamma\in \Z^{r}$  we will
denote by $\alpha_J*\gamma \in \Z^r$ the element given by
$pr_J(\alpha_J * \gamma) = \alpha_J$ and
$pr_{i}(\alpha_J*\gamma)=\gamma_i$ for all $i\notin J$.
Also, we will use $F(\alpha_J)$ to denote the fibre of $\alpha_J$ in $S$, that is,
$$
F(\alpha_J)=pr_J^{-1}(\alpha_J)\cap S = \{ \beta\in
S\st \beta_J=\alpha_J\}\; .
$$
One has that
$\alpha_J\in S_J$ if and only if $F(\alpha_J)\neq\emptyset$.

From the proof of the above proposition one easily deduces the equivalence of the
following conditions for an element $\alpha_J\in S_J$:
\begin{listin}
\item[1)] $\alpha_J*\delta\in S$.
\item[2)] There exists $\beta\in F(\alpha_J)$ such that
$\beta^J\ge \delta^J$.
\item[3)] $\{\beta\in \Zr r\st \beta_J=\alpha_J,
\beta^J\ge \delta^J\}\subset S$
\item[4)] There exists a sequence
$\{\alpha^n\}_{n=1}^{\infty}$ of elements of $S$ such that
$\alpha^n_{J}=\alpha_J$ for all
$n$ and  $\alpha^{n+1}_i>\alpha^n_i$ for all $n\ge 1$ and all $i\notin J$.
\end{listin}

\medskip

When the above conditions are satisfied for an element $\alpha_J\in S_J$ we will say
that
$\alpha_J$ has {\it infinite fibre\/} and we will indicate this by writing
$F(\alpha_J)=\infty$; otherwise, we will say that $\alpha_J$
has {\it finite fibre\/} ($F(\alpha_J)<\infty$).
We will denote by $\Fin(S_J)\subset S_J$ the set of elements of $S_J$ with finite fibre.
\end{notations}

\begin{remark}\label{complete}
Let $\Zr {}$ be equipped with the discrete topology and let
$\Zr {} \cup \{\infty\}$ be its one point compactification, with
$a+\infty =\infty$ for every $a\in \Z_{\ge 0}$. The completion $\b{S}$ of $S$ is
defined as the topological closure of $S$ in the product space
$(\Zr {} \cup \{\infty\})^r$. Obviously, the data $\b{S}$ and $S$ are equivalent to each
other.
Given $\alpha_J\in S_J$, let us denote by $\alpha_J*\infty\in (\Z\cup \{\infty\})^r$
the element such that  $pr_J(\alpha_J* \infty)=\alpha_J$ and
$(\alpha_J*\infty)_i = \infty$ for all $i\notin J$.
One also has that
$\alpha_J\in S_J$ has infinite fibre ($F(\alpha_J)=\infty$)
if and only if:
\medskip

\begin{listin}
\item[5)] $\alpha_J * \infty\in \b S$.
\end{listin}
\end{remark}

\begin{remark}\label{rem:fibre-set}
Note that the set $F(\alpha_J)$ may have an infinite number of elements while
the fibre of $\alpha_J$ is finite in the above sense.
As an example, let
$$
S =\{(0,0,0), (1,1,k), (1,k,1), (k,1,1) : k\ge 1\}\cup \{(2,2,2)+\Z_{\ge 0}^{3}\}
$$
be the semigroup of three smooth and pairwise transversal plane branches.
Let $J=\{1\}$ and  $1_1\in S_1$; one has that
$F(1_1) = \{(1,1,k), (1,k,1): k\ge 1\}$ is an infinite set.
However, $F(1_1)< \infty$ because $(1,k,m)\notin S$ for $k,m\ge 2$.
\end{remark}

\begin{notations}\label{not2}
Let $F\subset I :=\{1,\ldots,r\}$ and $\alpha\in \Z^r$. The following notations will be
widely used throughout the paper.
\begin{align*}
\Delta_F(\alpha) &= \{ \beta\in S | \beta_i=\alpha_i \ \forall i\in F\text{ and }
\beta_j > \alpha_j \ \forall j\notin F\}\; ,\\
\Delta^{F}(\alpha)
 &= \{ \beta\in S | \beta_i>\alpha_i \ \forall i\in F\text{ and }
\beta_j = \alpha_j \ \forall j\notin F\}
= \Delta_{I\setminus F}(\alpha)  \; ,
 \\
\w{\Delta}_F(\alpha) &= \{ \beta\in S | \beta_i=\alpha_i \ \forall i\in F\text{ and }
\beta_j\geq \alpha_j \ \forall j\notin F\} = {\bigcup}_{F\subset G}
\Delta_G(\alpha)
\\
\Delta(\alpha) &= {\bigcup}_{i=1}^r \Delta_{\{i\}}(\alpha) \; .
\end{align*}
When $F=\{i\}$ we simply write
$\Delta_{i}(\alpha)$ or
$\Delta^{i}(\alpha)$ instead of
$\Delta_{\{i\}}(\alpha)$ or
$\Delta^{\{i\}}(\alpha)$.
\end{notations}

\begin{definition}\label{def:maximal}
Let $S\subset \Zr r$, $r\ge 2$,  be a good semigroup and $\alpha\in S$. We will say
that $\alpha$ is a {\it maximal element} if $\Delta(\alpha)=\emptyset$.

If, moreover,
$\Delta_F(\alpha) = \emptyset$ for all $F\subset I$, $F\neq \emptyset, I$, we call
$\alpha$ an {\it absolute maximal}.

In a similar way, we will say that $\alpha$ is a {\it relative maximal}
if $\alpha$ is maximal and $\Delta_F(\alpha)\neq
\emptyset$ for all
$F\subset I$ with $\# F \ge 2$.
\end{definition}

\begin{remark}\label{rem:maximal}
If $r=1$ and, in the above definition,
we replace the condition $\alpha\in S$ by $\alpha\in
\Zr{}$,
one has that $\alpha$ is maximal if $\alpha$ is a gap (i.e.\ $\alpha\in
\Zr{}\setminus S$). However, it is more convenient to
use this convention for relative maximals only. So, if $r=1$ we will say that
$\alpha\in \Zr{}$ is a relative maximal if $\alpha$ is a gap. On the other hand,
we will say that $\alpha\in \Zr{}$ is an absolute maximal if $\alpha<c$ ($c$ being the
conductor of $S$) and $\alpha\in S$.
If $r=2$ the notions of maximal, absolute maximal and relative maximal coincide.

\medskip

If $\alpha$ is a maximal element, then $\alpha\in [\0,\tau]$:
let $\alpha\in S$ be such that $\alpha_i\ge \delta_i$ for some $i\in I$, and
let $\beta$ be such that $\beta_i=\alpha_i\ge \delta_i$ and $\beta_j\ge \delta_j$,
$\beta_j>\alpha_j$ for $j\neq i$. Then $\beta\ge \delta$, so $\beta\in S$, and obviously
$\beta\in \Delta_i(\alpha)$.

So, in particular, there is a finite number of
maximals.
Note that $\0\in S$ is {\bf always} an absolute maximal. For the Frobenius vector
$\tau$ one
can easily see that $\Delta(\tau)=\emptyset$; however, it is not true in general that
$\tau\in S$: the semigroup $\{(0,0)\}\cup (\delta+\Zr{2})$ is a trivial example.

Let $\gamma\in S$. If
$\sigma\in\Delta_i(\tau-\gamma)$, then $\sigma+\gamma\in\Delta_i(\tau)=\emptyset$.
Hence one has the following simple but relevant fact:
$$
\gamma\in S\Rightarrow\Delta(\tau-\gamma)=\emptyset\; .
$$
\end{remark}

The maximal elements play a central role in the inductive description of
good
semigroups (see \cite{D-87}, \cite{D-88}). One has the following

\begin{theorem}\label{generation}
Let $\alpha\in \Z_{> 0}^r$ and assume that $\w{\Delta}_i(\alpha)=\emptyset$ for some
index $i\in I$ (in particular, this is the case if
$\alpha\notin S$). Then, there exist
$J\subset I$ with $i\in J$  and a relative maximal $\beta_J$ of $S_J$ such that
$\alpha_i = \beta_i$ and $\alpha_j > (\beta_J*0)_j$ for $j\neq i$.

In particular, if $\alpha_K\in S_K$ for all $K\subset I$, $K\neq I$, then
$J=I$ and $\beta$ is a relative maximal of $S$.
\end{theorem}

\begin{proof}
The result is just Theorem (1.5) of
\cite{D-88}.

For the second statement, if $J\neq I$ and $\#J\ge 2$ then
$\alpha_J\in\Delta^{S_J}_i(\beta_J)$ (the latter being the set $\Delta_i$ computed in the semigroup
$S_J$), which contradicts the maximality of $\beta_J$. If $J=\{i\}$, then $\beta_i$ is a gap of
$S_i$ while $\alpha_i\in S_i$, against the hypothesis.

Notice also that the second statement is the Generation Theorem (1.5) of
\cite{D-87}, where an independent proof can be found.
\end{proof}

The next result will be used later in Section~\ref{sec:Symmetry}.

\begin{lemma}\label{lem:new}
Let $\alpha\in\Z^r$ be such that $\w\Delta_i(\alpha)=\emptyset$ for some $i\in I$.
Then there exists
$\beta\in\Z^r$ with $\beta_i=\alpha_i$, $\beta_j<\alpha_j$ for $j\ne i$ and
$\Delta(\beta)=\emptyset$.
\end{lemma}

\begin{proof}

Throughout the proof we use the following two elementary facts, which follow from
$S\subset\Z^r_{\ge0}$ and P4: if $\beta\in\Z^r$ and $\beta_k<0$ for some $k$, then
$\Delta_k(\beta)=\emptyset$; and if $\beta_k=0$, then $\Delta_k(\beta)\subset\{\0\}$.

\emph{Case 1: $\alpha_i<0$.} Let $\beta_i=\alpha_i$ and $\beta_j=\min\{\alpha_j,0\}-1$
for $j\neq i$. All the coordinates of $\beta$ are negative, so $\Delta(\beta)=\emptyset$,
and clearly $\beta_i=\alpha_i$ and $\beta_j<\alpha_j$ for $j\neq i$.

\emph{Case 2: $\alpha_i\ge0$.} Let $K=\{i\}\cup\{j\in I:\alpha_j>0\}$ and
$Z=I\setminus K$, so that $\alpha_j\le0$ for $j\in Z$. For $\sigma\in S$ the conditions
$\sigma_j\ge\alpha_j$, $j\in Z$, are automatically satisfied, hence
\[
\w\Delta_i(\alpha)=\{\sigma\in S : \sigma_i=\alpha_i,\ \sigma_j\ge\alpha_j \text{ for } j\in K\}
= (pr_K |_{S})^{-1}\big(\w\Delta_i^{S_K}(\alpha_K)\big)\; ,\tag{$*$}
\]
where $\w\Delta_i^{S_K}$ denotes the set $\w\Delta_i$ computed in the good
semigroup $S_K$. In particular $\w\Delta_i^{S_K}(\alpha_K)=\emptyset$. We will
construct $\beta$ with $\beta_j=\min\{\alpha_j,0\}-1<\alpha_j$ for $j\in Z$ (so that
$\Delta_k(\beta)=\emptyset$ for $k\in Z$, and the conditions $\sigma_j>\beta_j$,
$j\in Z$, are empty for $\sigma\in S$), which reduces the problem to the choice of
$\beta_K$: as in $(*)$, for $k\in K$ one has
$\Delta_k(\beta)=pr_K^{-1}\big(\Delta_k^{S_K}(\beta_K)\big)$, so it suffices to find
$\beta_K\in\Z^{K}$ with $\Delta^{S_K}(\beta_K)=\emptyset$, $\beta_i=\alpha_i$ and
$\beta_j<\alpha_j$ for $j\in K\setminus\{i\}$.

\emph{Case 2a: $K=\{i\}$.} Then $\w\Delta_i^{S_i}(\alpha_i)=\emptyset$ means
$\alpha_i\notin S_i$, and $\beta_K=\alpha_i$ works: $\Delta_i^{S_i}(\alpha_i)=\emptyset$
because $\alpha_i\notin S_i$. (This is the case $r=1$ of
Theorem~\ref{generation}, with the convention of Remark~\ref{rem:maximal}.)

\emph{Case 2b: $K\neq\{i\}$ and $\alpha_i=0$.} Take $\beta_j=0$ for all $j\in K$. For
$k\in K$, an element $\sigma\in\Delta_k^{S_K}(\beta_K)$ has $\sigma_k=0$, hence
$\sigma=\0$ by P4, but then $\sigma_j>\beta_j=0$ fails for $j\in K\setminus\{k\}\neq\emptyset$.
Thus $\Delta^{S_K}(\beta_K)=\emptyset$, and $\beta_i=0=\alpha_i$, $\beta_j=0<\alpha_j$ for
$j\in K\setminus\{i\}$.

\emph{Case 2c: $\alpha_i>0$.} Then $\alpha_K\in\Z^{K}_{>0}$ and
$\w\Delta_i^{S_K}(\alpha_K)=\emptyset$, so Theorem~\ref{generation} applied to the
good semigroup $S_K$ provides $J\subset K$ with $i\in J$ and a relative maximal
$\beta_J$ of $S_J=(S_K)_J$ with $\beta_i=\alpha_i$ and $\alpha_j>(\beta_J*\0)_j$ for
$j\in K\setminus\{i\}$. Put $\beta_K=\beta_J*\0\in\Z^{K}$. Let $k\in K$ and suppose
$\sigma\in\Delta_k^{S_K}(\beta_K)$. If $k\notin J$, then $\sigma_k=0$, so $\sigma=\0$ by
P4, against $\sigma_i>\beta_i=\alpha_i>0$. If $k\in J$ and $\#J\ge2$, then
$pr_J(\sigma)\in\Delta_k^{S_J}(\beta_J)$, against the maximality of $\beta_J$ in
$S_J$. If $J=\{k\}=\{i\}$, then $\sigma_i=\beta_i$ is a gap of $S_i$, which is
impossible for $\sigma\in S_K$. Hence $\Delta^{S_K}(\beta_K)=\emptyset$, and $\beta_K$
has the required properties.
\end{proof}

\section{The Poincar\'e series.}\label{Poincare S}

In this section we attach to a good semigroup $S$ a series
$P_S(\tt)\in\Z[[t_1,\dots,t_r]]$ which, when $S$ is the value semigroup of a reduced
curve singularity $C$, coincides with the Poincar\'e series of \cite{Duke}.
The construction uses only two combinatorial functions on $\Z^r$, the
\emph{colength} $\ell(\alpha)$ and the \emph{dimension} $\lambda(\alpha)$, which
play the role of $\dim\OO_C/J(\alpha)$ and $\dim J(\alpha)/J(\alpha+\e)$,
$J(\alpha)=\{g\in\OO_C: v(g)\ge\alpha\}$. These definitions closely follow the
construction of the Poincar\'e series of a curve given in \cite{Duke}.
The point here is
that all the ingredients can be written in terms of the semigroup alone.
This point of view is not new: a combinatorial Poincar\'e series for good
semigroups, and more generally for good semigroup ideals, was introduced by Tozzo
\cite{Tozzo}.
What is really new here is the description of the coefficients of $P_S$ as Euler
characteristics of simplicial complexes attached to the elements of $S$
(Theorem~\ref{coef-PS}).

\subsection{Colength and dimension}

\begin{definition}
Given $\alpha\in \Z^r$ and $i\in I$, let $\ell(\alpha,i)\in\{0,1\}$ be defined by
$\ell(\alpha,i)=0$ if and only if $\w{\Delta_i}(\alpha)=\emptyset$ (and
$\ell(\alpha,i)=1$ otherwise).
\end{definition}

\begin{lemma}\label{lem:commute}
Let $\alpha\in \Z^r$ and $i,j\in I=\{1,\ldots, r\}$. One has
$$
\ell(\alpha,i)+\ell(\alpha+\e_i,j) =
\ell(\alpha,j)+\ell(\alpha+\e_j,i)\; .
$$
\end{lemma}

\begin{proof}
One has the following disjoint decompositions:
$$
\w{\Delta}_i(\alpha) = \w{\Delta}_i (\alpha+\e_j) \cup \w\Delta_{i,j}(\alpha)
\text{ and }
\w{\Delta}_j(\alpha) = \w{\Delta}_j (\alpha+\e_i) \cup \w\Delta_{i,j}(\alpha) \; .
$$
As a consequence, if $\ell(\alpha,i)=0$ then $\ell(\alpha+\e_j,i)=0$ and, since
$\w\Delta_{i,j}(\alpha)=\emptyset$,
$\w{\Delta}_j(\alpha)=\emptyset$
if and only if $\w{\Delta}_j(\alpha+\e_i)=\emptyset$, i.e.\
$\ell(\alpha,j)=0$
if and only if
$\ell(\alpha+\e_i,j)=0$. Thus the equality holds in the case
$\ell(\alpha,i)=0$ (or
$\ell(\alpha,j)=0$).

Let us assume that $\ell(\alpha,i)=1$ and $\ell(\alpha+\e_i,j)=1$. Notice that in this
case one has $\ell(\alpha,j)=1$. Let
$\beta\in \w{\Delta}_i(\alpha)$ and
$\gamma \in \w{\Delta}_j(\alpha+\e_i)$. If
$\beta_j > \alpha_j$ then $\beta\in \w{\Delta}_i(\alpha+\e_j)$, so
$\ell(\alpha+\e_j,i)=1$ and the equality follows.
Let us suppose that $\beta_j = \alpha_j$; since $\gamma_j = \alpha_j$ one has that
$\beta_j = \gamma_j$ and $\beta_i<\gamma_i$, so the element
$\rho = \beta\vee_j \gamma$ satisfies
$\rho_j> \beta_j =\gamma_j = \alpha_j$, $\rho_i = \beta_i =\alpha_i$ and
$\rho_k \ge \min \{\beta_k,\gamma_k\}\ge \alpha_k$ for $k\neq i, j$.
This means that $\rho\in \w{\Delta}_i(\alpha+\e_j)\neq \emptyset$, so
$\ell(\alpha+\e_j,i)=1$.
The case $\ell(\alpha,i)=\ell(\alpha,j)=1$, $\ell(\alpha+\e_i,j)=0$
follows from the previous one by exchanging the roles
of $i$ and $j$ (if $\ell(\alpha+\e_j,i)$ were $1$, the case already treated would give
$\ell(\alpha+\e_i,j)=1$). So,
the proof is finished.
\end{proof}

Let $\alpha\in \Z_{\ge 0}^r$, and let $|\alpha| = \sum_{i=1}^r \alpha_i$. Let
$\ell(\alpha)$ be the nonnegative integer defined recursively on
$\abs{\alpha}$
by
$\ell(\0)=0$ and, if $\abs{\alpha}>0$,
$\ell(\alpha) = \ell(\alpha-\e_i) + \ell(\alpha-\e_i,i)$ for any $i$ such that
$\alpha_i>0$.  Notice that
\begin{align*}
\ell(\alpha) & = \ell(\alpha-\e_i) + \ell(\alpha-\e_i,i) =
\ell(\alpha-\e_i-\e_j) + \ell(\alpha-\e_i-\e_j,j)+\ell(\alpha-\e_i,i) \\
&= \ell(\alpha-\e_i-\e_j) + \ell(\alpha-\e_i-\e_j,i)+\ell(\alpha-\e_j,j) \\
&= \ell(\alpha-\e_j) + \ell(\alpha-\e_j,j)
\end{align*}
(by Lemma~\ref{lem:commute}), and so $\ell(\alpha)$ does not depend on the index $i$ chosen.

\begin{definition}
For $\alpha\in \Z_{\ge 0}^r$ the {\it colength\/} of $\alpha$ is the integer
$\ell(\alpha)$ and the {\it dimension\/} of $\alpha$ is defined as
$$
\lambda(\alpha):=\ell(\alpha+\e)-\ell(\alpha)\; .
$$
\end{definition}

\begin{remark}\label{rem:curves}
Notice that $\ell(\alpha)$ coincides with the {\it distance} $d(\0,\alpha)$ defined
in \cite{Barucci-2000}. The names \emph{colength} and \emph{dimension} come from the
case of curves. Let $S$ be the semigroup of values of a reduced curve singularity
$C$ with $r$ branches, let $\OO=\OO_C$ be its local ring and let
$v=(v_1,\dots,v_r)$ be the valuation, $v_i(g)$ being the order of $g$ along the
$i$-th branch. For $\alpha\in\Z^r$ consider the ideal
$$
J(\alpha)=\{g\in\OO : v(g)\ge\alpha\}\; ,
$$
where $v(g)\ge\alpha$ is understood coordinatewise and with the convention
$v_i(g)=\infty$ if $g$ vanishes on the $i$-th branch. Then
$J(\alpha)\supset J(\alpha+\e_i)\supset J(\alpha+\e)$, the quotient $J(\alpha)/J(\alpha+\e_i)$
has dimension at most $1$ (it embeds in $\C$ through the coefficient of
$t_i^{\alpha_i}$), and it is nonzero if and only if there exists $g\in\OO$ with
$v(g)\ge\alpha$ and $v_i(g)=\alpha_i$, that is, if and only if
$\w\Delta_i(\alpha)\neq\emptyset$. Hence $\ell(\alpha,i)=\dim J(\alpha)/J(\alpha+\e_i)$,
and by induction on $\abs\alpha$,
$$
\ell(\alpha)=\dim_{\C}\OO/J(\alpha)\qquad\text{and}\qquad
\lambda(\alpha)=\dim_{\C}J(\alpha)/J(\alpha+\e)
$$
for $\alpha\in\Z^r_{\ge0}$ (see \cite{D-87}, where $\dim\OO/J(\alpha)$ is computed
from $S$ in this way). In particular, the integers $\dim\OO/J(\alpha)$ are
invariants of the semigroup. For plane curves, the function
$\alpha\mapsto\dim\OO/J(\alpha)$ is the function $h$ of Gorsky and N\'emethi
\cite{GN}, from which they compute the lattice homology of the link of $C$.

\end{remark}

\begin{remarks}\label{colength}
{\bf (1)} The above definitions can be extended to any $\alpha\in \Z^r$ in the
following way. Let $\alpha^0\in \Z_{\ge 0}^r$ be defined by
$(\alpha^0)_i = \alpha_i$ if $\alpha_i\ge 0$ and
$(\alpha^0)_i = 0$ if $\alpha_i < 0$. Then we put $\ell(\alpha)= \ell(\alpha^0)$.
(Note that $\ell(\alpha,i) = \ell(\alpha^0,i)$ if $\alpha_i\ge0$, and
$\ell(\alpha,i)=0$ if $\alpha_i<0$; in particular,
$\ell(\alpha+\e_i)=\ell(\alpha)+\ell(\alpha,i)$ for every $\alpha\in\Z^r$.)
The
definition of $\lambda(\alpha)$ is the same:
$\lambda(\alpha)=\ell(\alpha+\e)-\ell(\alpha)$.

Note that, for $\gamma\in\Z^{r-1}$, one has
$\ell((a,\gamma),j)=\ell((0,\gamma),j)=\ell_{S^1}(\gamma,j)$ (the latter computed in
$S^1 = S_{I\setminus\{1\}}$) for $j\ge2$ and $a\le 0$. In particular
$\lambda((-1,\gamma))=\lambda_{S^1}(\gamma)$.

\medskip

{\bf (2)} Notice that $0\le \lambda(\alpha)\le r$ for all $\alpha\in \Z^r$,
and $\lambda(\alpha)=0$ if $\alpha_i < 0$ for all $i\in I$. If $\alpha\in S$ then
$\lambda(\alpha)\ge 1$.
Moreover, if $\alpha_j \ge \delta_j$ for
some
$j\neq i$, then $\ell(\alpha,i) = \ell(\w \alpha,i)$, where $\w \alpha$ is the element equal
to $\alpha$ but with $\delta_j$ instead of $\alpha_j$. Moreover,
$\ell(\alpha,i)=\ell(\delta,i)=1$ if $\alpha_i\ge \delta_i$. As a consequence,
$\lambda(\alpha)=\lambda(\delta)=r$ for any $\alpha\geq \delta$.

\medskip

{\bf (3)}
One has
$$
\ell(\alpha+\e)= \ell(\alpha) + \sum_{k=1}^r \ell(\alpha + \e - \e_{\{1,\ldots,k\}},k)
= \ell(\alpha) + \sum_{k=1}^r \ell(\alpha + \e_{\{k+1,\ldots,r\}},k)\; .
$$
As a consequence,
\begin{equation}\label{eqn1}
\lambda(\alpha) =
\sum_{k=1}^r \ell(\alpha + \e -  \e_{\{1,\ldots,k\}},k)
= \sum_{k=1}^r \ell(\alpha + \e_{\{k+1,\ldots,r\}},k)\; .
\end{equation}
The same formulae are valid for any reordering $\{\seq i1r\}$ of the set of
indices $\{1,\ldots, r\}$.

\medskip

{\bf (4)} Let $\alpha, \beta\in \Zr r$ with $\alpha \le \beta$ and let
$k = \abs{\beta-\alpha} = \sum_{i=1}^r (\beta_i-\alpha_i)$.
A {\it path\/} from $\alpha$ to $\beta$ is a sequence
$$
\alpha = \alpha^0 < \alpha^1 < \cdots < \alpha^k = \beta
$$
such that for  $i=1,\ldots, k$ one has
$\alpha^{i}-\alpha^{i-1} = \e_{j_i}$ for some
$j_i\in I$.
Then
\begin{equation}\label{path}
\ell(\beta) = \ell(\alpha) + \sum_{i=1}^k \ell(\alpha^{i-1},j_i)
= \ell(\alpha) + \sum_{i=1}^k \ell(\alpha+\e_{j_1}+\cdots+\e_{j_{i-1}},j_i)\; .
\end{equation}
\end{remarks}

\begin{lemma}\label{lem-2}
For any $\alpha\in \Z^r$ and $i\in I$ one has
\begin{enumerate}
\item
$\ell(\alpha,i) + \ell(\delta-\alpha-\e_i,i)\le 1$.
\item
$\lambda(\alpha)+\lambda(\delta-\alpha-\e)\le r$.
\end{enumerate}
\end{lemma}

\begin{proof}
Let us assume that $\ell(\alpha,i)= \ell(\delta-\alpha-\e_i,i)=1$; then one can pick
$\beta\in \w{\Delta}_i(\alpha)$ and $\gamma\in\w{\Delta}_i(\delta-\alpha-\e_i)$.
One has that $\beta_i+\gamma_i = \delta_i-1$ and
$\beta_j + \gamma_j\ge \alpha_j+\delta_j -\alpha_j > \delta_j-1$ for $j\neq i$. Thus
$\beta+\gamma\in \Delta_i(\tau)$, which contradicts the fact that
$\Delta(\tau)=\emptyset$.
The second inequality follows from the first one: take the path from $\alpha$ to
$\alpha+\e$ in the order $1,2,\ldots,r$ and the reversed one from
$\delta-\alpha-\e$ to $\delta-\alpha$. Then apply the inequality of the first
statement to the pairs
$\ell(\alpha+\e_{\{1,\ldots,k-1\}},k)$ and $\ell(\delta-\alpha-\e_{\{1,\ldots,k\}},k)$
for
$k=1,\ldots,r$.
\end{proof}

\subsection{The Poincar\'e series of \texorpdfstring{$S$}{S}}
Consider the polynomial ring
$\Z[\seq t1r]$ and let ${\cal L}$ be the set of Laurent series in
$\seq t1r$ (i.e.\ formal sums $\sum_{\alpha\in\Z^r}a(\alpha)\tt^{\alpha}$,
$a(\alpha)\in\Z$); ${\cal L}$ carries a canonical structure of module
over the ring $\Z[\seq t1r, t^{-1}_1, \ldots, t^{-1}_r]$.
For $\alpha\in \Z^r$ we denote by $\tt^{\alpha}$ the monomial
$t_1^{\alpha_1}t_2^{\alpha_2}\cdots t_r^{\alpha_r}$.

Let $S$ be a good semigroup.
Associated with the semigroup $S$ we can define the Laurent series
$$
L_S(\seq t1r) = \sum_{\alpha\in \Z^r} \lambda(\alpha) t_1^{\alpha_1}\cdots
t^{\alpha_r}_r \; ,
$$
and also
\begin{equation} \label{Poinc-1}
P'_S(\seq t1r) = (t_1-1)\cdots (t_r-1) L_S(\seq t1r)\; .
\end{equation}

If we write $P'_S(\seq t1r) = \sum_{\alpha \in \Z^r}
p'(\alpha)\tt^{\alpha}$, then for every $\alpha\in \Z^r$ one has
$$
p'(\alpha) =(-1)^r \sum_{p=0}^r (-1)^p \sum_{1\le i_1 <
\cdots<i_p\le r} \lambda(\alpha - \e_{i_1}-\cdots-\e_{i_p})=
(-1)^r \sum_{J\subseteq I} (-1)^{\# J} \lambda(\alpha-\e_J)\; .
$$
As in \cite{CDK}, one has that
$P'_S$ is a polynomial with
integer coefficients. Moreover,
if $p'(\alpha) \neq 0$ then $0\le \alpha \le \delta$ and $p'(\0)=-1$.

\begin{definition}\label{def:PS}
The Poincar\'e series of $S$ is the series
$$
P_S(\tt)  = \frac{P'_{S}(\tt)}{t_1 t_2 \cdots t_r -1} =
\sum_{\alpha\in \Zr r} p(\alpha)\tt^{\alpha} \in \Z[[\seq t1r]]\; .
$$
\end{definition}
Note that,
since $t_1\cdots t_r-1$ is a unit of $\Z[[\tt]]$, $P_S$ is well defined; moreover
$P_S(\0)=1$ (because $p'(\0)=-1$).

Let $\alpha\in \Z^r$ and fix an index $i\in I$. Define
$p_i(\alpha)$ as the integer
$$
p_i(\alpha) =(-1)^{r-1} \sum_{J\subseteq I\setminus \{i\}} (-1)^{\# J}
\ell(\alpha + \e -\e_i -\e_J,i)
= (-1)^r\sum_{i\in K\subseteq I} (-1)^{\# K} \ell(\alpha + \e -
\e_K,i)\;
$$
and let
$$
P_i(\seq t1r) = \sum_{\alpha\in \Z^r} p_i(\alpha) \tt^{\alpha}\; .
$$

\begin{proposition}\label{theo-1}
Let $\alpha\in \Z^r$ and $i\in I$. Then
$p'(\alpha) = -p_i(\alpha) + p_i(\alpha-\e)$ and thus
$P_i(\tt) = P_S(\tt)$.
\end{proposition}

\begin{proof}
Let us first prove that $P_i(\tt)$ is a polynomial.
If $\alpha_i<0$ then all the terms vanish. If $\alpha_k<0$ (or
$\alpha_k\ge\delta_k$) for some $k\neq i$, the terms indexed by $K\cup\{k\}$ and
$K\setminus \{k\}$ cancel, because $\ell(\cdot,i)$ does not see the $k$-th coordinate
when it is $\le0$ (resp.\ $\ge\delta_k$) (see Remarks~\ref{colength}(1),(2)).
Lastly, if $\alpha_i\ge\delta_i$, then all the $\ell$'s involved are $1$
and $\sum_{i\in K}(-1)^{\#K}=0$.
Hence $p_i(\alpha)\neq0$ implies $\0\le\alpha\le\tau$ and $P_i(\tt)$ is a polynomial.

\medskip

The coefficient $p'(\alpha)$ also admits the following expression:
\begin{equation}\label{Poinc-2}
p'(\alpha) = (-1)^r\sum_{i\notin J \subseteq I} (-1)^{\# J}
(\lambda(\alpha -\e_J)-\lambda(\alpha-\e_i-\e_J))\; .
\end{equation}
Let $i_1, \ldots, i_{r-1}$ be any ordering of the set
$I-\{i\}$. Then one can compute each summand above as
$$
\lambda(\alpha-\e_J) = \sum_{p=0}^{r-2} \ell(\alpha -\e_J +
\e_{i_1}+\cdots+\e_{i_p},i_{p+1}) + \ell(\alpha -\e_J+ \e_{I-\{i\}},i)
$$
$$
\lambda(\alpha-\e_J-\e_i) =  \ell(\alpha -\e_J-\e_i,i) +
\sum_{p=0}^{r-2} \ell(\alpha -\e_J +
\e_{i_1}+\cdots+\e_{i_p},i_{p+1})
$$
and as a consequence one has that
$$
p'(\alpha) = (-1)^r\sum_{i\notin J \subseteq I} (-1)^{\# J}
(\ell(\alpha -\e_J-\e_i+\e,i)-\ell(\alpha-\e_i-\e_J,i)) = - p_i(\alpha) +p_i(\alpha-\e)
\;.
$$
This equality says that $P'_S=(\tt^{\e}-1)P_i$. Since
$P_i\in\Z[[\tt]]$ and $\tt^{\e}-1$ is a unit of $\Z[[\tt]]$, we get the
required equality $P_i=P'_S/(\tt^{\e}-1)=P_S$.
\end{proof}

\begin{remark}
Notice that,
as a consequence, $P_i(\tt)$ is independent of $i$ and the series
$P_S(\seq t1r)$ is in fact a polynomial if $r\ge 2$.
Note also that, if $r=1$,
$P_S(t) = \sum_{s\in S}t^s$.
\end{remark}

\begin{remark}
Some coefficients of the
Poincar\'e series can be deduced
from the expression of the coefficients of $P'_S$.
As an example: if $p(\alpha)\neq 0$ then $\alpha$ is a maximal element.
In the next subsection we compute the coefficients in a different
way. In particular, the above facts will be also obtained as a consequence
(Proposition~\ref{prop:coef}).
\end{remark}

\subsection{Simplicial complexes and coefficients of the Poincar\'e series.}
\label{EulerCh}

\begin{definition}
Let $\alpha\in \Z^r$. The {\it configuration space} of $\alpha$ is
the subset of $\PP(I)$ defined~as
$$
\Sigma(\alpha) :=\{ F\subset I, F\neq I | \Delta^F(\alpha)\neq \emptyset \}\subset
\PP(I).
$$
\end{definition}

\begin{remarks}\label{config}
One has the following properties of $\Sigma(\alpha)$,
$\alpha\in \Z^{r}$:
\begin{enumerate}
\item
$\alpha\in S$ if and only if $\emptyset \in \Sigma(\alpha)$.
\item
Let $\alpha\in \Z^{r}$.
If $F, G\in \Sigma(\alpha)$ then $F\cap G\in \Sigma(\alpha)$. As a consequence,
$\Sigma(\alpha)$ is closed under finite intersections.
(Let $\beta\in \Delta^F(\alpha)$ and $\gamma\in \Delta^G(\alpha)$. Then
$\beta\wedge \gamma\in \Delta^{F\cap G}(\alpha)$.)
\item
Let $\alpha\in S$ and
let $F\in \Sigma(\alpha)$, $F\neq \emptyset$,  and $i\notin F$. Then there exists
$G\subset I\setminus F$
with $i\in G$ and $G\in \Sigma(\alpha)$. In particular, $\Sigma(\alpha)$ is a covering of
$I$
whenever $\Sigma(\alpha)$ contains some nonempty set.

Indeed,
let $\beta\in \Delta^F(\alpha)$; the element $\gamma = \alpha\vee_i \beta$ satisfies
$\gamma_i>\alpha_i$ and $\gamma_j = \alpha_j$ for $j\in F$. Thus
$\gamma\in \Delta^G(\alpha)$ for some $G\subset I\setminus F$ with $i\in G$.
Note also that $G\neq I$ because $F\neq\emptyset$.
Obviously, this implies that $\Sigma(\alpha)$ is a covering of $I$ if
$\Sigma(\alpha)\neq \{\emptyset\}$.

\item
$\alpha$ is maximal if and only if $\emptyset \in \Sigma(\alpha)$ and
$I\setminus \{i\}\notin \Sigma(\alpha)$ for all $i\in I$.
Moreover, $\alpha$ is an absolute maximal if and only if
$\Sigma(\alpha)=\{\emptyset\}$ and a relative one if and only if
$\Sigma(\alpha) = \{F :  \# {F} <r-1\}$.
\end{enumerate}
\end{remarks}

\begin{lemma}\label{lem:step}
Let $\alpha\in\mathbb{Z}^r$, $A\subseteq I$ and $j\in I\setminus A$. Then
$\ell(\alpha+\e_A,j)=1$ if and only if there exists $F\in\Sigma(\alpha)$ with
$A\subseteq F$ and $j\notin F$.
\end{lemma}

\begin{proof}
Suppose $\ell(\alpha+\e_A,j)=1$ and let $\beta\in\widetilde\Delta_j(\alpha+\e_A)$, i.e.\
$\beta\in S$, $\beta_j=\alpha_j$, $\beta_i\ge\alpha_i+1$ for $i\in A$ and
$\beta_k\ge\alpha_k$ for $k\notin A\cup\{j\}$. Put $F=\{i\in I:\beta_i>\alpha_i\}$. Then
$A\subseteq F$, $j\notin F$ (so $F\neq I$) and $\beta\in\Delta^F(\alpha)$, hence
$F\in\Sigma(\alpha)$.

Conversely, let $F\in\Sigma(\alpha)$ with $A\subseteq F$ and $j\notin F$, and let
$\beta\in\Delta^F(\alpha)$. Then $\beta_j=\alpha_j$, $\beta_i>\alpha_i$ for $i\in A\subseteq F$
and $\beta_k\ge\alpha_k$ for every $k$, so $\beta\in\widetilde\Delta_j(\alpha+\e_A)$ and
$\ell(\alpha+\e_A,j)=1$.
\end{proof}

\begin{proposition}\label{prop:chain} Let $\alpha\in \Zr r$.
The dimension $\lambda(\alpha)$ of $\alpha$ is the
maximal $k$ such that there exists a strict chain
$ F_1\subsetneq F_2\subsetneq\cdots \subsetneq F_{k}\subsetneq I$ in
$\Sigma(\alpha)$.

In particular, if $\alpha\in S$, then
the nonempty elements of $\Sigma(\alpha)$ form a partition of $I$ if and only if
$\lambda(\alpha)=2$.
\end{proposition}

\begin{proof}
Given a chain $F_1\subsetneq\cdots\subsetneq F_k$ in $\Sigma(\alpha)$, arrange $I$ by listing
first the elements of $F_1$, then those of $F_2\setminus F_1$, \dots, and finally those of $I\setminus F_k$, and
consider the path formula (\ref{path}) from $\alpha$ to $\alpha+\e$ which adds the
indices in this order.
For $m=1,\ldots,k$, the step
that adds the first element $j$ of $F_{m+1}\setminus F_m$ (with $F_{k+1}:=I$) has
$A=F_m$ and so, by Lemma~\ref{lem:step},
$\ell(\alpha+\e_A,j)=1$. Thus
there are at least $k$ summands with value $1$, and so $\lambda(\alpha)\ge k$.

Let us prove that \emph{$\lambda(\alpha)$ is the length of some chain.}
Fix an arbitrary ordering $(i_1,\dots,i_r)$ of $I$,
let $A_k=\{i_1,\dots,i_{k-1}\}$ and call
$\ell(\alpha+\e_{A_k},i_k)$ the $k$-th step of the path formula (see Remarks \ref{colength}, equation
(\ref{path})) for this ordering.
Let $s_1<s_2<\dots<s_m$ be the
steps which are equal to $1$, so that $m=\lambda(\alpha)$. If $m=0$
there is nothing to prove, so assume $m\ge1$. For $t=1,\dots,m$,
Lemma~\ref{lem:step} applied to $A_{s_t}$
provides $F^{(t)}\in\Sigma(\alpha)$ with
\[
A_{s_t}\subseteq F^{(t)}\qquad\text{and}\qquad i_{s_t}\notin F^{(t)} .
\]
Define
\[
G_t:=F^{(t)}\cap F^{(t+1)}\cap\dots\cap F^{(m)}\qquad(t=1,\dots,m).
\]
By Remark \ref{config}, part 2.,
$G_t\in\Sigma(\alpha)$, and obviously
$G_1\subseteq G_2\subseteq\dots\subseteq G_m$. We claim that these inclusions are strict.
Since $s_t<s_u$ implies $A_{s_t}\subseteq A_{s_u}\subseteq F^{(u)}$, we have
$A_{s_t}\subseteq F^{(u)}$ for all $u\ge t$, hence $A_{s_t}\subseteq G_t$. On the other
hand $i_{s_t}\notin F^{(t)}\supseteq G_t$, so $i_{s_t}\notin G_t$, whereas, for $t<m$,
$i_{s_t}\in A_{s_t+1}\subseteq A_{s_{t+1}}\subseteq G_{t+1}$. Therefore
$i_{s_t}\in G_{t+1}\setminus G_t$ and $G_t\subsetneq G_{t+1}$ for $t=1,\dots,m-1$. Thus
$G_1\subsetneq G_2\subsetneq\dots\subsetneq G_m$ is a chain of length $m=\lambda(\alpha)$
in $\Sigma(\alpha)$.

\smallskip
Finally, if $\alpha\in S$ then
$\emptyset\in\Sigma(\alpha)$ and $\Sigma(\alpha)$ is a covering of $I$ as soon as it
contains a nonempty set (see Remarks
\ref{config}). Thus, by the first part and Remarks~\ref{config}(2), $\lambda(\alpha)=2$ if and only if
$\Sigma(\alpha)\neq\{\emptyset\}$ and the intersection of two
distinct nonempty elements of $\Sigma(\alpha)$ is empty, that is, if and only
if the nonempty elements of $\Sigma(\alpha)$ form a partition of $I$.
\end{proof}

\begin{lemma}\label{lem:lambda-max}
Let $\alpha\in \Zr r$. If $\Delta(\alpha)=\emptyset$ then $\lambda(\alpha)< r$, and if
$\alpha\in S$ then $\lambda(\alpha)\ge 1$.
Moreover, if $\alpha$ is an absolute maximal then $\lambda(\alpha)=1$, and if
$\alpha$ is maximal then
$\alpha$ is a relative maximal if and only if $\lambda(\alpha)=r-1$.
\end{lemma}

\begin{proof}
Note that for $k\in I$ one has
$\w \Delta_k(\alpha+\e-\e_k) = \Delta_k(\alpha)$, and so
$\ell(\alpha+\e-\e_k,k)=0$ if $\Delta_k(\alpha)=\emptyset$.
Using this fact and equation (\ref{eqn1}) one easily proves the first two statements
of the Lemma.

The condition
$\Delta(\alpha)=\emptyset$ means $I\setminus\{i\}\notin\Sigma(\alpha)$ for
all $i$, so every chain in $\Sigma(\alpha)$ has at most $r-1$ members.
If $\alpha$ is an absolute maximal then
$\Sigma(\alpha)=\{\emptyset\}$ and so $\lambda(\alpha)=1$.
If $\alpha$ is a relative maximal then
$\Sigma(\alpha)=\{F:\#F\le r-2\}$ and then $\lambda(\alpha)=r-1$.

For the converse of the last
statement: if $\alpha$ is maximal and $\lambda(\alpha)=r-1$, then for \emph{every}
ordering the first $r-1$ steps of the path formula (\ref{path}) are equal to $1$; the $(r-1)$-th
step, with $A$ an arbitrary set of cardinality $r-2$ and $j\notin A$, gives by
Lemma~\ref{lem:step} some $F\in\Sigma(\alpha)$ with $A\subset F$, $j\notin F$, and
$\#F\le r-2$ by maximality, hence $F=A$. Thus all the subsets of cardinality $r-2$ are
in $\Sigma(\alpha)$ and, by Remarks~\ref{config}(2), so are all the subsets of
cardinality $\le r-2$, i.e.\ $\alpha$ is a relative maximal.
\end{proof}

\begin{definition}\label{def:nabla}
Let $\alpha\in \Z^r$ and fix an index $i\in I$.
Let $\Sigma_i(\alpha):= \{G\in \Sigma(\alpha) \st i\notin G\} =
\{G\subset I\setminus \{i\} \st \Delta^{G}(\alpha)\neq \emptyset \}$.
The {\it simplicial complex\/} of $\alpha$ (with respect to $i$) is
the abstract simplicial complex generated by $\Sigma_i(\alpha)$, i.e.\
$$
\nabla_{\alpha,i}:=\{F\subset I\setminus \{i\} \st \exists G\in \Sigma_i(\alpha)
\text{ with } F\subset G\} \; .
$$
Its  {\it combinatorial} Euler characteristic is defined by
$$
\chi ( \nabla_{\alpha,i}) = \sum_{X\in \nabla_{\alpha,i}} (-1)^{\# X}\; .
$$
\end{definition}

\begin{remark}
Note that
the condition $X\in \nabla_{\alpha,i}$ is equivalent to the existence of
$\beta\in S$  such that
$\beta_j>\alpha_j$ for all $j\in Y$ and $\beta_k=\alpha_k$ for $k\notin Y$, for some
$Y$ with $X\subset Y$, $i\notin Y$. So, it is also equivalent to the existence of
$\beta\in S$ such that $\beta_j\ge \alpha_j + 1$ for all $j\in X$,
$\beta_i=\alpha_i$ and $\beta_k\ge \alpha_k$ for all $k\notin X$. Thus, one has that
\begin{equation}\label{eqEuler}
X\in \nabla_{\alpha,i}\iff \w{\Delta}_i(\alpha+\e_{X})\neq \emptyset
\iff \ell(\alpha+\e_{X},i) =1 \;
\end{equation}
and as a consequence:
\begin{equation}\label{eqEuler-2}
\nabla_{\alpha,i}=\{F\subset I\setminus \{i\} \st \w \Delta_i(\alpha+\e_F)\neq
\emptyset\}
=\{F\subset I\setminus \{i\} \st \ell(\alpha+\e_F,i)=1\}\; .
\end{equation}
\end{remark}

Notice that  the reduced (or augmented) Euler characteristic of a realization of the
simplicial complex $\nabla_{\alpha,i}$ is the integer
$$
\w \chi ( \nabla_{\alpha,i}) = \sum_{X\in \nabla_{\alpha,i}} (-1)^{\# X-1} \; .
$$
Here we prefer to use the {\it combinatorial} Euler
characteristic, i.e.\
$
\chi ( \nabla_{\alpha,i}) =
- \w \chi (\nabla_{\alpha,i})\; .
$

\begin{theorem}\label{coef-PS}
Let $\alpha\in \Z^r$. Then one has
$p(\alpha) =  {\chi}(\nabla_{\alpha,i})$ for any index $i\in I$.
As a consequence (recall that $p(\alpha)=0$ if $\alpha\notin S$, see also
Proposition~\ref{prop:coef}(1) below),
$$
P_S(\tt) =  \sum_{\alpha\in S} {\chi}(\nabla_{\alpha,i}) \;  \tt^{\alpha}\; .
$$
(Note that, in particular, $\chi(\nabla_{\alpha,i})$ does not  depend on the index
$i$.)
\end{theorem}

\begin{proof}
For the coefficients of the Poincar\'e series one has, by Proposition~\ref{theo-1},
\begin{align*}
p(\alpha) & = p_i(\alpha) = (-1)^r \sum_{i\in K\subset I}(-1)^{\# K}
\ell(\alpha + \e_{I-K},i) \\
& = (-1)^{r} \sum_{X\subset I-\{i\}} (-1)^{\# (I-X)} \ell(\alpha+\e_X,i) =
\sum_{X\subset I-\{i\}} (-1)^{\# X}\ell(\alpha+\e_X,i)\; .
\end{align*}

Now, by (\ref{eqEuler}), for each subset $X$ of $I\setminus\{i\}$ one has
$\ell(\alpha+\e_X,i)=1$ if and only if $X\in \nabla_{\alpha,i}$, so that
$$
p(\alpha) = \sum_{X\in \nabla_{\alpha,i}} (-1)^{\# X} =
{\chi}(\nabla_{\alpha,i})\; .
$$
\end{proof}

Writing
$\ell(\alpha+\e_X,i)=\ell(\alpha+\e_X+\e_i)-\ell(\alpha+\e_X)$
and grouping the subsets $Y=X$ and $Y=X\cup\{i\}$ of $I$ in the expression
$p(\alpha)=\sum_{X\subset I-\{i\}} (-1)^{\# X}\ell(\alpha+\e_X,i)$ one gets the
following index-free formula (which makes the independence of $i$ evident):

\begin{corollary}\label{cor:indexfree}
For $\alpha\in \Z^r$ one has
\begin{equation}\label{eq:indexfree}
p(\alpha)=-\sum_{Y\subseteq I}(-1)^{\#Y}\,\ell(\alpha+\e_Y)\; .
\end{equation}
\end{corollary}

\begin{proposition}\label{prop:coef}
Let $\alpha\in \Zr r$. Then:
(1) If $\alpha\notin S$, then $p(\alpha)=0$. (2) If $\alpha\in S$ is not  maximal, then
$p(\alpha)=0$. (3) If $\alpha$ is an absolute maximal, then $p(\alpha)=1$. (4) If
$\alpha$ is a relative maximal, then $p(\alpha)=(-1)^r$.
\end{proposition}

\begin{proof}
{\bf (1)}
If $\alpha\notin S$ there exists $i$ such that
$\w \Delta_i(\alpha)=\emptyset$. For $F\subset I\setminus \{i\}$ one has
$\w \Delta_i(\alpha+\e_F)\subset \w \Delta_i(\alpha)$, so
$\w \Delta_i(\alpha+\e_F) = \emptyset$ and $F\notin \nabla_{\alpha,i}$. Thus
$\chi(\nabla_{\alpha,i})=0$.

{\bf (2)} Let $k$ be such that $\Delta_k(\alpha)\neq \emptyset$ and let
$\beta\in \Delta_k(\alpha)$.
Then $\beta_j> \alpha_j$ for $j\neq k$, and so $\beta\in \w \Delta_k(\alpha+\e_F)$ for
all $F\subset I\setminus \{k\}$. This means that $\nabla_{\alpha,k}$ is the full
simplex on $I\setminus\{k\}$, and $p(\alpha)=\chi (\nabla_{\alpha,k})=0$.

{\bf (3)} If $\alpha$ is an absolute maximal one has $\nabla_{\alpha,i}=\{\emptyset\}$,
and so $p(\alpha)=\chi(\nabla_{\alpha,i})=1$.

{\bf (4)}
If $\alpha$ is a relative maximal, then
$\nabla_{\alpha,i} = \{F \subset I\setminus \{i\} \st \# F\le r-2\}$, and so
$p(\alpha) = \chi(\nabla_{\alpha,i})= \sum_{k=0}^{r-2}(-1)^k\binom{r-1}{k}=(-1)^r$.
\end{proof}

\begin{example}
\label{ex:p-zero}
By the Proposition, the support of the Poincar\'e series is a subset of
the set of maximal elements; however, it is not equal to it, as this example shows.  Let
$S\subset\mathbb{Z}^4_{\ge0}$ be the value semigroup of the plane curve
\[
xy\,(y-x^2)\,(x-y^2)=0,
\]
with branches $C_1=(t,0)$, $C_2=(t,t^2)$, $C_3=(0,t)$,
$C_4=(t^2,t)$: two pairs of tangent smooth branches, the pairs being
transversal to each other. The intersection numbers are
$(C_1 \cdot C_2)=(C_3 \cdot C_4)=2$ and all the others equal $1$, so
$\delta=(4,4,4,4)$ and $\tau=(3,3,3,3)$. The elements of $S$ below $\tau$ are
\begin{gather*}
\underline0,\ (1,1,1,1),\ (1,1,2,2),\ (1,1,2,3),\ (1,1,3,2),\ (2,2,1,1),\\
(2,3,1,1),\ (3,2,1,1),\ (2,2,2,2),\ (2,2,3,3),\ (3,3,2,2),\ \tau ,
\end{gather*}
and the maximal elements, with their $\Sigma$ and their coefficient in $P_S$, are:
\begin{center}
\begin{tabular}{lll}
\toprule
$\alpha$ & $\Sigma(\alpha)$ & $p(\alpha)$\\
\midrule
$(0,0,0,0)$ & $\{\emptyset\}$ & $1$\\
$(1,1,1,1)$ & $\{\emptyset,\{1,2\},\{3,4\}\}$ & $0$\\
$(1,1,2,2)$ & $\{\emptyset,\{3\},\{4\},\{1,2\}\}$ & $-1$\\
$(2,2,1,1)$ & $\{\emptyset,\{1\},\{2\},\{3,4\}\}$ & $-1$\\
$(2,2,2,2)$ & $\{\emptyset,\{1,2\},\{3,4\}\}$ & $0$\\
$(3,3,3,3)$ & all $F$ with $\#F\le2$ & $1$\\
\bottomrule
\end{tabular}
\end{center}
Hence
\[
P_S(\underline t)=1-t_1t_2t_3^2t_4^2-t_1^2t_2^2t_3t_4+t_1^3t_2^3t_3^3t_4^3 \; .
\]
Thus,
the two maximal elements $(1,1,1,1)$ and
$(2,2,2,2)=\tau-(1,1,1,1)$ do not appear in the support of $P_S$.
Note that the semigroup is symmetric (see the next section), being the semigroup of a
plane curve.

Notice that in this case the semigroup $S$ is equivalent to the equisingularity type of
the plane curve (i.e.\ to its topological type) and
$P_S(\tt)$ coincides with the Alexander polynomial
of the link of the plane curve (see \cite{Duke}). A well-known theorem of Yamamoto (see
\cite{Yamamoto}) states that the Alexander polynomial determines the topological
type of the curve, so that, in this case, $S$ can be recovered from $P_S$; the key
ingredient of the proof is the projection formula of Torres (see \cite{Torres}). In
particular, the two missing elements $(1,1,1,1)$ and $(2,2,2,2)$ can be recovered from
$P_S$ by other means, e.g.\ through the projections $S_{\{1,2\}}$ and $S_{\{3,4\}}$.
We discuss the projection formula and the recovery of $S$ from $P_S$ in the context
of good semigroups
in Sections~\ref{Supersymmetric} and~\ref{sec:recovering}.
\end{example}

\begin{remark}\label{rem:euler}
Let $S$ be the value semigroup of a reduced curve
singularity $C$ with $r$ branches, and let $\OO$, $v$ and $J(\alpha)$ be as in
Remark~\ref{rem:curves}. The vector space
$V_\alpha=J(\alpha)/J(\alpha+\e)$ can be seen in a natural way as a vector subspace of
$\C^r$ of dimension $\lambda(\alpha)$, by mapping the class of $g$ to the vector of
the coefficients of $t_1^{\alpha_1},\dots,t_r^{\alpha_r}$ in the expansions of $g$
along the branches. Let
$F_\alpha\subset V_\alpha$ be the image of $\{g\in\OO : v(g)=\alpha\}$. Under the
above identification, $F_\alpha=V_\alpha\cap(\C^*)^r$, that is,
the complement in $V_\alpha$ of the union
of the $r$ subspaces
$V_j=J(\alpha+\e_j)/J(\alpha+\e)=V_\alpha\cap\{x_j=0\}$, $j\in I$ (each of them is
either a hyperplane of $V_\alpha$ or $V_\alpha$ itself).
The extended semigroup defined in
\cite{Extended} (see also \cite{Duke}) is
$\widehat S=\bigcup_\alpha\{\alpha\}\times F_\alpha\subset\Z^r\times(\C^*)^r$. Since $F_\alpha$ is
$\C^*$-invariant, its projectivization $\P F_\alpha$ is the complement of a hyperplane
arrangement in $\P(V_\alpha)\simeq\P^{\lambda(\alpha)-1}$, and one relevant
result (see \cite[Theorem 2]{Duke}) is that
$$
P_S(\tt)=\sum_{\alpha\in\Z^r_{\ge0}}\chi(\P F_\alpha)\,\tt^{\alpha}\; ,
$$
where $\chi$ is the topological Euler characteristic. By
inclusion--exclusion, using $\chi(\P^{d-1})=d$ and
$\bigcap_{j\in Y}V_j=J(\alpha+\e_Y)/J(\alpha+\e)$,
$$
\chi(\P F_\alpha)=\sum_{Y\subseteq I}(-1)^{\#Y}\big(\ell(\alpha+\e)-\ell(\alpha+\e_Y)\big)
=-\sum_{Y\subseteq I}(-1)^{\#Y}\,\ell(\alpha+\e_Y)\; ,
$$
which is exactly the formula (\ref{eq:indexfree}) for $p(\alpha)$. Thus Theorem~\ref{coef-PS} is the combinatorial
counterpart of the formula of \cite{Duke}: the simplicial complex
$\nabla_{\alpha,i}$ records which intersections $\bigcap_{j\in X}V_j$ are not
contained in $V_i$, and its Euler characteristic computes $\chi(\P F_\alpha)$.
It is also remarkable that
$P_S$ can be seen as the integral of $\tt^{v(g)}$ over $\P\OO$ with respect to
the Euler characteristic, see \cite{IJM-03}.
\end{remark}

\section{Symmetric Semigroups}\label{sec:Symmetry}

The symmetry of the
semigroup of values is the combinatorial counterpart of the Gorenstein property of
the ring of the curve (Remark~\ref{numerical}). In this section we characterize
symmetric good semigroups in terms of the colength and the dimension
(Proposition~\ref{symmetric}) and of an Alexander duality between the simplicial
complexes $\nabla_{\alpha,i}$ and $\nabla_{\tau-\alpha,i}$
(Proposition~\ref{teo-dualidad}). As a consequence, we obtain the functional
equation of the Poincar\'e series and the exchange of absolute and relative
maximals under $\alpha\mapsto\tau-\alpha$ (Theorem~\ref{thm:duality-maximals}), and
we study the numbers $\mu^i,\chi_i,\kappa_i$ which measure the fibres of the
projections $S\to S_i$ (Proposition~\ref{sym-1} and Theorem~\ref{sym-2}).

\medskip

Let $S\subset \Z_{\ge 0}^r$ be a good semigroup.
Let $c_i = c(S_i)$ denote the conductor of the numerical semigroup $S_i$
$(1\le i\le r)$, $\delta = (c_1+\mu^1, \ldots, c_r+\mu^r)$ the conductor of the
semigroup $S$ and $\tau =\delta-\e$ its Frobenius vector.

\begin{definition}\label{def:symmetric}
The semigroup $S\subset \Zr r$, $r\ge 2$, is said to be \emph{symmetric} if for
any $\alpha\in \Z^r$ one has
$$
\alpha\in S  \iff \Delta(\tau-\alpha)=\emptyset \; .
$$
\end{definition}

Note that if $\alpha\in S$ then $\Delta(\tau-\alpha)=\emptyset$ trivially, because
$\Delta(\tau)=\emptyset$ (if $s\in\Delta_i(\tau-\alpha)$ then $s+\alpha\in\Delta_i(\tau)$),
so the only new condition is the converse one.

\begin{remark}\label{numerical}
The above definition is the generalization of the symmetry of numerical semigroups.
If $S\subset \Zr {}$ is a numerical semigroup and $c$ is its conductor, $S$ is
symmetric  if for any $\alpha\in \Z$ one has that $\alpha\in S$ if and only if
$c-1-\alpha\notin S$. This means that the number of gaps of the
semigroup coincides with the number of elements of $S$ below $c$. In particular,
$c$ must be even.
A well-known theorem of Kunz (see \cite{Kunz}) shows that,
for an irreducible curve singularity, the ring of functions is
Gorenstein if and only if its semigroup $S$ is symmetric. The reason is that
$\ell(c) = \# [0,c-1]\cap S$ and symmetry is equivalent to $\ell(c)=c/2$.

This result is extended to reduced curves (i.e.\ $r\ge 2$)
in \cite{D-88} (see also \cite{CDK}), where it is proved that the curve is Gorenstein if
and only if $S$ is symmetric.
\end{remark}

The symmetry can also be characterized in terms of the colength and the dimension
of the elements of $S$. More precisely:

\begin{proposition}\label{symmetric}
The following conditions are equivalent:
\begin{enumerate}
\item
$S$ is symmetric.
\item
$\ell(\alpha,i) + \ell(\delta-\alpha-\e_i,i)=  1$ for all $\alpha\in \Z^r$ and $i\in
I$.
\item $\lambda(\alpha)+\lambda(\delta-\alpha-\e) = r$ for all $\alpha\in \Z^r$.
\item
$2 \ell(\delta) = \abs{\delta} = \delta_1+\cdots +\delta_r$\; .
\end{enumerate}
\end{proposition}

\begin{proof}
The equivalence of (2) and (3) is a trivial consequence of Lemma
\ref{lem-2}. Let us prove the equivalence of (1) and (2).

Let us assume that $S$ is symmetric and suppose that  $\ell(\alpha,i)=0$,
i.e.\ $\w{\Delta}_i(\alpha)=\emptyset$.
Then, by Lemma~\ref{lem:new}, there exists $\beta\in \Z^r$ such that
$\alpha_i=\beta_i$, $\beta_j<\alpha_j$ for $j\neq i$ and $\Delta(\beta)=\emptyset$.
Then $\gamma = \delta-\e-\beta\in S$. However, for the element $\gamma$ one has that
$\gamma_i = \delta_i-1-\beta_i = \delta_i-1-\alpha_i$ and, if $j\neq i$,
$\gamma_j = \delta_j -1 -\beta_j > \delta_j-1-\alpha_j$, i.e.\
$\gamma_j\ge\delta_j-\alpha_j$, and so
$\gamma\in \w{\Delta}_i(\delta-\alpha-\e_i)\neq \emptyset$, i.e.\
$\ell(\delta-\alpha-\e_i,i)=1$.

Assume that condition (2) holds and let $\alpha\in \Z^r$ be such that
$\Delta(\alpha)=\emptyset$. In this case one has
$\w{\Delta}_i(\alpha+\e-\e_i) =  \Delta_i(\alpha)=\emptyset$, which means that
$\ell(\alpha+\e-\e_i,i)=0$. So, $\ell(\delta-\alpha-\e+\e_i-\e_i,i) =
\ell(\tau-\alpha,i)=1$, and therefore
$\w{\Delta}_i(\tau-\alpha)\neq \emptyset$ for any $i\in I$. Now,
if $\beta^i\in \w{\Delta}_i(\tau-\alpha)$, one has that
$\tau-\alpha = \inf \{\beta^i : i=1,\ldots,r\}\in S$.

The equivalence of (4) and (2) can be proved following the proof of
Theorem (3.6) in \cite{CDK} (see also \cite{D-88}). The idea is to pair each step of
a path $\0 = \alpha^0 < \cdots < \alpha^{n}=\delta$, $n=\abs{\delta}$, with the
corresponding step of
the ``reflected'' (or symmetric) path $\beta^m = \delta-\alpha^{n-m}$ between $\delta$
and $\0$. In this way one has
$$
2 \ell(\delta) = \sum_{k=1}^{n} \left[\ell(\alpha^{k-1},j_k) +
\ell(\delta-\alpha^{k-1}-\e_{j_k},j_k)\right]
\le \abs{\delta}\; .
$$
Now, (2)$\Rightarrow$(4) is immediate from the displayed
formula. For (4)$\Rightarrow$(2), the equality $2\ell(\delta)=|\delta|$ forces every
pair in every path to add up to $1$ (Lemma~\ref{lem-2}); every pair $(\alpha,i)$ with
$\0\le\alpha\le\delta-\e_i$ occurs in some path, and for the remaining
$\alpha\in\Z^r$ one reduces to this case with Remarks~\ref{colength}(1),(2)
(if $\alpha_i<0$ or $\alpha_i\ge\delta_i$ the sum is trivially $1$; if
$\alpha_j<0$ or $\alpha_j>\delta_j$ for $j\neq i$, replacing $\alpha_j$ by $0$,
resp.\ $\delta_j$, changes neither $\ell(\alpha,i)$ nor $\ell(\delta-\alpha-\e_i,i)$).
\end{proof}

Let
$\alpha\in\Z^r$,  $i\in I$ and $I'=I\setminus\{i\}$. Let
\[
\nabla_{\alpha,i}=\{X\subset I':\ \ell(\alpha+\e_X,i)=1\}
=\{X\subset I':\ \tD_i(\alpha+\e_X)\neq\emptyset\}
\]
be the simplicial complex of $\alpha$ relative to $i$ (see (\ref{eqEuler-2}); for
$\alpha\in S$ this is the complex of Definition~\ref{def:nabla}).
Let us consider its \emph{Alexander dual complex}
\begin{equation}\label{dual}
\nabla^{*}_{\alpha,i}=\{X\subset I':\ I'\setminus X \notin \nabla_{\alpha,i}\}\; .
\end{equation}
Note that
\begin{align*}
\chi(\nabla^{*}_{\alpha,i}) &=
\sum_{X\in\nabla^*}(-1)^{\#X}
= \sum_{Y\notin\nabla}(-1)^{\#(I'\setminus Y)} \\
&= (-1)^{r-1}\big(\sum_{Y\subset I'}(-1)^{\#Y}
-\sum_{Y\in\nabla}(-1)^{\#Y}\big) =(-1)^{r}\chi(\nabla),
\end{align*}
since $I'\neq\emptyset$.

\begin{proposition}\label{teo-dualidad}
Let $S\subset \Zr r$ be a good semigroup, $r\ge 2$. The following conditions are
equivalent:
\begin{enumerate}
\item $S$ is symmetric.
\item $\nabla_{\tau-\alpha,i}={\nabla^{*}_{\alpha,i}}$ for all $\alpha\in\Z^r$ and
all $i\in I$.
\item $\nabla_{\tau-\alpha,i}={\nabla^{*}_{\alpha,i}}$ for all  $\alpha$ with
$\0\le\alpha\le\delta$ and all $i\in I$.
\end{enumerate}
\end{proposition}

\begin{proof}
Let $\rho\in \Z^{r}$ and $X\subset I'$, one has
$\delta-(\rho+\e_X)-\e_i=(\tau-\rho)+\e-\e_X-\e_i=(\tau-\rho)+\e_{I'\setminus X}$.
So, for $\alpha\in \Z^r$:
\[
X\in {\nabla^*_{\alpha,i}}\iff I'\setminus X\notin\nabla_{\alpha,i}
\iff \ell(\alpha+\e_{I'\setminus X},i)=0
\iff \ell\big(\delta-\beta-\e_i,\,i\big)=0\;.
\]
where $\beta=(\tau-\alpha)+\e_X$. On the other hand,
$$
X\in \nabla_{\tau-\alpha,i} \iff \ell(\tau-\alpha+\e_{X},i)=1\iff \ell(\beta,i)=1\;.
$$
So,
condition (2) for $\alpha$ is equivalent to saying that for any $X\subset I'$ one has
$\ell\big(\delta-\beta-\e_i,\,i\big)=0 \iff \ell(\beta,i)=1$ or, equivalently (by
Lemma~\ref{lem-2}), that
$\ell\big(\delta-\beta-\e_i,\,i\big) + \ell(\beta,i)=1$.

Now, it is clear that when $\alpha$ runs over $\Z^r$ (and $X$ over the subsets of
$I'$), $\beta$ also runs over $\Z^r$. That is,
(2) is equivalent to the symmetry of $S$ by Proposition \ref{symmetric}.

\medskip

(2) $\Rightarrow$(3) is trivial, so let us prove that (3)$\Rightarrow$(1).
Let $\alpha\in \Z^r$ with $\Delta(\alpha)=\emptyset$. We have to prove that
$\tau-\alpha\in  S$ (the converse is always true).
Let $\w \alpha\in \Zr r$ be such that $\w \alpha_i=\max\{0,\alpha_i\}$; one also has
that $\Delta(\w \alpha)=\emptyset$, and it is easy to see that if $\tau-\w \alpha\in S$
then $\tau-\alpha\in S$. So, we can assume that $\alpha\ge \0$.

If $\alpha_j\ge\delta_j$ for some $j$, the element $\beta$ with  $\beta_j=\alpha_j$ and
$\beta_k=\max(\alpha_k+1,\delta_k)$ for $k\neq j$ is in $\delta+\Zr r\subset S$
and also in $\Delta_j(\alpha)$, a contradiction; hence $\0\le\alpha\le\tau$.

For each
$i$, $\Delta_i(\alpha)=\emptyset$ implies that $I'\notin\nabla_{\alpha,i}$, i.e.\
$\emptyset\in\nabla^*_{\alpha,i}=\nabla_{\tau-\alpha,i}$.
Then $\ell(\tau-\alpha,i)=1$ and one can take
$\beta^i\in\w\Delta_i(\tau-\alpha)$ for each $i$.
So, $\tau-\alpha=\bigwedge_i\beta^i = \inf\{\beta^i : \ i\in I\}\in S$.
\end{proof}

\begin{example}\label{ex-symmetry}

By Proposition~\ref{teo-dualidad}, $S$ is symmetric if and only if
$\nabla_{\tau-\alpha,i}=(\nabla_{\alpha,i})^{*}$ for every $\alpha\in\mathbb{Z}^r$ and every
$i\in I$. The next two examples, both with $r=2$, show that this condition cannot be
weakened in two natural ways: it does not suffice to require it for a single fixed
index $i$ and every $\alpha$, nor for every $i$ but only for $\alpha\in S$.

\medskip
{\bf (1)}\ \emph{A single index is not enough.}
Let
\[
S=\{(0,0)\}\ \cup\ \{(2,k):k\ge3\}\ \cup\ \{(3,3)\}\ \cup\ \bigl((4,4)+\mathbb{Z}^2_{\ge0}\bigr)
\subset\mathbb{Z}^2_{\ge0}.
\]
$S$ is a good semigroup (see the left part of Figure \ref{fig-1}), with $\delta=(4,4)$
and $\tau=(3,3)$.
Moreover, $S$ is
{\bf not symmetric:} $(3,2)\notin S$ and $\Delta(\tau-(3,2))=\Delta(0,1)=\emptyset$.

It is not hard to see that the
duality (condition 2 of the Proposition) holds for $i=1$ and every
$\alpha\in\mathbb{Z}^2$. Notice that the
duality fails for $i=2$, for instance at $\alpha=(1,1)$ and at
$\alpha=(0,1)$.

\begin{figure}[h]
\centering
\unitlength=1.00mm
\begin{picture}(60,40)(7,7)
\thicklines
\put(10,10){\line(1,0){50}}
\put(10,10){\line(0,1){40}}
\put(6,6){$0$}
\put(25,6){$2$}
\put(33,6){$3$}
\put(41,6){$4$}
\put(6,24){$2$}
\put(6,32){$3$}
\put(6,40){$4$}

\thinlines
\put(18,10){\line(0,1){40}}
\put(26,10){\line(0,1){40}}
\put(34,10){\line(0,1){40}}
\put(42,10){\line(0,1){40}}

\put(10,18){\line(1,0){40}}
\put(10,26){\line(1,0){40}}
\put(10,34){\line(1,0){40}}
\put(10,42){\line(1,0){40}}

\put(10,10){\circle*{2}}
\put(26,34){\circle*{2}}
\put(26,42){\circle*{2}}
\put(26,50){\circle*{2}}
\put(34,34){\circle*{2}}

\put(42,42){\circle*{2}}

\thicklines
\put(42,42){\line(1,0){20}}
\put(42,42){\line(0,1){10}}
\end{picture}
\qquad \qquad
\begin{picture}(60,40)(7,7)
\thicklines
\put(10,10){\line(1,0){40}}
\put(10,10){\line(0,1){40}}
\put(6,6){$0$}
\put(16,6){$2$}
\put(25,6){$4$}
\put(29,6){$5$}
\put(6,16){$2$}
\put(6,24){$4$}
\put(6,28){$5$}

\thinlines
\put(14,10){\line(0,1){25}}
\put(18,10){\line(0,1){25}}
\put(22,10){\line(0,1){25}}
\put(26,10){\line(0,1){25}}
\put(30,10){\line(0,1){20}}

\put(10,14){\line(1,0){25}}
\put(10,18){\line(1,0){25}}
\put(10,22){\line(1,0){25}}
\put(10,26){\line(1,0){25}}
\put(10,30){\line(1,0){20}}

\put(10,10){\circle*{1.5}}
\put(18,18){\circle*{1.5}}
\put(26,26){\circle*{1.5}}

\put(30,30){\circle*{1.5}}

\thicklines
\put(30,30){\line(1,0){20}}
\put(30,30){\line(0,1){20}}
\end{picture}
\caption{The semigroups $S$ and $T$ of Example~\ref{ex-symmetry}
}
\label{fig-1}
\end{figure}
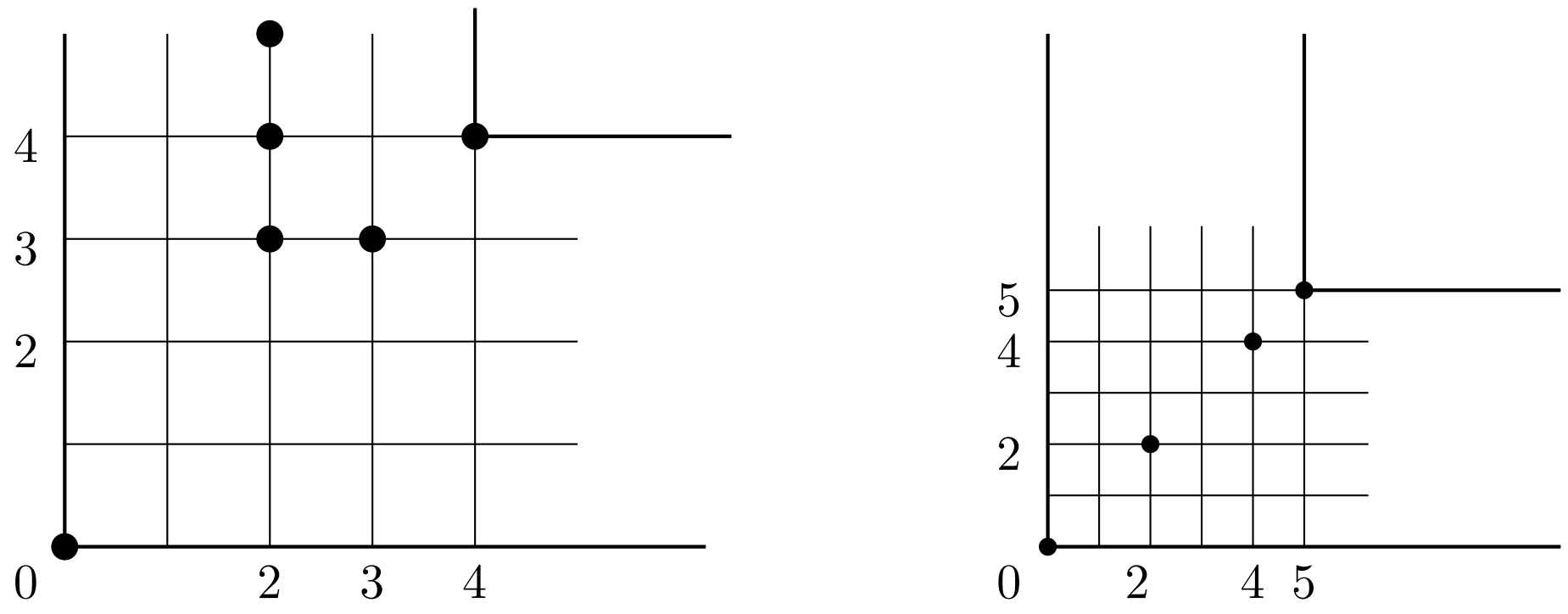

{\bf (2)}\ \emph{Elements of the semigroup are not enough.}
Let
$$
T=\{(0,0),(2,2),(4,4)\}\cup\bigl((5,5)+\mathbb{Z}^2_{\ge0}\bigr)\; .
$$
$T$ is a good
semigroup with $\tau=(4,4)$ (see the right part of Figure \ref{fig-1}). It is not
symmetric: $(1,1)\notin T$ and
$\Delta(\tau-(1,1))=\Delta(3,3)=\emptyset$.

However, $\nabla_{\tau-\alpha,i}=(\nabla_{\alpha,i})^{*}$ for \emph{every} $\alpha\in
T$ and $i=1,2$.
The duality fails, for instance, at $\alpha=(1,1)\notin T$:
$\nabla_{(1,1),i}=\nabla_{(3,3),i}=\emptyset$, whereas $\emptyset^{*}=\{\emptyset,I'\}$.

\end{example}

\begin{theorem}\label{thm:duality-maximals}
Let $S$ be a symmetric good semigroup and let $\alpha\in\mathbb{Z}^r$. Let
$P_S(\tt) = \sum p(\alpha)\tt^{\alpha}$ be its Poincar\'e series.
Then:
\begin{enumerate}
\item[(1)] $\alpha$ is a maximal element of $S$ if and only if $\tau-\alpha$ is a maximal
element of $S$.
\item[(2)] $\alpha$ is an absolute maximal if and only if $\tau-\alpha$ is a relative
maximal.
\item[(3)] $p(\tau-\alpha)=(-1)^r\,p(\alpha)$; equivalently:
$$
P_S(\underline t)=(-1)^r\,\underline t^{\,\tau}\,P_S(t_1^{-1},\dots,t_r^{-1})\; .
$$
\end{enumerate}
\end{theorem}

\begin{proof}
(1) is obvious from the definition of symmetry.

(2) One has that $\alpha$ is an absolute maximal if and only if
$\nabla_{\alpha,i} = \{\emptyset\}$ for every $i\in I$, and that $\tau-\alpha$ is a
relative maximal if and only if
$\nabla_{\tau-\alpha,i} = \{X\subset I\setminus \{i\} : \ \# X\le r-2\}$ for every
$i\in I$
As $\nabla_{\tau-\alpha,i} = \nabla^{*}_{\alpha,i}$ for every $i$,
the result follows from Proposition \ref{teo-dualidad}.

(3) If $\alpha$ and $\tau-\alpha$ both belong to $\mathbb{Z}^r_{\ge0}$, then by
Theorem \ref{coef-PS}
\[
p(\tau-\alpha)=\chi(\nabla_{\alpha,i}^{*})=(-1)^{r}\,
\chi(\nabla_{\alpha,i})=(-1)^r\,p(\alpha).
\]
Otherwise $\alpha\notin[\underline0,\tau]$ and $\tau-\alpha\notin[\underline0,\tau]$, and both
coefficients are equal to $0$. The functional equation for the Poincar\'e series is
just a restatement of the equality $p(\tau-\alpha)=(-1)^rp(\alpha)$ for all $\alpha$.
\end{proof}

\begin{remark}
The functional equation of
Theorem~\ref{thm:duality-maximals}(3) is not new: for Gorenstein curves it goes back
to \cite{CDK}, and is also a particular case of
\cite[Theorem 4]{Tozzo} proved in the more general context
of the duality of good semigroup
ideals (see
\cite{Korell-2019} and \cite{Pol}). The proof given here is different: it is a
consequence of the
Alexander duality between the complexes $\nabla_{\alpha,i}$ and $\nabla_{\tau-\alpha,i}$
(Proposition~\ref{teo-dualidad}), which is in fact equivalent to the symmetry of $S$,
and it gives at the same time the exchange of absolute and relative maximals,
statement (2).
\end{remark}

\begin{remark}
Recall that for $r=2$ all the maximals are absolute, and so
$P_S(\tt) = \sum_{\alpha \text{ absolute}} \tt^{\alpha}$. For $r=3$ absolute maximals
and relative maximals are the only types  of maximals, and
$$
P_S(\tt) = \sum_{\alpha   \text{ absolute}} \tt^{\alpha} - \sum_{\alpha \text{
relative}} \tt^{\alpha}\; .
$$

For $r\ge 4$ there are maximal elements
which are neither absolute nor relative.
For instance, for the plane curve with the four branches $(t,0)$, $(t,t^2)$, $(t^2,t^3)$, $(t,-t)$ one
gets $\tau=(5,5,9,3)$; the maximal elements are $\underline0$ (absolute), $(5,5,9,3)$
(relative), and the pair $(2,2,3,1)$, $(3,3,6,2)$, whose elements are of neither type and
add up to $\tau$. Notice that the last two are also interchanged by the symmetry.
\end{remark}

\begin{remark}\label{rem:PS-sym-not-S}
In Example \ref{ex-symmetry}(1) the set of maximals is symmetric (there are only $(0,0)$
and $\tau=(3,3)$) and the
Poincar\'e series is $P_S(\tt) = 1+\tt^{(3,3)}$,
which is symmetric (with respect to $\tau = (3,3)$). So, this is also an example
showing that
the symmetry of the Poincar\'e series does not imply the symmetry of the semigroup.

For $r=1$, however, the converse holds: if a numerical semigroup $\Gamma$ with
conductor $c$ satisfies $P_\Gamma(t)=-t^{\gamma}P_\Gamma(t^{-1})$ for some
$\gamma\in\Z$, then $\Gamma$ is symmetric (in particular $\gamma$ is the Frobenius
number $c-1$).
Indeed, since $P_\Gamma(t)=\sum_{n\in\Gamma,\,n<c}t^{n}+t^{c}/(1-t)$, one has
$$
-t^{\gamma}P_\Gamma(t^{-1})=\sum_{m\ge\gamma-c+1}t^{m}\;-\sum_{n\in\Gamma,\,n<c}t^{\gamma-n}.
$$
If $\gamma<c-1$, the coefficient of $t^{\gamma-c+1}$ on the right is $1$ (because
$c-1\notin\Gamma$), which is impossible since $\gamma-c+1<0$. If $\gamma>c-1$, the
coefficient of $t^{0}$ on the right is $0$ (because $\gamma\ge c$), whereas it is $1$ on
the left. Hence $\gamma=c-1$, and comparing the coefficients of $t^{m}$, $m\ge0$, gives
$m\in\Gamma$ if and only if $c-1-m\notin\Gamma$.
\end{remark}

\begin{example}
\label{ex:Sprime}
Let $L$ be the value semigroup of the plane curve $y(y-x^2)(y-2x^2)$ (three smooth
branches with pairwise contact equal to 2). It is symmetric, with $\delta(L)=(4,4,4)$, and
\[
L=\{\0,\ \e,\ \ul2,\ \ul3 \}\ \cup\ \{(2,2,k)^\sigma: k\ge3\}\ \cup\
\{(3,3,k)^\sigma: k\ge4\}\ \cup\ (\ul4 +\Zr 3),
\]
where $(a,a,k)^\sigma$ stands for its three permutations. Its maximal elements are
$\0,\e$ (absolute) and $\ul 2,\ul 3$ (relative), and they satisfy the duality of
Theorem~\ref{thm:duality-maximals}. Define
\[
S=\{\0\}\cup\big((L\setminus\{\ul0\})+\e\big)
=\{\0,\ul2,\ul3,\ul4 \}
\cup\{(3,3,k)^\sigma:k\ge4\}\cup\{(4,4,k)^\sigma:k\ge5\}\cup(\ul5+\Zr 3).
\]
$S$ is a (non-symmetric) good semigroup with $\delta=\ul 5=(5,5,5)$ and
$\tau=\ul4$. The maximal elements of $S$ are:
\begin{center}
\begin{tabular}{cccc}\toprule
$\alpha$ & type & $\tau-\alpha$ & type of $\tau-\alpha$\\\midrule
$\0$ & absolute & $\ul 4$ & relative\\
$\ul2$ & \textbf{absolute} & $\ul2$ & \textbf{absolute (not relative)}\\
$\ul3$ & \textbf{relative} & $\e$ & $\notin S$\\
$\ul4$ & relative & $\0$ & absolute\\\bottomrule
\end{tabular}
\end{center}
Hence $\ul 2$ is absolute maximal but $\tau-\ul2$ is not
relative maximal, and $\ul 3$ is relative maximal while $\tau-\ul 3$ does not
even belong to $S$.
Moreover, one gets $P_{S}=1+\tt^{\ul 2}-\tt^{\ul3}-\tt^{\ul 4}$. Here
$p(\ul2)=1\ne(-1)^3p(\ul2)$ and $p(\ul3 )=-1\neq(-1)^3p(\e)=0$.
Thus, if $S$ is not symmetric, all the implications of
Theorem~\ref{thm:duality-maximals} may fail.

\medskip
The construction $L\mapsto\{\0\}\cup((L\setminus\{\0\})+\e)$ is a simple way of
breaking the symmetry while preserving the structure of the maximal elements: the nonzero
maximal elements of $L$ are translated, but $\tau$ increases by $\e$, so the duality no
longer matches.
\end{example}

\medskip

Theorem \ref{thm:duality-maximals} allows us to use the absolute maximals of $S$
instead of the relative ones in
order to describe the semigroup $S$ inductively. Moreover, if $\alpha$ is an absolute
maximal and $\alpha = \beta +\gamma$ with $\beta, \gamma\in S$, then so are $\beta$ and
$\gamma$: indeed,
if $F\in\Sigma(\beta)$ and $\sigma\in\Delta^F(\beta)$, then
$\sigma+\gamma\in\Delta^F(\alpha)$, so
$\Sigma(\beta)\subset\Sigma(\alpha)=\{\emptyset\}$.

We will say that an element $\alpha$ of $S$ is irreducible if for any
decomposition $\alpha = \beta +\gamma$ with $\beta, \gamma\in S$ one has $\beta=\0$
or $\gamma=\0$.
Thus, as the next Theorem shows (see also \cite{D-87}), one can use only the set
of irreducible absolute maximals of $S$:

\begin{theorem}\label{AbsoluteMaximals}
Let $S$ be a symmetric good semigroup.
Let $B=\{\beta^1, \ldots , \beta^k\}$ be the set of irreducible absolute maximals
of $S$, $A = \langle \sseq{\beta}{1}k \rangle_{\Zr{}} \cap [\0,\tau]$ and
$R = \tau-A = \{\tau - \alpha \st \alpha\in A\}$.
Let $\alpha\in \Zr r$ be such that $\alpha_J\in S_J$ for all
$J\subset I, J\neq I$. Then, $\alpha\notin S$ if and only if
there exist $i\in I$ and $\rho\in R$ with $\alpha_i=\rho_i$ and $\alpha_j>\rho_j$
for $j\neq i$.
\end{theorem}

\begin{proof}
The set $A$ contains all the absolute maximals of $S$, and so the set $R$ contains
all the relative maximals (Theorem~\ref{thm:duality-maximals}). Note that some elements of $R$ may not
be in $S$; however, $\Delta(\rho)=\emptyset$ for all $\rho\in R$ (see Remark~\ref{rem:maximal}).
Now, ($\Rightarrow$) is the second statement of
Theorem~\ref{generation} (if $r\ge3$ the hypothesis forces $\alpha\in\Z^r_{>0}$
unless $\alpha=\0$; if $r=2$, an element $\alpha=(0,a)\notin S$ is covered by
$\rho=\0\in R$). On the other hand, if  $\alpha\in S$ then
$\alpha\in\Delta_i(\rho)\neq\emptyset$, against $\Delta(\rho)=\emptyset$, and the Theorem is proved.
\end{proof}

\begin{remark}
If we assume that $S_J$ is also symmetric for all $J\subset I$
(in the next Section these semigroups will be called
completely symmetric) then
one can eliminate the condition $\alpha_J\in S_J$ for all $J\subset I$, $J\neq I$ in
the Theorem by using the irreducible absolute maximals of all the
projections. More precisely, for $\emptyset\neq J\subset I$ let $A(S_J)$ be the set $A$ of the
Theorem for the symmetric semigroup $S_J$, i.e.\ the set of sums of irreducible
absolute maximals of $S_J$ which are $\le\tau(S_J)$ (for $J=\{j\}$, $A(S_j)=S_j\cap [0,c_j-1]$),
and let
$R'=\bigcup_{\emptyset\neq J\subset I}\{(\tau(S_J)-a)*\0 : a\in A(S_J)\}$.
Then one has that, for $\alpha \in \Z_{>0}^r$,
$\alpha\notin S$ if and only if there exist $i\in I$ and $\rho\in R'$ with
$\alpha_i=\rho_i$ and $\alpha_j>\rho_j$ for $j\neq i$.
The proof is the same:
Theorem~\ref{generation} gives $\rho=\rho_J*\0$ with $\rho_J$ a relative maximal of
$S_J$, and $\rho_J\in\tau(S_J)-A(S_J)$ by the Theorem applied to $S_J$; conversely
$\Delta(\rho_J*\0)=\emptyset$ (Case 2c in the proof of Lemma~\ref{lem:new}).
\end{remark}

\begin{remark}\label{rem:maxcontact}
For the semigroup of values of a plane curve $C$ with branches
$C_1,\dots,C_r$, the irreducible absolute maximals have a geometric meaning.
By
\cite[Theorem (3.18)]{D-87}, $\alpha$ is an irreducible absolute maximal of $S$ if
and only if it is a (finite) maximal contact value of $C$. (For several branches
this notion was first defined by Lejeune~\cite{Lejeune}; essentially, it is the
value of an irreducible curve
whose intersection multiplicities with the branches $C_j$ of $C$ are as large as
possible for a given multiplicity.) This generalizes, to some extent, the well-known fact that,
for an irreducible curve $C$, the maximal contact values form the minimal system of
generators of the (numerical) semigroup of $C$.
\end{remark}

\medskip

Let $S\subset S_1\times \cdots \times S_r\subset \Z_{\ge 0}^r$
be a good semigroup.
Recall that $\Fin(S_i)$ denotes the set of elements of $S_i$ with finite fibre
(Notations~\ref{not1}).
Let us denote by $\chi_i$ the minimum element of $S_i$ with infinite fibre and by
$\kappa_i$ the number of elements of $S_i$ with finite fibre:
\begin{equation}\label{keynumbers}
\chi_i := \min \{\alpha_i\in S_i | F(\alpha_i) = \infty\}\; , \
\kappa_i := \# \Fin(S_i)\; .
\end{equation}

Recall that for a numerical semigroup
$\Gamma\subset  \Z_{> 0}$ and $m\in \Zr{}$, the
Ap\'ery set of $\Gamma$ with respect to $m$ is the set
$$
A_m(\Gamma) =  \{a\in \Gamma | a-m\notin \Gamma\}\; .
$$
The Ap\'ery set $A_m(\Gamma)$
is a finite set and its maximum is $c+m-1$,
$c$ being the conductor of $\Gamma$.
Moreover, one has that
$\# A_m(\Gamma) \ge m$, with equality if and only if $m\in \Gamma$.

The numbers
$\mu^i, \chi_i, \kappa_i$, $1\le i\le r$, play an important role in the description of
the semigroup $S$. The next lemma shows some relations among them.

\begin{lemma} \label{keynumbers-1}
For $i\in I$, one has
$\kappa_i \le \chi_i \quad \text{and} \quad  \mu^i\le \chi_i \;$.
\end{lemma}

\begin{proof}
If $\alpha_i\in S_i$ has finite fibre, then
$\alpha_i-\chi_i\notin S_i$  (otherwise, adding a lift of $\alpha_i-\chi_i$ to
$\chi_i*\delta\in S$ one gets $\alpha_i*\delta'\in S$ with $\delta'\ge\delta$, i.e.\
$F(\alpha_i)=\infty$)
and so $\Fin(S_i)\subset A_{\chi_i}(S_i)$. In particular,
$\alpha_i\le c_i+\chi_i-1$ and $\kappa_i\le \chi_i$.
On the other hand, if
$\mu^i>0$, then $\tau_i= c_i + \mu^i-1\in S_i$ and, as
$\Delta(\tau)=\emptyset$,
it has finite fibre ($\tau_i$ is the biggest element of $S_i$ with finite
fibre). In particular,
$\tau_i\le c_i+\chi_i-1$ and, since $\tau_i=c_i+\mu^i-1$,  one has that
$\mu^i\le \chi_i$.
\end{proof}

\begin{proposition}\label{sym-1}
Let $S$ be a symmetric semigroup. Then,
$\tau\in S$ is a relative maximal and $\mu^i=\chi_i\in S_i$ for all $i =1,\ldots,r$.
\end{proposition}

\begin{proof}
The fact that $\tau=\delta-\e$ is a relative maximal is already known
(Theorem~\ref{thm:duality-maximals} applied to $\alpha=\0$).
Note that $\tau\in S$ implies $\mu^i\ge1$ for all $i$, since
$\tau_i=c_i+\mu^i-1\in S_i$ and $c_i-1\notin S_i$.

Let us prove that
$\mu^i\in S_i$ and $F(\mu^i)=\infty$.
Assume that $\mu^i\notin S_i$. In this case
$\Delta(\mu^i*\0) = \emptyset$, thus $\beta=\tau-\mu^i*\0\in S$. However, in this
case one has that
$\beta_i= c_i-1 + \mu^i-  \mu^i=c_i-1 \in S_i$, which is a contradiction.
If $F(\mu^i)< \infty$ then, by the characterization of the elements with
infinite fibre, $\w{\Delta}_i(\mu^i*\delta)=\emptyset$ and, by Lemma
\ref{lem:new}, there exists
$\alpha\in pr_i^{-1}(\mu^i)$ such that $\Delta(\alpha)=\emptyset$.
The symmetry implies that $\tau-\alpha\in S$, but this is a contradiction
because $(\tau-\alpha)_i = c_i-1\notin S_i$.

Now, since $F(\mu^i)=\infty$, the minimality of $\chi_i$ implies $\mu^i\ge \chi_i$,
and so $\chi_i=\mu^i$.
\end{proof}

\begin{theorem}\label{sym-2}
Let us fix $i\in I$ and
assume that $S_i$ is symmetric. Then all the elements of the Ap\'ery set
$A_{\mu^i}(S_i)$ have finite fibre; as a consequence, $\mu^i \le \kappa_i\le \chi_i$.

If, moreover, $S$ is symmetric, then the elements with finite fibre are exactly
those of
$A_{\mu^i}(S_i)$, so $\mu^i=\kappa_i=\chi_i$.
\end{theorem}

\begin{proof}
If $\mu^i=0$ the result is trivial.
Let $\alpha_i\in A_{\mu^i}(S_i)$ and  assume that $F(\alpha_i)=\infty$.
First, note that
$\beta_i = c_i-1-\alpha_i+\mu^i\in S_i$, because $\alpha_i-\mu^i\notin S_i$ and $S_i$
is symmetric.
On the other hand, $F(\alpha_i)=\infty$ is equivalent to saying that
$\alpha_i * \delta\in S$, and so
$\Delta(\tau-\alpha_i*\delta)=\emptyset$. The $j$-th coordinate of
$\tau- \alpha_i*\delta$ is $-1$ if $j\neq i$ and $\beta_i$ if $j=i$ and, since
$\beta_i\in S_i$, one gets a contradiction because
$$
\emptyset = \Delta(\tau-\alpha_i*\delta) \supset \Delta_i(\beta_i * (-1)) =
F(\beta_i) \neq \emptyset\; .
$$
Thus $A_{\mu^i}(S_i)\subset \{\alpha_i\in S_i :  F(\alpha_i)<\infty\}$.
Now, $\# A_{\mu^i}(S_i)\ge \mu^i$ gives the inequality $\mu^i\le \kappa_i$ (note that
we do not assume that $\mu^i\in S_i$).

Now, let us also assume that $S$
is symmetric. Then  $\mu^i=\chi_i\in S_i$, so
$\kappa_i\le\chi_i=\mu^i=\#A_{\mu^i}(S_i)\le\kappa_i$ and the inclusion
$\Fin(S_i)\subset A_{\chi_i}(S_i)$ of Lemma~\ref{keynumbers-1} is an equality.
\end{proof}

\medskip

Note that the first statement of the Theorem holds for any good semigroup $S$
(only the symmetry of $S_i$ is used).
The next example shows
that the hypothesis ``$S_i$ symmetric'' cannot be dropped in either statement of
Theorem~\ref{sym-2}.

\begin{example}\label{ex4}
Let
$$
S=\{\0,(3,3),(4,4),(7,7)\}\cup\{(5,5+k),(5+k,5),(6,6+k),(6+k,6):k\ge0\}\cup\big((8,8)+\Zr2\big)\; .
$$
One checks that $S$ is a symmetric good semigroup (see Figure~\ref{fig2})
with $\delta=(8,8)$,
$\tau=(7,7)$, $P_S=1+\tt^{(3,3)}+\tt^{(4,4)}+\tt^{(7,7)}$ and
$S_1=S_2=\langle3,4,5\rangle$, which is not symmetric. Here $c_i=3$,
$\mu^i=\chi_i=5$ (as predicted by Proposition~\ref{sym-1}), but
$\Fin(S_i)=\{0,3,4,7\}$, so $\kappa_i=4<\mu^i$ and
$\Fin(S_i)\subsetneq A_{\mu^i}(S_i)=\{0,3,4,6,7\}$. Note that
$6\in A_{\mu^i}(S_i)$ has infinite fibre, so the first statement of
Theorem~\ref{sym-2} fails as well.

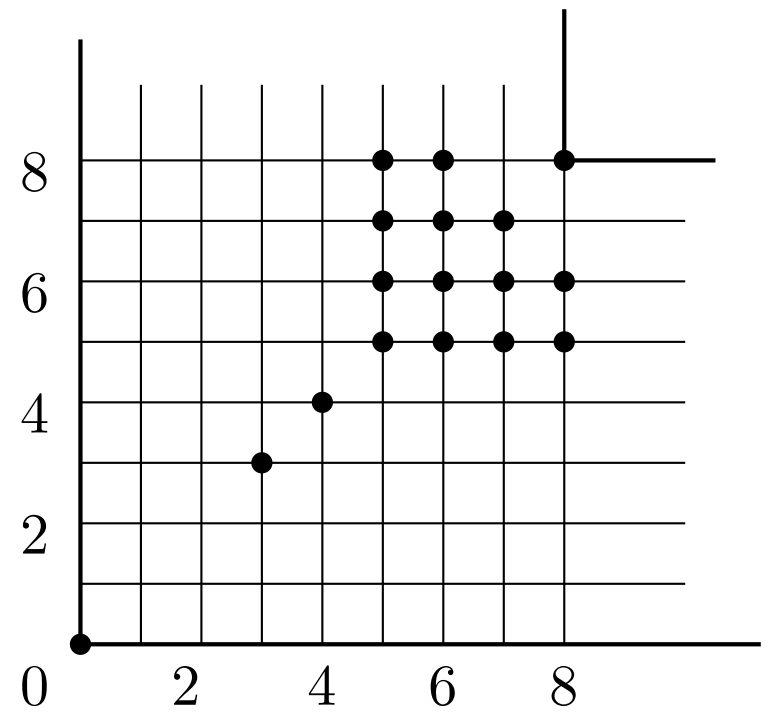
\begin{figure}[h]
\centering
\unitlength=1.00mm
\begin{picture}(50,50)(8,8)
\thicklines
\put(10,10){\line(1,0){45}}
\put(10,10){\line(0,1){40}}
\put(6,6){$0$}
\put(16,6){$2$}
\put(25,6){$4$}
\put(33,6){$6$}
\put(41,6){$8$}
\put(6,16){$2$}
\put(6,24){$4$}
\put(6,32){$6$}
\put(6,40){$8$}

\thinlines
\put(14,10){\line(0,1){37}}
\put(18,10){\line(0,1){37}}
\put(22,10){\line(0,1){37}}
\put(26,10){\line(0,1){37}}
\put(30,10){\line(0,1){37}}
\put(34,10){\line(0,1){37}}
\put(38,10){\line(0,1){37}}
\put(42,10){\line(0,1){37}}

\put(10,14){\line(1,0){40}}
\put(10,18){\line(1,0){40}}
\put(10,22){\line(1,0){40}}
\put(10,26){\line(1,0){40}}
\put(10,30){\line(1,0){40}}
\put(10,34){\line(1,0){40}}
\put(10,38){\line(1,0){40}}
\put(10,42){\line(1,0){40}}

\put(10,10){\circle*{1.5}}
\put(22,22){\circle*{1.5}}
\put(26,26){\circle*{1.5}}
\put(30,30){\circle*{1.5}}
\put(30,34){\circle*{1.5}}
\put(30,38){\circle*{1.5}}
\put(30,42){\circle*{1.5}}
\put(34,30){\circle*{1.5}}
\put(38,30){\circle*{1.5}}
\put(42,30){\circle*{1.5}}

\put(34,34){\circle*{1.5}}
\put(34,38){\circle*{1.5}}
\put(34,42){\circle*{1.5}}
\put(38,34){\circle*{1.5}}
\put(42,34){\circle*{1.5}}

\put(38,38){\circle*{1.5}}

\put(42,42){\circle*{1.5}}

\thicklines
\put(42,42){\line(1,0){10}}
\put(42,42){\line(0,1){10}}
\end{picture}
\caption{The semigroup $S$ of Example~\ref{ex4}}
\label{fig2}
\end{figure}

\end{example}

\section{Beyond symmetry}\label{Supersymmetric}

In this section we consider good semigroups all of whose
projections are symmetric; we call them \emph{completely symmetric}. The main
examples are the value semigroups of plane curves: any union of branches of a plane
curve is again a plane curve, hence Gorenstein, so by \cite{D-88} all the
projections of its value semigroup are symmetric. For such semigroups we prove
two main results. First, the conductor is computed from the projections
(Theorem~\ref{thm:supersymmetry}): $\delta_i=c_i+\sum_{j\neq i}\xi_{ij}$, where
$\xi_{ij}$ is the number of maximal elements of $S_{\{i,j\}}$. Moreover, for
completely symmetric semigroups the three numbers attached to each index $i$ in
Section~\ref{sec:Symmetry} (the difference $\mu^i=\delta_i-c_i$, the minimal element
$\chi_i$ of $S_i$ with infinite fibre and the number $\kappa_i$ of elements with
finite fibre) coincide, and are equal to $\sum_{j\neq i}\xi_{ij}$. More generally,
the elements of each projection $S_J$ with finite fibre are described by an Ap\'ery
set (Proposition~\ref{prop:Apery-fibres}). This means that the numbers $\xi_{ij}$
are the combinatorial analogue of the intersection multiplicities of two plane
branches.
Second, we prove a projection formula
for the Poincar\'e series (Theorem~\ref{projection} and
Corollary~\ref{cor:projection1}), the analogue of the formula of Torres for the
Alexander polynomial, which holds as soon as $S$ and $S^{i}$ are symmetric, and fails in
general for symmetric $S$ (Example~\ref{ex:two-branch} and Proposition~\ref{thm:two}).
We also compare completely symmetric semigroups with the semigroups of plane curves
(Subsection~\ref{sec:CSplane}): the former class is strictly larger, but the two
classes coincide when all the projections $S_i$ are equal to the semigroup $M$ of a
plane branch and the two-dimensional projections are ``diagonal''
(see Definition~\ref{def:diagonal} and Theorem~\ref{thm:diagonal-plane}). This is the
case, for instance, when $S$ is completely symmetric and all the projections $S_i$ are
equal to $\Zr{}$ (Corollary~\ref{cor:smooth-plane}).

\subsection{Completely symmetric Semigroups}

\begin{definition}\label{def:completely}
Let $S\subset \Zr r$ be a good semigroup. We will say that $S$ is \emph{completely
symmetric} if
$S_J$ is symmetric for all $J\subset I$, $J\neq\emptyset$.
\end{definition}

\begin{remark}\label{rmk r=2}
Let us fix $r=2$ and suppose that $S_1, S_2\subset \Zr{}$ are numerical
semigroups. Assume that $S\subset \Zr 2$ is a good semigroup such that $S_1$ and $S_2$
are its projections.
There is a natural
bijection between maximals of $S$ and $\Fin(S_1)$ (the elements of $S_1$ with finite
fibre):
if $a\in\Fin(S_1)$ the maximum $(a,b)$ of its fibre is a maximal element (an element
$(a',b)\in S$ with $a'>a$ would give, by P3 at the second index, an element $(a,b')$
with $b'>b$), and conversely if $(a,b)$ is maximal then $a\in\Fin(S_1)$ and $b$ is the
maximum of the fibre. Obviously the same is true for $S_2$ and $\Fin(S_2)$.

Let us suppose now that $S_1$ and $S_2$ are symmetric and
let $\xi = \# \Fin(S_i)$, $i=1,2$, be the number of maximals of $S$.
So, $\xi=\kappa_1=\kappa_2$. By Theorem \ref{sym-2}, if $S$ is symmetric then
$\xi = \kappa_i = \mu^i$, $i=1,2$.

Let us see that in fact
$$
S \text{ is symmetric } \iff \mu^1 = \mu^2 =\xi \; .
$$

Indeed, let $\0 = \alpha^0 < \cdots < \alpha^{\abs{\delta}}$ be  the path from
$\0 =  (0,0)$ to $\delta= (c_1+\mu^1, c_2+\mu^2)$ going first horizontally from $\0$
to $(\delta_1,0)$ and then vertically up to $\delta$.
By using equation (\ref{path}) and Remark \ref{numerical} one sees that
$$
\ell(\delta) = \frac{c_1}2 + \mu^1 + \frac{c_2}2 + \mu^2 -\xi =
\frac{\abs{\delta}}2 + \frac{\mu^1+\mu^2 -2\xi}{2} \; .
$$
By Proposition~\ref{symmetric},
$S$ is symmetric if and only if $\mu^1+\mu^2=2\xi$ and, since $\mu^i\le\kappa_i=\xi$
by
the first part of Theorem~\ref{sym-2} (which only uses the symmetry of $S_i$), this
happens if and only if $\mu^1=\mu^2=\xi$.

Note that, if $S$ is symmetric, then $\xi\in S_1\cap S_2$, the conductor of $S$ is
$(c_1+\xi, c_2+\xi)$ and $\xi$ is the minimum element of $S_1$ (or of $S_2$) with
infinite fibre in $S$.
Thus, the number $\xi$ plays the role of the intersection multiplicity in the
classical setting of semigroups of plane branches.

For the semigroup $T$ of
Example~\ref{ex-symmetry}(2) one has $T_1=T_2=\langle2,5\rangle$ (symmetric),
$\xi=3$ (the maximals are $\0,(2,2),(4,4)$), $c_i=4$ and $\delta=(5,5)$, so
$\mu^i=1<\xi$ and $T$ is not symmetric, as we already know. Note that
$\chi_i=5>\xi$.
\end{remark}

\medskip

If $r\ge 2$, we will denote by
$\xi_{ij}$ the number of maximals of the semigroup $S_{\{i,j\}}$.
Note that, if $S_{\{i,j\}}$ is completely symmetric, then
the conductor of $S_{\{i,j\}}$ is $\delta(S_{\{i,j\}}) = (c_i+\xi_{ij},
c_j+\xi_{ij})$.

\begin{theorem}\label{thm:supersymmetry}
Let $S\subset \Zr r$, $r\ge 2$,  be a good semigroup such that $S_J$ is symmetric for
all
$J\subset I$, $J\neq I$. Then, $S$ is symmetric if and only if
$\mu^i = \sum_{j=1, j\neq i}^r \xi_{ij} = \kappa_i$ for all $i=1,\ldots, r$.
\end{theorem}

\begin{proof}
The proof is by induction on $r$. The case $r=2$ has already been proved in Remark
\ref{rmk r=2} (for $r=2$ one has $\kappa_i=\xi_{12}$ for $i=1,2$).
Let $\delta = (c_1+\mu^1, \ldots, c_r+\mu^r)\in S$ be the conductor of $S$
and denote $\xi^i = \sum_{j\neq i}\xi_{ij}$, $i=1,\ldots, r$.
Let us fix the index $i=1$ and denote
$S^1 = S_{I\setminus\{1\}}$, $\delta^1 = \delta(S^1)$,
$\eta = (\eta^2,\ldots, \eta^r) = pr_{I\setminus\{1\}}(\delta) - \delta^1$.
Let us first assume only that $S$ is a good semigroup and that $S_1$ is symmetric.
In this case one can compute
$\ell(\delta)$ going through $\delta^1*0$ in the  following way:

\begin{align*}
\ell(\delta) & = \ell ((\delta^1+\eta)*0) +
\sum_{k=0}^{c_1+\mu^1-1} \ell((\delta^1+\eta)*0 + k \e_1, 1) \\
& = \ell_{S^1}(\delta^1) + \sum_{j=2}^r \eta^j +
\sum_{k=0}^{c_1+\mu^1-1} \ell((\delta^1+\eta)*0 + k \e_1, 1)\; .
\end{align*}
The last equality holds because $\ell(\delta^1*0)=\ell_{S^1}(\delta^1)$ (the latter
computed in $S^1$) and each step from $(0,\delta^1)$ to
$(0,\delta^1+\eta)$ contributes $1$, $\delta^1$ being the conductor of $S^1$.

Now, one has that
$\ell((\delta^1+\eta)*0 + k \e_1, 1) = 1$ if and only if
$\w{\Delta}_1((\delta^1+\eta)*0 + k \e_1) \neq \emptyset$, i.e.\
$\ell((\delta^1+\eta)*0 + k \e_1, 1) = 1$ if and only if
$F(k) =\infty$ (for $k\in S_1$). The number of elements $n$ of $S_1$ with
$0\le n< c_1+\mu^1$ and
$F(n)=\infty$ coincides with the number of elements of
$S_1\cap [0,c_1+\mu^1)$ minus $\kappa_1$,
the number of elements of $S_1$ with finite fibre; thus the last summand of the above
formula is equal to $\#(S_1\cap [0,c_1+\mu^1))-\kappa_1$.
Since $S_1$ is symmetric,
$\#(S_1\cap [0,c_1+\mu^1)) = c_1/2 +\mu^1$, and as a consequence
\begin{equation}\label{eq1}
\ell(\delta)
= \ell_{S^1}(\delta^1) + \sum_{j=2}^r \eta^j +
c_1/2 + \mu^1 -\kappa_1\;
= \ell_{S^1}(\delta^1) + \abs{\eta} +
c_1/2 + \mu^1 -\kappa_1\; .
\end{equation}

Now, let us assume that $S^1$ is completely symmetric; in particular, $S^1$ is symmetric.
By the induction hypothesis,
the $i$-th coordinate of $\delta^1$ is
$$
\delta^1_i = c_i + \sum_{j=2, j\neq i}^r \xi_{ij} = c_i + \xi^i -\xi_{1i}\; .
$$
Thus,
$$
\abs{\delta^1} = \sum_{i=2}^r (c_i + \xi^i- \xi_{1i}) =
\sum_{i=2}^r ( c_i + \xi^i)  - \xi^1\;
$$
and
$$
\abs{\eta} = \sum_{i=2}^r ( c_i + \mu^i) - \sum_{i=2}^r (c_i + \xi^i) + \xi^1
= \sum_{i=2}^r (\mu^i -\xi^i)  + \xi^1 \; .
$$

The semigroup $S^1$ is symmetric, so $\ell_{S^1}(\delta^1)= (1/2)\abs{\delta^1}$ and
substituting the above computations of $\abs{\eta}$ and $\abs{\delta^1}$ in
equation~(\ref{eq1}) and
rearranging the sums one has that:
\begin{align*}
\ell(\delta)
&= (1/2)\left( \sum_{i=2}^r ( c_i + \xi^i)  - \xi^1 \right) +
\sum_{i=2}^r (\mu^i -\xi^i)  + \xi^1 + c_1/2 +\mu^1 -\kappa_1 \\
&= (1/2) \sum_{i=1}^r  c_i - (1/2)\sum_{i=2}^r \xi^i + \sum_{i=1}^r \mu^i
-(1/2)\xi^1 + \xi^1 -\kappa_1   \\
&= (1/2) \sum_{i=1}^r  (c_i + \mu^i) + (1/2) \sum_{i=1}^r \mu^i
- (1/2)\sum_{i=1}^r \xi^i + \xi^1 -\kappa_1\\
&= (1/2)\abs{\delta} + (1/2) \sum_{i=1}^r (\mu^i - \xi^i)  + \xi^1
-\kappa_1
\end{align*}

Assume that $S_J$ is symmetric for all $J\subset I$, $J\neq I$. Then, for each
$j=1,\ldots, r$ we have the same equation with $1$ replaced by $j$, i.e.,
\begin{equation}\label{eq2}
\ell(\delta) = \frac 12 \abs{\delta} + \frac 12 \sum_{i=1}^r (\mu^i-\xi^i) + \xi^j
-\kappa_j \qquad \qquad j=1, \ldots, r \; .
\end{equation}

If $\xi^j=\mu^j=\kappa_j$ for all $j=1,\ldots, r$, then equation~(\ref{eq2}) gives
$\ell(\delta) = (1/2)\abs{\delta}$ and,  by Proposition \ref{symmetric},   $S$
is a symmetric semigroup.

Conversely, let us assume that $S$ is symmetric, that is,
$\ell(\delta)=(1/2)\abs{\delta}$. By
Theorem \ref{sym-2} we know that $\kappa_j=\mu^j$ for $j\in I$, so equation
(\ref{eq2}) now reads
$$
0 =  \frac 12 \sum_{i=1}^r (\mu^i-\xi^i) + \xi^1 -\mu^1 =
\frac 12  \left( \sum_{i=2}^r \mu^i -\mu^1 - \sum_{i=2}^r \xi^i + \xi^1 \right) \; .
$$
Thus we have the following linear system of equations relating
$(\mu^1,\ldots, \mu^r)$ and $(\xi^1,\ldots, \xi^r)$:
\begin{equation}\label{eq3}
\sum_{i=1, i\neq j}^r \mu^i -\mu^j = \sum_{i=1, i\neq j}^r \xi^i -  \xi^j \; ,
\text{ for } j=1,\ldots, r \; .
\end{equation}
Note that the vector $(\xi^1,\ldots, \xi^r)$, i.e.\ the right hand side of the system
(\ref{eq3}), depends only on the projections $S_{\{i,j\}}$, so one can see
the system as a linear system in the unknowns $(\mu^1,\ldots, \mu^r)$.
The matrix of the system has
eigenvalues $r-2$ and $-2$, so it is regular exactly for $r\ge3$ (which is why the
case
$r=2$ has to be treated separately in Remark~\ref{rmk r=2}).
So,
the unique solution is $\mu^j = \xi^j$ for $1\le j \le r$, and the proof is finished.
\end{proof}

\begin{notations}\label{not3}
Let $\beta\in S$. The Ap\'ery set of $S$ with respect to $\beta$ is the set
$A_{\beta}(S)$ of elements
$\alpha\in S$ such that  $\alpha-\beta\notin S$.

For a nonempty subset $K\subseteq I$ we denote by $\delta(S_K)$ the conductor of
the good semigroup $S_K=pr_K(S)\subset\mathbb{Z}^{\#K}_{\ge 0}$ and by
$\tau(S_K)=\delta(S_K)-pr_K(\e)$ its Frobenius vector; $\Delta^{S_K}_i(-)$ and
$\Delta^{S_K}(-)=\bigcup_{i\in K}\Delta_i^{S_K}(-)$ denote the sets of
Notations~\ref{not2} computed in the semigroup $S_K$. Note that $(S_K)_i=S_i$ for
$i\in K$, so the numbers $c_i$ do not depend on $K$, and that
$(S_K)_{\{i,j\}}=S_{\{i,j\}}$, so the numbers $\xi_{ij}$ do not depend on $K$
either. Since $pr_K(\delta+\mathbb{Z}^r_{\ge 0})=pr_K(\delta)+\mathbb{Z}^{\#K}_{\ge 0}
\subset S_K$, the minimality of the conductor (Remark~\ref{rem:conductor}) gives
\begin{equation}\label{eq:deltaK-le}
\delta(S_K)\le pr_K(\delta).
\end{equation}
If $A\subset K$, we say that $\alpha_A\in S_A$ \emph{has infinite fibre in}
$S_K$ if it satisfies the equivalent conditions 1)--4) of Notations~\ref{not1} with respect
to the good semigroup $S_K$; for instance, this is the case if there exists
$\beta\in S_K$ with $pr_A(\beta)=\alpha_A$ and $\beta_k\ge\delta(S_K)_k$ for all
$k\in K\setminus A$. For $K=I$ this is the notion $F(\alpha_A)=\infty$ used so far.
For $i\in K$ we put
\[
\chi_i(S_K)=\min\{a\in S_i \;:\; a \text{ has infinite fibre in } S_K\},
\]
so that $\chi_i=\chi_i(S)$.
\end{notations}

Let $S$ be a completely symmetric semigroup and $K\subseteq I$, $i\in K$,  with
$\#K\ge 2$. The projection $S_K$ is also completely symmetric. Thus, by
Theorem~\ref{thm:supersymmetry} and Proposition~\ref{sym-1} applied to $S_K$ one has
\[
\chi_i(S_K)=\mu^i(S_K)=\sum_{j\in K\setminus\{i\}}\xi_{ij},
\qquad\text{hence}\qquad
\delta(S_K)_i=c_i+\sum_{j\in K\setminus\{i\}}\xi_{ij}.
\]
In particular, for every $J$ with $\emptyset\neq J\subsetneq I$,
\begin{equation}\label{eq:prJ-delta}
pr_J(\delta)=\delta(S_J)+\xi^J,
\qquad\text{where}\quad
(\xi^J)_i=\sum_{j\in I\setminus J}\xi_{ij}\quad(i\in J).
\end{equation}

\begin{proposition}\label{prop:Apery-fibres}
Let $S\subset\mathbb{Z}^r_{\ge 0}$, $r\ge 2$, be a completely symmetric semigroup,
let $J$ be a subset with $\emptyset\neq J\subsetneq I$, and let
$\xi^J\in\mathbb{Z}^{\#J}_{\ge 0}$ be defined by
$(\xi^J)_i=\sum_{j\notin J}\xi_{ij}$ for $i\in J$. Then:
\begin{enumerate}
\item[(i)]
Every $\alpha_J\in S_J$ with infinite fibre satisfies
$\alpha_J\ge\xi^J$.
\item[(ii)]
$\xi^J\in S_J$ and $\xi^J$ has infinite fibre, i.e.\ $F(\xi^J)=\infty$.
As a consequence, $\xi^J$ is the minimum (for the product order) of the
set of elements of $S_J$ with infinite fibre.
\item[(iii)] For $\alpha_J\in S_J$ one has $F(\alpha_J)<\infty$ if and only if
$\alpha_J\in A_{\xi^J}(S_J)$, that is,
\[
\Fin(S_J)=A_{\xi^J}(S_J)=\{\alpha_J\in S_J \;:\; \alpha_J-\xi^J\notin S_J\}.
\]
\end{enumerate}
\end{proposition}

\begin{proof}
For $i\in J$, write $K_i=\{i\}\cup (I\setminus J)$. Since
$\#K_i\ge 2$, one has
\begin{equation}\label{eq:chi-Ki}
\chi_i(S_{K_i})=\sum_{j\in J^c}\xi_{ij}=(\xi^J)_i\qquad(i\in J).
\end{equation}

\emph{Step 1.} We claim: \emph{if $\beta\in S$ satisfies $\beta_k\ge\delta_k$ for
all $k\notin J$, then $pr_J(\beta)\ge\xi^J$.} Let $i\in J$. The element
$pr_{K_i}(\beta)\in S_{K_i}$ has $i$-th coordinate $\beta_i$ and, by
\eqref{eq:deltaK-le}, $k$-th coordinate $\beta_k\ge\delta_k\ge\delta(S_{K_i})_k$ for
every $k\notin J$. Hence $\beta_i$ has infinite fibre in $S_{K_i}$ and, by
\eqref{eq:chi-Ki}, $\beta_i\ge\chi_i(S_{K_i})=(\xi^J)_i$.

If $\alpha_J\in S_J$ has infinite fibre, condition~2) in Notations~\ref{not1} provides
$\beta\in S$ with
$pr_J(\beta)=\alpha_J$ and $\beta_k\ge\delta_k$ for $k\notin J$; by the claim,
$\alpha_J\ge\xi^J$. This proves (i).

\medskip
\emph{Step 2: proof of (ii).} Let $i\in J$. By \eqref{eq:chi-Ki}, $(\xi^J)_i$ has
infinite fibre in $S_{K_i}$. By condition~3) in Notations~\ref{not1} applied to
$S_{K_i}$ and
\eqref{eq:deltaK-le}, the element of $\mathbb{Z}^{\#K_i}$ with $i$-th coordinate
$(\xi^J)_i$ and $k$-th coordinate $\delta_k$ for $k\notin J$ belongs to $S_{K_i}$.
Let $\gamma^i\in S$ be such that $pr_{K_i}(\gamma^i)$ is this element, i.e.\
$\gamma^i_i=(\xi^J)_i$ and $\gamma^i_k=\delta_k$ for $k\notin J$. By the claim of
Step~1,
$\gamma^i_j\ge(\xi^J)_j$ for every $j\in J$. Therefore, by property P2,
\[
\gamma=\bigwedge_{i\in J}\gamma^i = \inf \{\gamma^i : \ i\in J\}\in S
\]
satisfies $\gamma_i=\gamma^i_i=(\xi^J)_i$ for $i\in J$ and $\gamma_k=\delta_k$ for
$k\notin J$,
that is, $\gamma=\xi^J*\delta\in S$. Thus $\xi^J\in S_J$ and, by condition~1) in
Notations~\ref{not1},
$\xi^J$ has infinite fibre.

\medskip
\emph{Step 3: $\Fin(S_J)\subset A_{\xi^J}(S_J)$.} Let $\alpha_J\in S_J$ with
$\alpha_J-\xi^J\in S_J$, and let $\rho\in S$ with $pr_J(\rho)=\alpha_J-\xi^J$. By
Step~2, $\rho+\xi^J*\delta\in S$; this element projects onto $\alpha_J$ and its
$k$-th coordinate is $\ge\delta_k$ for every $k\notin J$. By condition~2) in
Notations~\ref{not1},
$\alpha_J$ has infinite fibre. Hence, if $F(\alpha_J)<\infty$, then
$\alpha_J-\xi^J\notin S_J$, i.e.\ $\alpha_J\in A_{\xi^J}(S_J)$.

\medskip
\emph{Step 4: $A_{\xi^J}(S_J)\subset\Fin(S_J)$.} Let $\alpha_J\in S_J$ with
infinite fibre; we show that $\alpha_J-\xi^J\in S_J$. By condition~1) in
Notations~\ref{not1},
$\alpha_J*\delta\in S$, which gives
$\Delta(\gamma)=\emptyset$ for $\gamma=\tau-\alpha_J*\delta$ (see the remark after
Definition~\ref{def:symmetric}). Here
$\gamma_k=\tau_k-\delta_k=-1$ for $k\notin J$ and
$pr_J(\gamma)=pr_J(\tau)-\alpha_J$.
Let $i\in J$. Every $\sigma\in S$ satisfies $\sigma_k>-1$ for $k\notin J$, hence
\[
\Delta_i(\gamma)=\{\sigma\in S \;:\; pr_J(\sigma)\in\Delta_i^{S_J}(pr_J(\gamma))\} =
\emptyset
\quad\text{and}\quad
\Delta_i^{S_J}(pr_J(\gamma))=pr_J\big(\Delta_i(\gamma)\big)=\emptyset .
\]
Thus $\Delta^{S_J}(pr_J(\gamma))=\emptyset$. On the other hand, by
\eqref{eq:prJ-delta},
\[
pr_J(\tau)=pr_J(\delta)-pr_J(\e)=\delta(S_J)+\xi^J-pr_J(\e)=\tau(S_J)+\xi^J,
\]
so that $pr_J(\gamma)=\tau(S_J)-(\alpha_J-\xi^J)$. Since $S_J$ is symmetric,
Definition~\ref{def:symmetric} applied to $S_J$ yields $\alpha_J-\xi^J\in S_J$, that is,
$\alpha_J\notin A_{\xi^J}(S_J)$.
(If $\#J=1$, $S_J$ is a symmetric numerical
semigroup and the symmetry is the usual one: $a\in S_J$ if and only if
$\tau(S_J)-a\notin S_J$, see Remark~\ref{numerical}.)

Steps~3 and~4 prove (iii).
\end{proof}

\begin{remark}
For $J=\{i\}$ one has $\xi^{\{i\}}=\sum_{j\neq i}\xi_{ij}=\mu^i$ by
Theorem~\ref{thm:supersymmetry}, and Proposition~\ref{prop:Apery-fibres} gives
$\chi_i=\mu^i$ and
$\Fin(S_i)=A_{\mu^i}(S_i)$, which is the second part of Theorem~\ref{sym-2}.

The Ap\'ery set $A_\beta(S)$ of the value semigroup of a curve with
respect to an element $\beta\in S$ was considered in \cite{D-87}. For $r\ge 2$ it is an
infinite set; however, it is controlled by a finite number of special elements
(in a similar way
to the maximals), called the \emph{Ap\'ery corners}, which are used to describe the
semigroup inductively. For arbitrary good semigroups it has been studied by
D'Anna, Guerrieri and Micale \cite{DGM}, who show in particular that, for $r\ge2$,
$A_\beta(S)$ can be partitioned into a finite number of ``levels'' (see also
\cite{GMM-21,GMM-23}), a partition which, in the case of curves, is closely related to
the ring.

Proposition~\ref{prop:Apery-fibres}
shows that, for completely symmetric semigroups, the finiteness of the fibres of
the projection $S\to S_J$ is governed by the Ap\'ery set of $S_J$ with respect to
the single element $\xi^J$. For $J=\{i\}$ this is the Ap\'ery set of the numerical
semigroup $S_i$ with respect to $\mu^i$, which has exactly $\mu^i$ elements. For
the value semigroup of a plane curve this recovers the description of the finite
fibres given in \cite{D-87} in terms of the intersection multiplicities
$(\xi^J)_i=\sum_{j\notin J}(C_i\cdot C_j)$.
\end{remark}

\subsection{Completely symmetric versus plane semigroups}\label{sec:CSplane}

The previous results (and those of Subsection~\ref{sec:projection} and
Section~\ref{sec:recovering}) show that completely
symmetric semigroups behave, from many points of view, like the semigroups of plane
curves. It is therefore natural to ask how far the former are from the latter.
The following example of a space curve shows that the class of completely
symmetric semigroups is strictly larger than that of semigroups of plane curves.

\begin{example}\label{ex:spacecusps}
Let $S\subset\Zr2$ be the semigroup of values of the space curve with branches
$C_1:\ (t^2,t^3,0)$ and $C_2:\ (0,s^3,s^2)$. One finds
$$
S=\{(0,0)\}\cup\{(2,k),(k,2): k\ge2\}\cup\{(3,3)\}\cup\big((4,4)+\Zr2\big)\; ,
$$
with $S_1=S_2=\langle2,3\rangle$, $\delta=(4,4)$, $\tau=(3,3)$ and maximal elements
$\0$ and $(3,3)$, so that $\xi_{12}=2$ and $P_S(\tt)=1+t_1^3t_2^3$. The semigroup $S$
is completely symmetric (in agreement with Remark~\ref{rmk r=2}:
$\mu^1=\mu^2=\xi_{12}=2$), and in particular the curve is Gorenstein. However, $S$ is
not the semigroup of values of a plane curve: two plane branches of multiplicities
$n_1,n_2$ have intersection multiplicity at least $n_1n_2$, while here
$n_1=n_2=2$ and $\xi_{12}=2<4$. Note also that $S$ and the non-symmetric semigroup
of Example~\ref{ex-symmetry}(1) have the same Poincar\'e series (see
Remark~\ref{rem:PS-sym-not-S}).
\end{example}

On the other hand, we show now that complete symmetry characterizes the semigroups
of plane curves all of whose branches have the same semigroup and whose pairwise
contacts are large enough, in the sense that all the two-dimensional projections
are \emph{diagonal} (Theorem~\ref{thm:diagonal-plane}). This includes the case in
which all the projections $S_i$ are equal to $\Zr{}$, i.e.\ for curves, the case in
which all the branches are smooth (Corollary~\ref{cor:smooth-plane}). The rest of
the subsection is devoted to the proof of these facts.

\medskip
\noindent{\bf Setting.}
In this subsection, for simplicity, we will assume that the plane curves are defined
over $\C$. Given the results in \cite{Cam}, there is no difference in replacing $\C$
with an
algebraically closed field of arbitrary characteristic.
A direct treatment of constructions similar to those below, using
Hamburger--Noether expansions (as presented in \cite{Cam}) instead of Puiseux
expansions, can be found in \cite[Section~5]{Rendiconti}.

In the rest of this subsection $M\subset\Zr{}$ is the
semigroup of values of a plane branch, $m_0<m_1<\dots<m_g$ is its minimal system of
generators (the maximal contact values $\bar\beta_0,\dots,\bar\beta_g$ of
\cite[(1.1.2)]{D-94}) and $c$ is its conductor; $g=0$ means $M=\Zr{}$ (or $m_0=1$),
i.e.\ the
branch is smooth. For $0\le i\le g$ we put $e_i=\gcd(m_0,\dots,m_i)$ and, for $1\le
i\le g$, $N_i=e_{i-1}/e_i$. Recall (see \cite[(1.1.2)]{D-94}) that
$N_im_i\in\langle m_0,\dots,m_{i-1}\rangle$ for $1\le i\le g$, that $N_im_i<m_{i+1}$
for $1\le i<g$, and that $c-1=\sum_{i=1}^g(N_i-1)m_i-m_0$. We put
\[
\nu=\nu(M)=e_{g-1}m_g=N_gm_g\quad\text{if } g\ge1,\qquad\qquad \nu(\Zr{})=1 .
\]
One has that
$c<\nu$, $m_i<\nu$ for $i\le g$ if $g\ge 1$, and every integer $\ge\nu$
belongs to $M$.
In terms of Puiseux expansions, if $C_1$ and $C_2$ are two plane branches with
semigroup $M$, then $(C_1\cdot C_2)\ge\nu$ if and only if (suitable conjugates of)
their expansions coincide in all the terms of degree smaller than the last
characteristic exponent.

\begin{definition}\label{def:diagonal}
For an integer $\xi\ge\nu(M)$ we put
\[
C(M,\xi):=\{(a,a)\;:\;a\in M\}\ \cup\ \{(a,b)\in M\times M\;:\;a\neq b,\
\min\{a,b\}-\xi\in M\}\subset\Zr 2 .
\]
For $M=\Zr{}$ we simply write
$C(\xi):=C(\Zr{},\xi)=\{(a,a)\;:\;0\le a<\xi\}\cup\bigl((\xi,\xi)+\Zr 2\bigr)$; in this case
$P_{C(\xi)}=1+t_1t_2+\dots+(t_1t_2)^{\xi-1}$.
We will say that $C(M,\xi)$ is a \emph{diagonal} semigroup; by
Lemma~\ref{lem:diagonal} it is the semigroup of values of the diagonal curves with two
branches of \cite[(3.1.1)]{D-94}.
\end{definition}

We will say that the arrangement of positive integers $(\xi_{ij} : i,j\in I; i\neq
j)$ is \emph{ultrametric} if it is symmetric (i.e.\ $\xi_{ij}=\xi_{ji}$) and
$$
\xi_{ij}\ge\min\{\xi_{ik},\xi_{jk}\}
$$
for pairwise distinct $i,j,k$ (no condition
if
$r=2$).  Equivalently, the minimum of $\xi_{ij},\xi_{ik},\xi_{jk}$ is attained at least twice.
We first show that the plane curves needed in the sequel exist.

\begin{lemma}\label{lem:realization}
Let $(\xi_{ij} : i,j\in I; i\neq j )$ be ultrametric with $\xi_{ij} \ge \nu(M)$ for every
$i\neq j$.
Then there exists a plane curve $C=C_1\cup\dots\cup C_r$ whose branches have
semigroup $M$ and such that $(C_i\cdot C_j)=\xi_{ij}$ for $i\neq j$.
\end{lemma}

\begin{proof}
Put $d_{ij}=\xi_{ij}-\nu\ge0$ and $d=\max d_{ij}$. Since $(\xi_{ij})$ is
ultrametric, for each $s\ge0$ the relation ``$i=j$ or $d_{ij}>s$'' is an equivalence
relation on $I$. Choose complex numbers $a_{i,s}$, $i\in I$, $0\le s\le d$, such that
$a_{i,s}=a_{j,s}$ if and only if $i=j$ or $d_{ij}>s$, and such that the numbers
$a_{i,0}$ are positive real numbers. In particular, for $i\neq j$ one has
$a_{i,s}=a_{j,s}$ for $s<d_{ij}$ and $a_{i,d_{ij}}\neq a_{j,d_{ij}}$.

If $g=0$, let $C_i$ be the smooth branch $y=\sum_{s=0}^d a_{i,s}x^{1+s}$; then
$(C_i\cdot C_j)=\ord_x\big(\sum_s(a_{i,s}-a_{j,s})x^{1+s}\big)=1+d_{ij}=\xi_{ij}$.

Let $g\ge1$, $n=m_0$, and let $b_1<\dots<b_g$ be the characteristic exponents of a
branch with semigroup $M$ (the integers $\beta_1,\dots,\beta_g$ of
\cite[(1.1.2)]{D-94}), so that $e_i=\gcd(n,b_1,\dots,b_i)$, $m_1=b_1$ and
$m_{i+1}=N_im_i+b_{i+1}-b_i$ for $1\le i<g$. Multiplying the last equality by $e_i$
and adding up one gets
\begin{equation}\label{eq:nu-puiseux}
\nu=e_{g-1}m_g=\sum_{l=1}^{g-1}(e_{l-1}-e_l)\,b_l+e_{g-1}b_g .
\end{equation}
Let $C_i$ be the branch with parametrization $x=u^n$, $y=y_i(u)$, where
\[
y_i(u)=\sum_{l=1}^{g-1}u^{b_l}+\sum_{s=0}^{d}a_{i,s}\,u^{b_g+s}.
\]
Since $a_{i,0}\neq0$ and $\gcd(n,b_1,\dots,b_g)=1$, the characteristic exponents of
$C_i$ are $b_1,\dots,b_g$, so $C_i$ has semigroup $M$. For $i\neq j$ one has (see
e.g.\ \cite{Cam} or \cite{D-94})
\[
(C_i\cdot C_j)=\sum_{\omega^n=1}\ord_u\big(y_i(u)-y_j(\omega u)\big).
\]
Let $\omega$ be an $n$-th root of unity. If $\omega^{b_l}\neq1$ for some $l\le g-1$
and $l$ is the smallest such index, the terms of degree $<b_l$ of $y_i(u)$ and
$y_j(\omega u)$ coincide and the coefficient of $u^{b_l}$ in the difference is
$1-\omega^{b_l}\neq0$, so the order is $b_l$; there are $e_{l-1}-e_l$ such roots
$\omega$. If $\omega^{e_{g-1}}=1$ and $\omega\neq1$, then $\omega^{b_g}\neq1$ (because
$\gcd(e_{g-1},b_g)=e_g=1$), so $a_{i,0}-\omega^{b_g}a_{j,0}\neq0$ (the numbers
$a_{i,0},a_{j,0}$ are positive and $\omega^{b_g}$ is not a positive real number) and the
order is $b_g$; there are $e_{g-1}-1$ such
roots. Finally, for $\omega=1$ the order is $b_g+d_{ij}$. By
\eqref{eq:nu-puiseux},
$(C_i\cdot C_j)=\sum_{l=1}^{g-1}(e_{l-1}-e_l)b_l+(e_{g-1}-1)b_g+b_g+d_{ij}=\nu+d_{ij}
=\xi_{ij}$.
\end{proof}

\begin{lemma}\label{lem:diagonal}
Let $\xi\ge\nu(M)$.
Then $C(M,\xi)$ is the semigroup of values of any plane curve with two
branches $C_1$, $C_2$ with semigroup $M$ and $(C_1\cdot C_2)=\xi$ (such curves exist
by Lemma~\ref{lem:realization}; they are the \emph{diagonal} curves with two
branches of \cite[(3.1.1)]{D-94}). As a consequence one has:
\begin{itemize}
\item[(1)] $C(M,\xi)$ is a symmetric good semigroup with projections $M$ and
conductor
$(c+\xi,c+\xi)$. Its maximal elements are the elements $(a,a)$ with $a$ in the
Ap\'ery set $A_\xi(M)$; in particular $C(M,\xi)$ has exactly $\xi$ maximal elements.
\item[(2)] For $a\in M$ one has $\min\{b\;:\;(a,b)\in C(M,\xi)\}=\min\{a,\xi\}$, and
if $a<\xi$ then $(a,a)$ is the only element of $C(M,\xi)$ with first coordinate
$a$. Moreover, $(\xi,\xi+1)\in C(M,\xi)$.
\end{itemize}
Since $C(M,\xi)$ is invariant under $(a,b)\mapsto(b,a)$, (2) also holds with the
coordinates exchanged.
\end{lemma}

\begin{proof}
Let $T$ be the semigroup of values of $C=C_1\cup C_2$. By \cite{D-88} it is
symmetric and, by \cite{D-87} (see also \cite{D-94}), $T_1=T_2=M$, its conductor
is $(c+\xi,c+\xi)$,
the maximal
contact values (i.e.\ the irreducible maximals) of $C$ are the elements $(m_l,m_l)$,
$0\le l\le g$ (for $g=0$
$(1,1)$ is the only maximal contact value) and every
maximal of $T$ is a sum of maximal contact values, i.e.\ it is of the form $(a,a)$
with $a\in M$.
Conversely, $(a,a)\in T$ for every $a\in M$, because the maximal
contact values $(m_l,m_l)$ belong to $T$ and $\{m_i: i=0,\ldots, g\}$ is the minimal
set of generators of $M$.
As a consequence
the maximal elements of $T$ are the elements of the form $(a,a)$ such that
$a\in A_{\xi}(M)$ (Remark~\ref{rmk r=2} and Theorem~\ref{sym-2}).

Now, the proof of
Theorem~\ref{AbsoluteMaximals} only uses that the set $A$ is contained in
$T\cap[\0,\tau]$ and contains all the absolute maximals of $T$, so, putting
$\tau_0=c+\xi-1$ (so that $\tau(T)=(\tau_0,\tau_0)$), the Theorem holds
for $T$ with $A=\{(a,a)\;:\;a\in M,\ a\le\tau_0\}$ and
$R=\{(\tau_0-a,\tau_0-a)\;:\;a\in M,\ a\le\tau_0\}$. Let $(x,y)\in M\times M$. By the
Theorem, $(x,y)\notin T$ if and only if there is $a\in M$, $a\le\tau_0$, with
$x=\tau_0-a<y$ or $y=\tau_0-a<x$; that is, if and only if $x\neq y$ and
$\tau_0-\min\{x,y\}\in M$. Since $M$ is symmetric (Remark~\ref{numerical}), for every
$z\in\Z$ one has $\tau_0-z=(c-1)-(z-\xi)\in M$ if and only if $z-\xi\notin M$.
Therefore $(x,y)\in T$ if and only if $x=y$ or $\min\{x,y\}-\xi\in M$, i.e.\
$T=C(M,\xi)$.

Obviously (1) and (2) are a consequence of the equality $T=C(M,\xi)$.
\end{proof}

\begin{lemma}\label{lem:Cn}
Let $S\subset\Zr 2$ be a symmetric good semigroup with $S_1=S_2=\Zr{}$
and with $\xi$ maximal elements. Then $S=C(\xi)$.
\end{lemma}

\begin{proof}
Here $c_1=c_2=0$. By Remark~\ref{rmk r=2} and Theorem~\ref{sym-2},
$\mu^1=\mu^2=\xi$, so $\delta=(\xi,\xi)$ and
$\Fin(S_1)=\Fin(S_2)=A_\xi(\Zr{})=\{0,\dots,\xi-1\}$.

For $a\in\Fin(S_1)$ let $\varphi(a)$ be the largest second coordinate of an element of
the (finite, nonempty) fibre of $a$. By Remark~\ref{rmk r=2}, $(a,\varphi(a))$ is a maximal
element, $\varphi(a)\in\Fin(S_2)$ and $\varphi$ is injective (if
$\varphi(a)=\varphi(a')=b$ with $a<a'$, then
$(a',b)\in\Delta_2(a,b)$). Since
$\Fin(S_1)$ and $\Fin(S_2)$ both have $\xi$ elements, $\varphi$ is a bijection of
$\{0,\dots,\xi-1\}$ onto itself.

By P4, $\varphi(1)\ge1$ (if $\xi\ge2$). As $S$ is closed under addition,
$a\cdot(1,\varphi(1))=(a,a\,\varphi(1))\in S$, hence $\varphi(a)\ge a\,\varphi(1)\ge a$ for $0\le a<\xi$. A
bijection $\varphi$ of $\{0,\dots,\xi-1\}$ with $\varphi(a)\ge a$ for all $a$ is the
identity. Thus the maximal elements of $S$ are the $(a,a)$, $0\le a<\xi$.

Let $(a,b)\in S$ with $a<\xi$. Then $b\le \varphi(a)=a$; if $b<a$, then $b\in\Fin(S_2)$ and
$(a,b)$ lies in the fibre of $b$, whose largest first coordinate is
$\varphi^{-1}(b)=b<a$, a
contradiction. Hence $S\cap\{\alpha_1<\xi\}=\{(a,a):a<\xi\}$ and, exchanging the roles
of the coordinates, $S\cap\{\alpha_2<\xi\}$ is the same set. Every other element of $S$ is
$\ge(\xi,\xi)=\delta$, and $\delta+\Zr 2\subset S$. Therefore $S=C(\xi)$.
\end{proof}

\begin{remark}\label{rem:Cn-general}
We do not know whether Lemma~\ref{lem:Cn} extends to an arbitrary semigroup $M$ of
a plane branch, i.e.\ whether every symmetric good semigroup $S\subset\Zr2$ with
$S_1=S_2=M$ and with $\xi\ge\nu(M)$ maximal elements is equal to $C(M,\xi)$. Such a
semigroup $S$ is determined by the bijection $\varphi$ of $A_\xi(M)$ defined by its
maximal elements $(a,\varphi(a))$ (Remark~\ref{rmk r=2} and Theorem~\ref{sym-2}),
and an exhaustive computer search over these bijections gives an affirmative answer
for $M=\langle2,3\rangle$, $\langle2,5\rangle$, $\langle2,7\rangle$,
$\langle3,4\rangle$, $\langle3,5\rangle$, $\langle3,7\rangle$, $\langle4,5\rangle$,
$\langle4,6,13\rangle$ and $\langle4,6,15\rangle$, and $\nu\le\xi\le\nu+3$. The bound
$\nu$ cannot be lowered in general: for $M=\langle2,7\rangle$, $\nu=14$, the semigroup
of the plane curve $(u^2,u^7)\cup(v^2,v^6+v^7)$ is symmetric, its projections are
$M$ and it has $12$ maximal elements, but it is not $C(M,12)$ (for instance, $(6,7)$
is one of its maximal elements), while the set given by the formula of
Definition~\ref{def:diagonal} with $\xi=12$ is also a symmetric good semigroup with
projections $M$.
\end{remark}

\begin{lemma}[Ultrametric inequality]\label{lem:ultra}
Let $S\subset\Zr r$ be a good semigroup such that $S_{\{i,j\}}=C(M,\xi_{ij})$, with
$\xi_{ij}\ge\nu(M)$, for all $i\neq j$ (for instance, $S$ completely symmetric with
$S_i=\Zr{}$ for all $i$, by Lemma~\ref{lem:Cn}). Then,
$(\xi_{ij} : i\neq j)$ is ultrametric.
\end{lemma}

\begin{proof}
Assume $a:=\xi_{ij}<\min\{\xi_{ik},\xi_{jk}\}$. By Lemma~\ref{lem:diagonal}(2),
$(a,a+1)\in C(M,\xi_{ij})=S_{\{i,j\}}$; let $\alpha\in S$ with $\alpha_i=a$,
$\alpha_j=a+1$. Since $(\alpha_i,\alpha_k)\in C(M,\xi_{ik})$ and $a<\xi_{ik}$,
Lemma~\ref{lem:diagonal}(2) gives $\alpha_k=a$. Then
$(\alpha_j,\alpha_k)=(a+1,a)\in C(M,\xi_{jk})$ has distinct coordinates, so
$a-\xi_{jk}\in M$, which forces $a\ge\xi_{jk}>a$, a contradiction.
\end{proof}

\begin{remark}
Lemma~\ref{lem:ultra} says that, for each $n$, the relation `$i=j$ or $\xi_{ij}> n$'
is an equivalence relation on $I$: the numbers $\xi_{ij}$ have the hierarchical
(tree-like) structure of the intersection multiplicities of the branches of a plane
curve whose branches have the same semigroup $M$ and pairwise intersection
multiplicities $\ge\nu(M)$.
It is a consequence of the axioms
of good semigroups alone. In particular, not every family of semigroups
$\{C(M,\xi_{ij})\}$ is the family of two-dimensional projections of a good
semigroup.
\end{remark}

\begin{lemma}\label{lem:maxsym}
Let $S\subset S'\subset\Zr r$ be good semigroups with the same conductor
$\delta$. If $S$ is symmetric, then $S=S'$.
\end{lemma}

\begin{proof}
Let $\alpha\in S'$. Since $\tau$ is
also the Frobenius vector of $S'$, one has $\Delta^{S'}(\tau-\alpha)=\emptyset$
(Remark~\ref{rem:maximal}, applied to $S'$). As $S\subset S'$,
$\Delta^S_i(\beta)\subset\Delta^{S'}_i(\beta)$ for all $\beta$ and $i$, so
$\Delta^S(\tau-\alpha)=\emptyset$ and, $S$ being symmetric, $\alpha\in S$
(Definition~\ref{def:symmetric}). Hence $S'=S$. (This is the argument of the
canonical ideal: $S\subset K_S$ $=\{\alpha\in\Z^r:\Delta(\tau-\alpha)=\emptyset\}$,
see e.g.\ \cite{Korell-2019}, for every good semigroup, with equality if and only if
$S$ is symmetric.)
\end{proof}

\begin{remark}
For $r=1$ this is the classical fact that a symmetric numerical semigroup is
maximal among the numerical semigroups with the same Frobenius number
(Fr\"oberg, Gottlieb and H\"aggkvist~\cite{FGH}). The general case can be seen as a
consequence of more general results on relative ideals, see e.g.\ \cite{Good}.
We include the short proof for completeness.
\end{remark}

\begin{lemma}\label{lem:diagonal-elements}
Let $S\subset\Zr r$ be a good semigroup such that $S_1=M$ and
$S_{\{1,j\}}=C(M,\xi_{1j})$, with $\xi_{1j}\ge\nu(M)$, for all $j\neq 1$. Then
$\ul{a} \in S$ for every $a\in M$.
\end{lemma}

\begin{proof}
Since $S$ is closed under addition and every $a\in M$ is a sum of the generators
$m_l$, we may assume that $a=m_l$; then $a\le\nu\le\xi_{1j}$ for all $j\neq1$ (see the
Setting). If $\beta\in S$ and $\beta_1=a$, then, for $j\neq1$,
$(a,\beta_j)\in C(M,\xi_{1j})$, so $\beta_j\ge\min\{a,\xi_{1j}\}=a$ by
Lemma~\ref{lem:diagonal}(2). On the other hand, $(a,a)\in C(M,\xi_{1j})$, so there is
$\beta^j\in S$ with $\beta^j_1=\beta^j_j=a$. By P2,
$\rho=\bigwedge_{j\neq1}\beta^j\in S$; it satisfies $\rho_1=a$ and $\rho_j\le a$ for all
$j\neq1$, hence $\rho=a\,\e$.
\end{proof}

\begin{theorem}\label{thm:diagonal-plane}
Let $M$ be the semigroup of a plane branch and let $S\subset\Zr r$, $r\ge2$, be a
completely symmetric good semigroup such that $S_i=M$ for all $i\in I$ and
$S_{\{i,j\}}=C(M,\xi_{ij})$, with $\xi_{ij}\ge\nu(M)$, for all $i\neq j$. Then the
matrix $(\xi_{ij})$ is ultrametric, and $S$ is the semigroup of values of any plane
curve $C=C_1\cup\dots\cup C_r$ whose branches have semigroup $M$ and with
$(C_i\cdot C_j)=\xi_{ij}$ for all $i\neq j$; such curves exist. Conversely, the semigroup
of values of such a curve is completely symmetric, with $S_i=M$ and
$S_{\{i,j\}}=C(M,\xi_{ij})$.
\end{theorem}

\begin{proof}
The converse follows from Lemma~\ref{lem:diagonal}(1), since every union of
branches of a plane curve is a plane curve, hence Gorenstein, and so its semigroup
of values is symmetric \cite{D-88}. By Lemma~\ref{lem:ultra} the matrix $(\xi_{ij})$
is ultrametric, so a curve $C$ as in the statement exists by
Lemma~\ref{lem:realization}.

Let $T$ be the semigroup of values of $C$. We prove $S=T$ by induction on $r$; for
$r=2$ we have $S=S_{\{1,2\}}=C(M,\xi_{12})=T$ by Lemma~\ref{lem:diagonal}. Let
$r\ge3$. Every projection $S_J$, $J\neq I$, $\#J\ge2$, is completely symmetric and
satisfies the hypotheses of the Theorem with the numbers $\xi_{ij}$, $i,j\in J$, and
$T_J$ is the semigroup of the curve $\bigcup_{j\in J}C_j$; by induction, $S_J=T_J$ for
all $J\neq I$ (for $\#J=1$, $S_J=T_J=M$). By the converse, $T$ is completely symmetric
with the same numbers $\xi_{ij}$, so, by Theorem~\ref{thm:supersymmetry}, $S$ and $T$
have the same conductor $\delta$, with $\delta_i=c+\sum_{j\neq i}\xi_{ij}$, and the
same Frobenius vector $\tau=\delta-\e$.

From \cite{D-87}, the maximal
contact values of $C$ (i.e. the irreducible absolute maximals of $T$) are the
elements $m_k\,\e$, $0\le k\le g$ (for
$g=0$, we agreed to take just $\e$ as the only one, even if $\e$ is not maximal) and, as
every absolute maximal of $T$ is a sum of maximal contact values, all of them are of
the form $a\,\e$ with $a\in M$. (In terms of the dual graph, see~\cite{D-94}, this
corresponds to the fact that all the dead arcs are common to all the branches of $C$.)

We claim that $S\subset T$. Let $\alpha\in\Zr r\setminus T$. If $\alpha_J\notin T_J=S_J$ for
some $J\neq I$, then $\alpha\notin S$. Otherwise, by Theorem~\ref{AbsoluteMaximals}
applied to $T$ (as in the proof of Lemma~\ref{lem:diagonal}, we may take
$A=\{a\,\e\;:\;a\in M\}\cap[\0,\tau]$, which is contained in $T$ by
Lemma~\ref{lem:diagonal-elements} applied to $T$), there exist $i\in I$ and $a\in M$ such that
$\alpha_i=\tau_i-a$ and $\alpha_j>\tau_j-a$ for $j\neq i$. By Lemma~\ref{lem:diagonal-elements},
$a\,\e\in S$, hence $\Delta^S(\tau-a\,\e)=\emptyset$ (Remark~\ref{rem:maximal}). If
$\alpha$ were in $S$, it would belong to $\Delta^S_i(\tau-a\,\e)$, a contradiction. Thus
$\alpha\notin S$, which proves the claim. Since $S$ is symmetric and $S\subset T$
have the same conductor, Lemma~\ref{lem:maxsym} gives $S=T$.
\end{proof}

For $M=\Zr{}$ the hypothesis on the two-dimensional projections is automatic, by
Lemma~\ref{lem:Cn} (and $\nu(\Zr{})=1$), so we get:

\begin{corollary}\label{cor:smooth-plane}
Let $S\subset\Zr r$, $r\ge2$, be a completely symmetric good semigroup with
$S_i=\Zr{}$ for all $i\in I$. Then the matrix $(\xi_{ij})$ is ultrametric, and $S$ is
the semigroup of values of any plane curve with $r$ smooth branches $C_1,\dots,C_r$
with $(C_i\cdot C_j)=\xi_{ij}$ for all $i\neq j$; such curves exist. Conversely, the
semigroup of values of a plane curve with smooth branches is completely symmetric
with all its projections $S_i$ equal to $\Zr{}$.
\end{corollary}

\begin{proof}
By Lemma~\ref{lem:Cn}, $S_{\{i,j\}}=C(\xi_{ij})=C(\Zr{},\xi_{ij})$ and
$\xi_{ij}\ge1=\nu(\Zr{})$, so Theorem~\ref{thm:diagonal-plane} applies.
\end{proof}

\begin{remark}\label{rem:diagonal-plane}
(1) The hypothesis $\xi_{ij}\ge\nu(M)$ cannot be dropped: the semigroup of
Example~\ref{ex:spacecusps} is given by the formula of Definition~\ref{def:diagonal}
with $M=\langle2,3\rangle$ and $\xi=2$; it is completely symmetric, but
$\xi=2<\nu=6$ and it is not the semigroup of a plane curve. Note also that, for
$\xi<\nu(M)$, the semigroup of two plane branches with semigroup $M$ and intersection
multiplicity $\xi$ is not given, in general, by that formula
(Remark~\ref{rem:Cn-general}).

(2) Theorem~\ref{thm:diagonal-plane} says that, among good semigroups whose
one-dimensional projections are equal to the semigroup $M$ of a plane branch and
whose two-dimensional projections are of the form $C(M,\xi_{ij})$ with
$\xi_{ij}\ge\nu(M)$, complete symmetry characterizes the semigroups of plane curves.
If the answer to the question in Remark~\ref{rem:Cn-general} is affirmative, the
second hypothesis can be replaced by the condition $\xi_{ij}\ge\nu(M)$ on the
numbers of maximal elements of the projections $S_{\{i,j\}}$, as in
Corollary~\ref{cor:smooth-plane}.
\end{remark}

\subsection{Projection formula.}   \label{sec:projection}

\begin{theorem}[{\bf Projection formula}]\label{projection}
Let $S\subset \Zr r$, $r\ge 2$, be a symmetric semigroup. Let
$P_{S^1}(t_2,\ldots, t_r)$ be the Poincar\'e series of the semigroup
$S^1 = S_{I\setminus\{1\}}$.
Assume that $S^1$ is also symmetric and
let $\eta = (\seq{\eta}2r)= pr^1(\delta(S))-\delta(S^1)\in\Zr{r-1}$.
Then one has that
\begin{equation}\label{proj-2}
P_S(1, \seq t2r) = (1-t_2^{\eta_2} \cdots t_r^{\eta_r}) P_{S^1}(\seq t2r)\; .
\end{equation}
The same formula holds for $P'_{S}$ and $P'_{S^1}$.
\end{theorem}

\begin{proof}
It is enough to prove the formula for the series $P'_S(\tt)= \sum
p'(\alpha)\tt^{\alpha}$ (resp. $P'_{S^1}$)
instead of $P_S(\tt)$ (resp. $P_{S^1}$).
The expression of the coefficients $p'(\alpha)$ is (see (\ref{Poinc-2}))
$$
p'(\alpha) = (-1)^{r}\sum_{1\notin J\subset I} (-1)^{\# J} \left( \lambda(\alpha-\e_J)
- \lambda(\alpha-\e_1-\e_J) \right)\; .
$$
Let us fix the following notations: $\eta= (\seq{\eta}{2}{r})\in \Z^{r-1}$,
$S^{1} = pr^1 (S)\subset \Zr{r-1}$ and
$P'_{S^1}(t_2, \ldots, t_r) = \sum_{\beta\in \Z^{r-1}} p'_{S^1}(\beta) \tt^{\beta}$
the corresponding series for
$S^1$ as in equation
(\ref{Poinc-1}) (note that here $\tt^{\beta}=t_2^{\beta_2}\cdots t_r^{\beta_r}$,
where $\beta = (\seq{\beta}2r)\in \Z^{r-1}$).

The left hand side of the equation in the statement of the Theorem (for $P'$) is
$$
P'_S(1,t_2, \ldots, t_r) = \sum_{\beta\in \Z^{r-1}} \left( \sum_{pr^1(\alpha)=\beta,
\alpha\in \Z^r} p'(\alpha) \right)   t_2^{\beta_2}\cdots t_r^{\beta_r}\; ,
$$
and the right hand side is
$$
\sum_{\beta\in \Z^{r-1}} \left( p'_{S^1}(\beta)-p'_{S^1}(\beta-\eta)\right)
t_2^{\beta_2}\cdots t_r^{\beta_r}\; .
$$

Let $\beta = (\seq{\beta}2r)\in \Z^{r-1}$. We know that
$p'(\alpha_1,\seq{\beta}2r) = 0$ if $\alpha_1 > \delta_1$ or $\alpha_1 < 0$, so
proving the equality (\ref{proj-2}) is equivalent to proving that
\begin{equation}\label{proj-1}
\sum_{\alpha_1=0}^{\delta_1} p'(\alpha_1, \seq{\beta}2r) = p'_{S^1}(\beta)-
p'_{S^1}(\beta-\eta) \; .
\end{equation}

For $\alpha\in \Z^r$ one has that
\begin{align*}
p'(\alpha)+p'(\alpha-\e_1)
&= (-1)^{r}\sum_{1\notin J\subset I} (-1)^{\# J} \left( \lambda(\alpha-\e_J)
- \lambda(\alpha-\e_1-\e_J) \right) \\
& \qquad  + (-1)^{r}\sum_{1\notin J\subset I} (-1)^{\# J} \left(
\lambda(\alpha-\e_1-\e_J)
- \lambda(\alpha-2 \e_1-\e_J) \right)\\
&= (-1)^{r}\sum_{1\notin J\subset I} (-1)^{\# J} \left( \lambda(\alpha-\e_J)
- \lambda(\alpha-2 \e_1-\e_J) \right) \; .
\end{align*}
As a consequence (the sum is telescopic),
\begin{equation}\label{eqn10}
\sum_{\alpha_1=0}^{\delta_1} p'(\alpha_1,\beta)
= (-1)^{r}\sum_{1\notin J\subset I} (-1)^{\# J} \left( \lambda((\delta_1,\beta)-\e_J)
- \lambda((-1,\beta) -\e_J) \right) \; .
\end{equation}
Now, $\lambda((-1,\beta)-\e_J) = \lambda_{S^1}(\beta-\e_J)$ by Remarks~\ref{colength}(1), and so
\begin{align*}
(-1)^{r+1}\sum_{1\notin J\subset I} (-1)^{\# J}  \lambda((-1,\beta) -\e_J) & =
(-1)^{r+1}\sum_{J\subset I\setminus\{1\}} (-1)^{\# J}  \lambda_{S^1}(\beta -\e_J) =\\
&= (-1)^{r-1}\sum_{J\subset I\setminus\{1\}} (-1)^{\# J}  \lambda_{S^1}(\beta -\e_J) =
p'_{S^1}(\beta)  \; .
\end{align*}

The symmetry of  $S$ (Proposition~\ref{symmetric}(3)) implies that
$$
\lambda(\delta_1,\beta) = r- \lambda(\delta-(\delta_1,\beta)-\e) =
r- \lambda((0,pr^1(\delta)-\beta)-\e)
$$
and also
$$
\lambda((\delta_1,\beta)-\e_{J}) =
r- \lambda((0,pr^1(\delta)-\beta)-\e_{I\setminus J})\; .
$$
Thus,
\begin{align*}
& (-1)^{r}\sum_{1\notin J\subset I} (-1)^{\# J} \lambda((\delta_1,\beta)-\e_J) =
(-1)^{r}\sum_{1\notin J\subset I} (-1)^{\# J} \left(
r- \lambda((0,pr^1(\delta)-\beta)-\e_{I\setminus J}) \right) = \\
& \qquad  =
(-1)^{r}\sum_{1\notin J\subset I} (-1)^{\# J} r -
(-1)^{r}\sum_{1\notin J\subset I} (-1)^{\# J}
\lambda((0,pr^1(\delta)-\beta)-\e_{I\setminus J}) = \\
& \qquad  =
0 +
(-1)^{r-1}\sum_{1\in K\subset I} (-1)^{r- \# K} \lambda((0,pr^1(\delta)-\beta)-\e_{K}) = \\
&  \qquad =
\sum_{K\subset I\setminus \{1\}} (-1)^{\# K} \lambda((-1,pr^1(\delta)-\beta)-\e_{K}) =
\sum_{K\subset I\setminus \{1\}} (-1)^{\# K}
\lambda_{S^1}(pr^1(\delta)-\beta-\e_{K}) = \\
& \qquad =
\sum_{K\subset I\setminus \{1\}} (-1)^{\# K}
\lambda_{S^1}(\delta(S^1) + \eta-\beta-\e_{K}) \qquad \text{(by the definition of $\eta$)} \\
& \qquad = \sum_{K\subset I\setminus \{1\}} (-1)^{\# K}
\lambda_{S^1}(\delta(S^1) - (\beta- \eta) -\e_{K}) = (*)
\end{align*}
Now, by the symmetry of the semigroup $S^1$ one has
\begin{align*}
& \lambda_{S^1}(\delta(S^1) - (\beta- \eta) -\e_{K}) =
(r-1) - \lambda_{S^1}(\beta- \eta -\e_{(I\setminus \{1\})\setminus K}) =\\
& =
(r-1) - \lambda_{S^1}(\beta- \eta -\e_{L})
\end{align*}
where $L = \{2,\ldots, r\}\setminus K$; note also that
$\# L = (r-1)-\# K$ and so $(-1)^{\# K} = (-1)^{r-1}(-1)^{\# L}$.
Thus, for the equation $(*)$ above we have
\begin{align*}
(*) & =
\sum_{K\subset I\setminus \{1\}} (-1)^{\# K} (r-1) -
\sum_{L\subset I\setminus \{1\}} (-1)^{r-1}(-1)^{\# L}
\lambda_{S^1}((\beta- \eta) -\e_{L}) =  \\
& = 0 - (-1)^{r-1}
\sum_{L\subset I\setminus \{1\}} (-1)^{\# L}
\lambda_{S^1}((\beta- \eta) -\e_{L}) =  - p'_{S^1}(\beta-\eta) \; .
\end{align*}
Altogether we have proved that in equation (\ref{eqn10}) we have
$$
\sum_{\alpha_1=0}^{\delta_1} p'(\alpha_1,\beta) = - p'_{S^1}(\beta-\eta) + p'_{S^1}(\beta)
$$
and the proof is finished (see equation (\ref{proj-1})).
\end{proof}

If $S$ is completely symmetric, then
$pr^1(\delta)-\delta(S^1) = (\seq{\xi}{1,2}{1,r})$ by
Theorem~\ref{thm:supersymmetry} (see (\ref{eq:prJ-delta})).
Moreover, one can replace the index $1$ by any $i\in I$. Thus one has:

\begin{corollary}\label{cor:projection1}
Let $S$ be a completely symmetric semigroup. Then, for any $i\in I$, one has
\begin{equation}\label{proj-3}
P_S(\seq t1{r})|_{t_i=1} = (1- \prod_{j\neq i}t_j^{\xi_{j,i}}) P_{S^i}(\{\seq t1r
\}\setminus \{t_i\})\; .
\end{equation}
\end{corollary}

\noindent
By iterating Corollary~\ref{cor:projection1} (or, more precisely, equation
(\ref{proj-3})) one gets:

\begin{corollary}\label{cor:projection2}
Let $S$ be completely symmetric and
$A = \{ \seq i1a\}\subseteq I$. Then one has
$$
P_S(\tt) |_{\{t_i=1 : i\notin A\}} =
P_{S_A}(\seq t{i_1}{i_a})
\prod_{j\notin A} (1-\tt_{A}^{\xi_{jA}})\; ,
$$
where
$\tt_{A}^{\xi_{jA}}= \prod_{i\in A} t_i^{\xi_{ij}}$
and $P_{S_A}(\seq t{i_1}{i_a})$ is the Poincar\'e series of the
semigroup $S_A$.
\end{corollary}
\begin{proof}
Removing the index $j_1$ produces the factor
$(1-\prod_{i\neq j_1}t_i^{\xi_{j_1i}})$, and the later substitution $t_{j_2}=1$ turns
it into $(1-\prod_{i\neq j_1,j_2}t_i^{\xi_{j_1i}})$, and so on.
\end{proof}

\begin{remark}
The projection formula is the generalization of the well-known projection formula of
Torres (see \cite{Torres}) for the Alexander polynomial of the link of a plane curve
singularity. Given the equality of the Alexander polynomial and the Poincar\'e
series for plane curve singularities (see \cite{Duke}), Theorem \ref{projection}
establishes the intrinsic properties of the semigroup $S$ which are behind the projection
formula.

In \cite{Moyano} the formula~(\ref{proj-3}) is stated for Gorenstein curve
singularities, i.e.\ essentially for symmetric good semigroups. The proof there
follows the same idea as the above proof, based on the computation of the coefficients
of $P_S(\tt)$ or $P'_S(\tt)$ given in (\ref{Poinc-2}), which was described for
the first time in \cite{Duke}. However, the statement is not true under the sole
assumption that the ring is Gorenstein: it requires the same condition on all the
projection curves, that is, the requirement that $S$ be completely symmetric. Let
us show an example.
\end{remark}

\begin{example}\label{ex:two-branch}

{\bf A two-branch counterexample.}

Let $C_1\subseteq \C^3$ be the monomial curve $t\mapsto(t^3,t^4,t^5)$, let
$C_2\subset
\C^3$ be the $z$-axis, and let $C = C_1\cup C_2$ be the reduced space curve
germ with
components $C_1$ and $C_2$. The semigroup of $C_1$ is the non-symmetric semigroup
$S_1=\langle 3,4,5 \rangle\subset \Zr{}$, so $C_1$ is not Gorenstein.
Obviously, the semigroup of $C_2$ is $S_2=\Zr{}$.

Let $\OO = \C[[x,y,z]]/{\mathfrak I}$ denote the local ring
of functions of $C$. The normalization morphism of $\OO$ is given by the map
$\varphi\colon\OO \to \C[[t]]\times \C[[u]]$,
$g\mapsto\bigl(g(t^3,t^4,t^5),\,g(0,0,u)\bigr)$. The complete semigroup $\b{S}$ of
$C$ (see Remark~\ref{complete}) is $\{(\ord_t (g(t^3,t^4,t^5)), \ord_u(g(0,0,u)))\}$
for $g\in \C[[x,y,z]]$, $g\notin\mathfrak I$ (we assume $\ord(0) = \infty$).

The ideal of $C_2$ is $(x,y)$, so one has that
$\varphi: xg(x,y,z)\mapsto (t^3g(t^3,t^4,t^5), 0)$ and, as a consequence,
$(3+m,\infty)\in \b{S}$ for all $m\in S_1$.
Also,  $yg\in \C[[x,y,z]]$ gives $(4+m,\infty)\in \b{S}$ for all
$m\in S_1$.
In a similar way,
$z^2-x^2y$ vanishes on $C_1$, so $(z^2-x^2y)g\mapsto (0,u^2g(0,0,u))$ and
$(\infty, 2+m)\in \b S$ for all $m\in S_2$.
In particular, we have
$\{(3,b), (4,b): b\ge 1\}\cup ((6,2)+\Zr 2)\subset S$.
(The elements $(3,b)$ and $(4,b)$ are the values of $x+z^b$ and $y+z^b$.)

Since $z \mapsto (t^5, u)$, one has  $(5,1)\in S$.
Now, if $(5,b)\in S$ for some $b\ge 2$, then $(5,\infty)\in\b S$ (see
Notations~\ref{not1} and Remark~\ref{complete}). Let $h \in \C[[x,y,z]]$ be such that
its value is $(5,\infty)$. Since $h(0,0,u)=0$, $h$ belongs to the ideal $(x,y)$ of
$C_2$, so we can write $h=xh_1+yh_2$. One has $\ord_t(xh_1)\in 3+S_1$ and
$\ord_t(yh_2)\in 4+S_1$: if these two orders are different, $\ord_t(h)$ is their
minimum; if they coincide, their common value belongs to
$(3+S_1)\cap(4+S_1)=\{7,8,\dots\}$, and then $\ord_t(h)\ge 7$. In both cases
$\ord_t(h)\neq5$, because $5\notin(3+S_1)\cup(4+S_1)$, a contradiction. Thus
$\Delta_1(5,1)=\emptyset$ and, by P3 (see Remark~\ref{rmk r=2}), also
$\Delta_2(5,1)=\emptyset$, i.e.\ $\Delta(5,1)=\emptyset$.

This means that
$(6,2)$ is the conductor of $S$ and $\tau=(5,1)$ is maximal. The semigroup is the
symmetric semigroup (see also Figure~\ref{fig3}):
$$
S= \{(0,0)\}\cup \{(3,b), (4,b): b\ge 1\}\cup \{(5,1)\}\cup ((6,2)+\Zr 2)\subset
\Zr 2\; .
$$

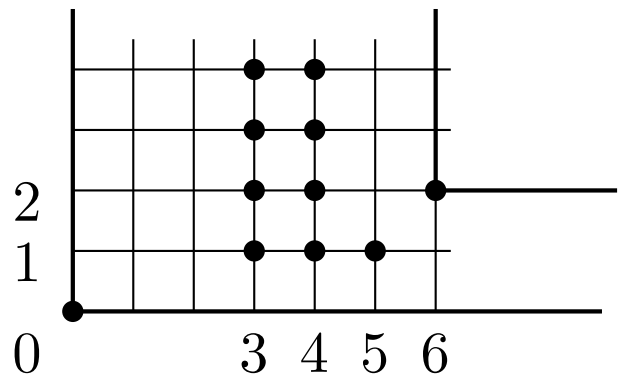
\begin{figure}[h]
\begin{center}
\unitlength=1.00mm
\begin{picture}(60,20)(0,10)
\thicklines
\put(10,10){\line(1,0){35}}
\put(10,10){\line(0,1){20}}
\put(6,6){$0$}
\put(21,6){$3$}
\put(25,6){$4$}
\put(29,6){$5$}
\put(33,6){$6$}
\put(6,12){$1$}
\put(6,16){$2$}
\thinlines
\put(14,10){\line(0,1){18}}
\put(18,10){\line(0,1){18}}
\put(22,10){\line(0,1){18}}
\put(26,10){\line(0,1){18}}
\put(30,10){\line(0,1){18}}
\put(34,10){\line(0,1){18}}
\put(10,14){\line(1,0){25}}
\put(10,18){\line(1,0){25}}
\put(10,22){\line(1,0){25}}
\put(10,26){\line(1,0){25}}
\put(10,10){\circle*{1.5}}
\put(22,14){\circle*{1.5}}
\put(22,18){\circle*{1.5}}
\put(22,22){\circle*{1.5}}
\put(22,26){\circle*{1.5}}
\put(26,14){\circle*{1.5}}
\put(26,18){\circle*{1.5}}
\put(26,22){\circle*{1.5}}
\put(26,26){\circle*{1.5}}
\put(30,14){\circle*{1.5}}
\put(34,18){\circle*{1.5}}
\thicklines
\put(34,18){\line(1,0){12}}
\put(34,18){\line(0,1){12}}
\end{picture}
\caption{The semigroup $S$ of Example~\ref{ex:two-branch}}
\label{fig3}
\end{center}
\end{figure}

As a consequence, $\OO$ is Gorenstein; in fact, $C$ is a complete intersection: the
ideal of $C$ is ${\mathfrak I}=(xz-y^2,x^3-yz)\subset \C[[x,y,z]]$.

\medskip

The Poincar\'e series of $S$ is $P_S(\tt) = 1+t_1^5t_2$ and that of $S_1$ is
$P_{S_1} = \sum_{a\in S_1}t_1^a =(1-t_1+t_1^3)/(1-t_1)$.
If there exists $\alpha\in \Z$ such that
$$
P_S(t_1,1) = 1+t_1^5 = (1-t_1^\alpha) \frac{1-t_1+t_1^3}{1-t_1} =
(1-t_1^\alpha)P_{S_1}(t_1)\;
$$
then,
evaluating at $t_1=1$ forces $\alpha=2$, and
$(1+t_1)(1-t_1+t_1^3)=1-t_1^2+t_1^3+t_1^4\neq1+t_1^5$ which is a contradiction.

Note that for the other projection the
formula does hold, in the form of Theorem~\ref{projection}:
$P_S(1,t_2)=1+t_2=(1-t_2^{2})P_{\Zr{}}(t_2)$ with $\eta_2=\delta_2-\delta(S_2)=2$,
because $S$ and $S_2=\Zr{}$ are symmetric.
Here $\mu^1=3\neq\kappa_1=\xi=2$: as in Example~\ref{ex4}, the symmetry of
$S_1$ cannot be dropped in the second part of Theorem~\ref{sym-2}, nor the
symmetry of the projections in Theorem~\ref{thm:supersymmetry}.

\end{example}

Let us now consider the converse of Theorem~\ref{projection}. A first consequence
of the projection formula is the following.

\begin{corollary}\label{cor:proj-FE}
Let $S\subset\Zr r$, $r\ge2$, be a symmetric semigroup and assume that
$$
P_S(1,t_2,\dots,t_r)=(1-t_2^{\eta_2}\cdots t_r^{\eta_r})\,P_{S^1}(t_2,\dots,t_r)
$$
for some $\eta\in\Z^{r-1}\setminus\{\0\}$ (not necessarily equal to
$pr^1(\delta(S))-\delta(S^1)$). Then, with $\tt=(t_2,\dots,t_r)$,
\begin{equation}\label{eq:fe-proj}
P_{S^1}(\tt)=(-1)^{r-1}\,\tt^{\,pr^1(\tau)-\eta}\,P_{S^1}(\tt^{-1}) ,
\end{equation}
that is, $P_{S^1}$ satisfies a functional equation as in
Theorem~\ref{thm:duality-maximals}(3). The same holds with any index $i\in I$
instead of $1$.
\end{corollary}

\begin{proof}
Substituting $t_1=1$ in the functional equation of
Theorem~\ref{thm:duality-maximals}(3) one gets
$P_S(1,\tt)=(-1)^r\,\tt^{\,pr^1(\tau)}P_S(1,\tt^{-1})$. Replacing both sides by the
projection formula and using $1-\tt^{-\eta}=-\tt^{-\eta}(1-\tt^{\eta})$, this becomes
$$
(1-\tt^{\eta})\,P_{S^1}(\tt)=(-1)^{r-1}\,\tt^{\,pr^1(\tau)-\eta}\,(1-\tt^{\eta})\,
P_{S^1}(\tt^{-1}),
$$
an equality of rational functions. Since $\eta\neq\0$, the factor $1-\tt^{\eta}$ is
not zero and can be cancelled.
\end{proof}

\noindent
If $\eta=pr^1(\delta(S))-\delta(S^1)$, as in Theorem~\ref{projection}, then
$pr^1(\tau)-\eta=\delta(S^1)-pr^1(\e)=\tau(S^1)$, and (\ref{eq:fe-proj}) is exactly the
functional equation of Theorem~\ref{thm:duality-maximals}(3) for $S^1$, as it must
be when $S^1$ is symmetric.

For $r=2$ the projection $S^1$ is a numerical semigroup and, by
Remark~\ref{rem:PS-sym-not-S}, the functional equation (\ref{eq:fe-proj}) is
equivalent to its symmetry. Hence, in this case, the converse of
Theorem~\ref{projection} holds:

\begin{proposition}\label{thm:two}
Let $S\subset \Zr 2$ be a symmetric good semigroup (that is, in the case of curves,
$C$ is Gorenstein) and let $c_1$
be the conductor of the numerical semigroup $S_{1}$. There exists $\xi\in\Z$ with
\[
P_S(t,1)=\bigl(1-t^{\xi}\bigr)P_{S_{1}}(t)
\]
if and only if $S_{1}$ is symmetric (that is, if and only if the branch
$C_{1}$ is Gorenstein). In that case $\xi=\delta_1-c_1$. The same holds
after exchanging the two branches.
\end{proposition}

\begin{proof}
($\Leftarrow$) It is Theorem~\ref{projection} with the index $2$ removed
($S$ and $S_1$ symmetric, $\eta=\delta_1-c_1$).

\noindent
($\Rightarrow$) For $r=2$ one has $P_S(t,1)=\sum_{a\in\Fin(S_1)}t^a$
(Proposition~\ref{prop:coef} and Remark~\ref{rmk r=2}). In particular its constant
term is $1$, so $\xi\neq0$. By Corollary~\ref{cor:proj-FE} (with the index $2$
instead of~$1$),
$$
P_{S_1}(t)=-\,t^{\,\tau_1-\xi}\,P_{S_1}(t^{-1}),
$$
and by Remark~\ref{rem:PS-sym-not-S} the semigroup $S_1$ is symmetric and
$\tau_1-\xi=c_1-1$, i.e.\ $\xi=\delta_1-c_1$.
\end{proof}

\noindent
In that case $\xi=\delta_1-c_1=\mu^1$ and, since $S$ and $S_1$ are symmetric,
Proposition~\ref{sym-1} and Theorem~\ref{sym-2} give $\xi\in S_1$ and
$\Fin(S_1)=A_\xi(S_1)$; in particular $\xi=\kappa_1$ is the number of maximal elements
of $S$ (Remark~\ref{rmk r=2}). Since $\xi\in S_1$,
$(1-t^{\xi})P_{S_1}(t)=\sum_{a\in A_\xi(S_1)}t^a$, so the formula of
Proposition~\ref{thm:two} is just the equality $\Fin(S_1)=A_{\xi}(S_1)$ of
Theorem~\ref{sym-2}. If, moreover, $S_2$ is symmetric, then $S$ is completely
symmetric and $\xi=\xi_{12}$ (Theorem~\ref{thm:supersymmetry}); note that this need not
be the case: in Example~\ref{ex:two-branch} the formula holds for the projection onto
the second coordinate, but $S_1$ is not symmetric.

For $r\ge3$, however, the functional equation (\ref{eq:fe-proj}) does not imply that
$S^1$ is symmetric (Remark~\ref{rem:PS-sym-not-S}), and the converse of
Theorem~\ref{projection} fails, even with $\eta=pr^1(\delta(S))-\delta(S^1)$, as the
following example shows.

\begin{example}\label{ex:proj-converse}
Let $C\subset\C^3$ be the curve with branches $C_1:(0,0,t)$, $C_2:(s^2,s^3,s^4)$ and
$C_3:(u,0,0)$, and let $S$ be its semigroup of values. One finds $\delta=(2,5,3)$ and
$\tau=(1,4,2)$; the maximal elements of $S$ are $\0$ and $\tau$, and
$2\ell(\delta)=10=\abs\delta$, so $S$ is symmetric (and $C$ is Gorenstein), with
$P_S(\tt)=1-t_1t_2^4t_3^2 .$
The projection $S^1$ is the semigroup of $C_2\cup C_3$:
$$
S^1=\{(0,0),(2,1)\}\cup\big((3,2)+\Zr2\big),\qquad \delta(S^1)=(3,2),\qquad
P_{S^1}=1+t_2^2t_3 .
$$
It is not symmetric: $(1,1)\notin S^1$, but $\Delta(\tau(S^1)-(1,1))=\Delta(1,0)=\emptyset$.
Nevertheless, with $\eta=pr^1(\delta)-\delta(S^1)=(2,1)$,
$$
P_S(1,t_2,t_3)=1-t_2^4t_3^2=(1-t_2^2t_3)(1+t_2^2t_3)=(1-t_2^{\eta_2}t_3^{\eta_3})\,
P_{S^1}(t_2,t_3).
$$
So, for $r\ge3$, the symmetry of $S^1$ is sufficient but not necessary for the
projection formula. Note that $P_{S^1}=1+t_2^2t_3$ does satisfy the functional
equation (\ref{eq:fe-proj}), as in Remark~\ref{rem:PS-sym-not-S}.
\end{example}

\section{Recovering the semigroup from its Poincar\'e series}
\label{sec:recovering}

This section is devoted to the problem of determining to what extent the Poincar\'e
series $P_S(\tt)$ of a good semigroup is equivalent to the semigroup itself.
This is the combinatorial counterpart of the theorem of Yamamoto \cite{Yamamoto}:
the Alexander
polynomial determines the topological type of a plane curve singularity. As we
will see
(Theorem~\ref{thm:recons}),
the key step is how (or whether) one can recover the numbers $\xi_{ij}$ from $P_S$.
The reason is that, once these numbers are known, the projection formula
(Corollary~\ref{cor:projection2}) allows one to compute
the Poincar\'e series of all the projections and then to proceed by induction.

If $r=2$ everything is elementary, because
$\xi_{12}=P_S(1,1)$ is the number of maximal elements (i.e.\ the number of monomials in
the support of $P_S$). For $r=3$,
Theorem~\ref{thm:xi-r3} gives a positive answer (see also
Corollary~\ref{cor:PS-semigroup}). However, in the general case $r\ge 4$ the problem
remains open.

\begin{lemma}[Reconstruction]\label{lem:recons}
A good semigroup $S$ is
determined by its Poincar\'e series $P_S(\tt)$, by its conductor $\delta$ and by the
set $S\setminus[\0,\tau]$. Moreover, $S\setminus[\0,\tau]$ is determined by the
projections $S_J$, $J\neq I$, and their elements with infinite fibre.
\end{lemma}

\begin{proof}
Let $\alpha\in[\0,\tau]$ and assume that we already know
which elements $\beta\in\Z^r_{\ge0}$ with $\beta>\alpha$ belong to $S$ (we go
through $[\0,\tau]$ in decreasing order along a linear extension of the product
order; the elements outside $[\0,\tau]$ are known by hypothesis). For $i\in I$ let
$K_i(\alpha)=\{\beta\in S:\beta\neq\alpha,\ \beta_i=\alpha_i,\ \beta\ge\alpha\}$. If
$K_i(\alpha)\neq\emptyset$ for all $i$, pick $\beta^i\in K_i(\alpha)$; then
$\alpha=\bigwedge_i\beta^i\in S$ by P2. If $K_i(\alpha)=\emptyset$ for some $i$, then
$\w\Delta_i(\alpha+\e_X)\subset K_i(\alpha)=\emptyset$ for every $\emptyset\neq X\subset
I\setminus\{i\}$ and $\w\Delta_i(\alpha)=\{\alpha\}\cap S$, so that
$\nabla_{\alpha,i}=\{\emptyset\}$ if $\alpha\in S$ and $\nabla_{\alpha,i}=\emptyset$
otherwise; by Theorem~\ref{coef-PS} (for arbitrary $\alpha\in\Z^r$),
$\alpha\in S$ if and only if $p(\alpha)=1$.

Let us see that $S\setminus[\0,\tau]$ is determined by the
projections and their finite (or infinite) fibres.
If $\beta\notin[\0,\tau]$ and
$K=\{k:\beta_k\ge\delta_k\}\neq\emptyset$, then, by the proof of
Proposition~\ref{prop:small}, $\beta\in S$ if and only if $\beta_{I\setminus K}*\delta\in S$,
i.e.\ if and only if $\beta_{I\setminus K}\in S_{I\setminus K}$ has infinite fibre
(for $K=I$, $\beta\ge\delta$ is always in $S$).
\end{proof}

\begin{theorem}\label{thm:recons}
Let $S\subset \Zr r$ be a
completely symmetric semigroup.
Then, the pair consistingf of its Poincar\'e series $P_S(\tt)$ and the
arrangement of integers $(\xi_{ij})$ determines $S$.
\end{theorem}

\begin{proof}
Since all the numbers $\{\xi_{ij} \st i\neq j\}$ are known, by the projection formula
(Corollary~\ref{cor:projection2}) one can compute the Poincar\'e series $P_{S_J}$ for all $J\subset I$,
$J\neq\emptyset$. Inductively
(the case $r=1$ being trivial, since $P_{S_i}(t)=\sum_{s\in S_i}t^s$) one can assume that $S_J$ is known for all
$J\subset I$, $J\neq I$. Then
Theorem~\ref{thm:supersymmetry} gives the conductor $\delta$, and the elements with infinite fibre of each $S_J$
are given by Proposition~\ref{prop:Apery-fibres}, so $S\setminus[\0,\tau]$ is known
(Notations~\ref{not1}, condition 3), and Lemma~\ref{lem:recons} concludes the
proof.
\end{proof}

\begin{notations}
For a numerical semigroup $\Gamma\subset\mathbb{Z}_{\ge 0}$ with (finite) set of
gaps $G(\Gamma)=\mathbb{Z}_{\ge 0}\setminus\Gamma$ we have
\[
P_\Gamma(t)=\sum_{s\in\Gamma}t^s=\frac{1}{1-t}-\sum_{g\in G(\Gamma)}t^g ,
\]
so that
\[
Q_\Gamma(t):=(1-t)\,P_\Gamma(t)=1-(1-t)\sum_{g\in G(\Gamma)}t^g\in\mathbb{Z}[t]
\qquad\text{and}\qquad Q_\Gamma(1)=1 .
\]
For an integer $n\ge 1$ we write $[n]_t=1+t+\dots+t^{n-1}$, so that
$1-t^n=(1-t)[n]_t$ and $[n]_1=n$. Recall also that, if $S$ is a completely symmetric
semigroup then $\xi_{ij}\ge 1$ for
$i\neq j$, because $\underline{0}$ is always a maximal element of~$S_{\{i,j\}}$.
\end{notations}

\begin{lemma}\label{lem:one-variable}
Let $S\subset\mathbb{Z}^r_{\ge 0}$, $r\ge 2$, be a completely symmetric semigroup,
let $i\in I$, and let $f_i(t)\in\mathbb{Z}[t]$ be the polynomial obtained from
$P_S(\underline{t})$ by substituting $t_i=t$ and $t_j=1$ for every $j\neq i$.
Then
\[
f_i(t)=(1-t)^{r-2}\,Q_{S_i}(t)\prod_{j\neq i}[\xi_{ij}]_t .
\]
In particular $(1-t)^{r-2}$ divides $f_i(t)$, and the polynomial
$g_i(t)=f_i(t)/(1-t)^{r-2}$ satisfies
\[
g_i(1)=M_i:=\prod_{j\neq i}\xi_{ij},
\qquad\text{equivalently}\qquad
M_i=\frac{(-1)^r}{(r-2)!}\,f_i^{(r-2)}(1).
\]
(Here $f_i^{(r-2)}$ denotes the $(r-2)$-th derivative of $f_i$.)
\end{lemma}

\begin{proof}
By Corollary~\ref{cor:projection2} with $A=\{i\}$,
\[
f_i(t)=P_{S_i}(t)\prod_{j\neq i}\bigl(1-t^{\xi_{ij}}\bigr)
=\frac{Q_{S_i}(t)}{1-t}\,(1-t)^{r-1}\prod_{j\neq i}[\xi_{ij}]_t ,
\]
an identity in $\mathbb{Z}[[t]]$ between polynomials, hence in $\mathbb{Z}[t]$.
Evaluating $g_i(t)=Q_{S_i}(t)\prod_{j\neq i}[\xi_{ij}]_t$ at $t=1$ gives
$g_i(1)=M_i$. Finally, $f_i(t)=(-1)^{r-2}(t-1)^{r-2}g_i(t)$, so the
$(r-2)$-th derivative of $f_i$ at $t=1$ is $(-1)^{r}(r-2)!\,g_i(1)$.
\end{proof}

\begin{theorem}\label{thm:xi-r3}
Let $S\subset\mathbb{Z}^3_{\ge 0}$ be a completely symmetric semigroup.
Let
\[
M_i=-\frac{\partial P_S}{\partial t_i}(1,1,1)\qquad(i=1,2,3).
\]
Then $M_i=\xi_{ij}\,\xi_{ik}$ for $\{i,j,k\}=\{1,2,3\}$, and
\[
\xi_{12}=\sqrt{\frac{M_1M_2}{M_3}},\qquad
\xi_{13}=\sqrt{\frac{M_1M_3}{M_2}},\qquad
\xi_{23}=\sqrt{\frac{M_2M_3}{M_1}} .
\]
In particular, the numbers $\xi_{ij}$ are determined by $P_S$.
\end{theorem}

\begin{proof}

For $r=3$  and $i\in I$, the chain rule gives
$f_i'(1)=\frac{\partial P_S}{\partial t_i}(1,1,1)$, so Lemma~\ref{lem:one-variable}
yields $M_i=-f_i'(1)=\prod_{j\neq i}\xi_{ij}=\xi_{ij}\xi_{ik}$. Hence,
$M_1M_2/M_3=(\xi_{12}\xi_{13})(\xi_{12}\xi_{23})/\xi_{13}\xi_{23}=\xi_{12}^2$, and
$\xi_{12}$ is the positive square root because $\xi_{12}\ge 1$. The other two
formulas are obtained in the same way.
\end{proof}

\begin{corollary}\label{cor:projections-r3}
Let $S\subset\mathbb{Z}^3_{\ge 0}$ be completely symmetric. Then $P_S$ determines
$P_{S_A}$ for every nonempty $A\subsetneq I$, and therefore the numerical semigroups
$S_1,S_2,S_3$. Explicitly, for $\{i,j,k\}=\{1,2,3\}$,
\[
P_{S_i}(t)=\frac{f_i(t)}{(1-t^{\xi_{ij}})(1-t^{\xi_{ik}})},
\qquad
P_{S_{\{i,j\}}}(t_i,t_j)=\frac{P_S(\underline{t})\big|_{t_k=1}}{1-t_i^{\xi_{ki}}\,t_j^{\xi_{kj}}},
\]
where the quotients are taken in $\mathbb{Z}[[t]]$ and $\mathbb{Z}[[t_i,t_j]]$.
\end{corollary}

\begin{proof}
The $\xi_{ij}$ are known by Theorem~\ref{thm:xi-r3}. The formulas are Corollary~\ref{cor:projection2}
with $A=\{i\}$ and Corollary~\ref{cor:projection1} (with the index $k$ removed), and the divisions are
possible because the denominators have constant term $1$. Finally, a numerical
semigroup $\Gamma$ is determined by $P_\Gamma(t)=\sum_{s\in\Gamma}t^s$.
\end{proof}

\begin{corollary}\label{cor:PS-semigroup}
Let $S\subset \Zr 3$ be a completely symmetric semigroup. Then
$S$ is completely determined by its Poincar\'e series $P_S(\tt)$.
\end{corollary}

\begin{proof}
By Theorem~\ref{thm:xi-r3} one knows $\{\xi_{12}, \xi_{13}, \xi_{23}\}$.
So, the result is Theorem~\ref{thm:recons} applied to $r=3$ and these numbers.
\end{proof}

\begin{example}
Let $S$ be the semigroup of values of the plane curve with branches
$C_1:(t^2,t^3)$, $C_2:(t^3,t^2)$ and $C_3:(t,t)$, so that $\xi_{12}=4$ and
$\xi_{13}=\xi_{23}=2$. One finds
\[
P_S(\underline{t})=(1-t_1^2t_2^2t_3)(1+t_1^2t_2^3t_3)(1+t_1^3t_2^2t_3).
\]
Since the first factor vanishes at $(1,1,1)$, only its derivative contributes, and
\[
M_1=2\cdot 2\cdot 2=8,\qquad M_2=2\cdot 2\cdot 2=8,\qquad M_3=1\cdot 2\cdot 2=4,
\]
which give $\xi_{12}=\sqrt{8\cdot 8/4}=4$ and $\xi_{13}=\xi_{23}=\sqrt{8\cdot4/8}=2$.
Moreover $f_1(t)=(1-t)\,(1-t+t^2)\,[4]_t\,[2]_t$, where $1-t+t^2=Q_{S_1}(t)$ with
$S_1=\langle 2,3\rangle$, as predicted by Lemma~\ref{lem:one-variable}.
\end{example}

\begin{remark}\label{rem:r4}
Lemma~\ref{lem:one-variable} holds for every $r\ge 2$, but it only provides the
$r$ products $M_i=\prod_{j\neq i}\xi_{ij}$, while there are $\binom r2$ numbers
$\xi_{ij}$. For $r\ge 4$ the one-variable specializations of $P_S$ do not determine
the $\xi_{ij}$ in general, even for semigroups of plane curves. Let $S$ be the semigroup
of Example~\ref{ex:p-zero}, for which
$(\xi_{12},\xi_{13},\xi_{14},\xi_{23},\xi_{24},\xi_{34})=(2,1,1,1,1,2)$, and let $S'$ be
the semigroup obtained from $S$ by exchanging the indices $2$ and $3$ (i.e.\ the
semigroup of the same curve with the branches $C_2$ and $C_3$ renumbered), for which
$(\xi'_{12},\dots,\xi'_{34})=(1,2,1,1,2,1)$. All the projections $S_i$, $S'_i$ are equal
to $\Zr{}$ and, for each $i$, the numbers $\xi_{ij}$ and $\xi'_{ij}$, $j\neq i$, are
$2,1,1$ up to order; hence, by Lemma~\ref{lem:one-variable},
$f_i(t)=f'_i(t)=(1-t)^2(1+t)$ for $i=1,\dots,4$. However, by
Corollary~\ref{cor:projection2} with $A=\{1,2\}$,
\[
P_S(t_1,t_2,1,1)=(1+t_1t_2)(1-t_1t_2)^2\neq(1-t_1^2t_2)(1-t_1t_2^2)=P_{S'}(t_1,t_2,1,1).
\]
Hence, for $r\ge 4$, an analogue of Theorem~\ref{thm:xi-r3} requires more information
from $P_S$ than its one-variable specializations; the two-variable specializations
$P_S|_{t_k=1,\,k\neq i,j}$, which by Corollary~\ref{cor:projection2} are equal to
$P_{S_{\{i,j\}}}(t_i,t_j)\prod_{k\neq i,j}(1-t_i^{\xi_{ki}}t_j^{\xi_{kj}})$, are the
natural candidates.
\end{remark}

\medskip
\noindent{\bf Open problem.}
Is a completely symmetric semigroup $S\subset\Zr r$ always determined by its
Poincar\'e series $P_S(\tt)$? By Theorem~\ref{thm:recons} the question is whether
$P_S$ determines the numbers $\xi_{ij}$. The answer is affirmative for $r\le3$
(Theorem~\ref{thm:xi-r3}) and, through the identification of $P_S$ with the
Alexander polynomial \cite{Duke} and of $S$ with the topological type, within the
class of semigroups of values of plane curves (Yamamoto's theorem
\cite{Yamamoto}); by Theorem~\ref{thm:diagonal-plane}, this class contains all the
completely symmetric semigroups with $S_i=M$ for every $i$, $M$ the semigroup of a
plane branch, and diagonal two-dimensional projections (in particular, those with
$S_i=\Zr{}$ for every $i$, Corollary~\ref{cor:smooth-plane}). For $r\ge4$ the
question is open in general: Remark~\ref{rem:r4} shows that the one-variable
specializations of $P_S$ are not enough, and Example~\ref{ex:spacecusps} shows that,
when the projections $S_i$ are not trivial, one cannot reduce the problem to plane
curves.
For instance, it is not known whether a completely
symmetric semigroup with some nontrivial projection $S_i$ can have the same
Poincar\'e series as the semigroup of a plane curve with smooth branches; note that
the one-variable specializations cannot distinguish them, since
$Q_{\langle2,3\rangle}(t)\,[2]_t\,[3]_t=[6]_t$ (compare Lemma~\ref{lem:one-variable}:
$S_1=\langle2,3\rangle$ with $(\xi_{12},\xi_{13},\xi_{14})=(2,3,x)$ and $S_1=\Zr{}$ with
$(6,x,1)$ give the same polynomial $f_1$).

\subsection*{Declaration on the use of generative AI}
During the preparation of this work the author used Claude (Anthropic), a generative
AI model, as a research and editing assistant, through several rounds of review of
successive versions of the manuscript (October and November 2026). Its main
contributions were the following.
(a) A detailed check of the proofs of each version, which pointed out some errors
and gaps and proposed corrections (for instance, the present proof of
Lemma~\ref{lem:new} or the final proof of Proposition~\ref{prop:chain} and
Proposition~\ref{prop:Apery-fibres}).
(b) Proofs, or first versions of proofs, which the author checked and adopted:
the index-free formula for
the coefficients (Corollary~\ref{cor:indexfree});
the characterization of
symmetry by the Alexander duality of the complexes $\nabla_{\alpha,i}$
(Proposition~\ref{teo-dualidad}) and
and Example~\ref{ex:Sprime};
the discussion of the converse of the projection
formula (Remark~\ref{rem:PS-sym-not-S},
Corollary~\ref{cor:proj-FE},
the proof of Proposition~\ref{thm:two} and Example~\ref{ex:proj-converse}).
In Subsection 4.2:
the ultrametric inequality, the
Lemmas~\ref{lem:Cn} and \ref{lem:ultra}, together with a first
proof of Corollary~\ref{cor:smooth-plane}.
In Section 5:
the reconstruction lemma (Lemma~\ref{lem:recons})
and the nice proof of the recovery problem for $r=3$
(Lemma~\ref{lem:one-variable}, Theorem~\ref{thm:xi-r3} and Remark~\ref{rem:r4}).
(c) Several examples (among them Example~\ref{ex:spacecusps}), the computer
verification of the examples of the paper, and the computer search reported in
Remark~\ref{rem:Cn-general}.
(d) Language editing, proposals and revisions for the abstract and the introduction.
All the content produced with the help of this tool has been checked by the author,
who takes full responsibility for the content of the publication.

\end{document}